\documentclass[12pt,letterpaper,titlepage]{amsart}
\usepackage{amsmath, amssymb, amsthm, amsfonts,amscd,xr}
\newcommand{\N}{\mathbb{N}}

\newcommand{\Z}{\mathbb{Z}}
\newcommand{\R}{\mathbb{R}}
\newcommand{\C}{\mathbb{C}}

\def\beq{\begin{equation}}
\def\eeq{\end{equation}}

\def\arr{\hbox to 20pt{\rightarrowfill}}
\def\mapright#1{\smash{\mathop{\arr}\limits^{#1}}}

\def\a{\alpha}
\def\b{\beta}

\def\d{\delta}
\def\e{\varepsilon}

\def\sl{\mathfrak {sl}}
\def\so{\mathfrak {so}}
\def\sp{\mathfrak {sp}}
\def\su{\mathfrak {su}}

\def\H{\mathcal H}
\def\O {\mathcal O}
\def\M{\mathcal M}

\newenvironment{res}
               {\begin{equation}
\begin{minipage}{0.85\textwidth}}
               { \end{minipage}\end{equation} }
\def\ber{\begin{res} }
\def\eer{\end{res}}

\numberwithin{equation}{section}
\newtheorem{thm}{Theorem}[section]
\newtheorem{theorem}[thm]{Theorem}

\newtheorem{lemma}[thm]{Lemma}
\newtheorem{lem}[thm]{Lemma}

\newtheorem{cor}[thm]{Corollary}
\newtheorem{ex}[thm]{Example}
\newtheorem{prop}[thm]{Proposition}

\newtheorem{dfn}[thm]{Definition}
\newtheorem{rem}[thm]{Remark}
\newtheorem{problem}[thm]{Problem}

\makeatletter
\def\section{\@startsection {section}{1}{\z@}{3.5ex plus 1ex minus
    .2ex}{2.3ex plus .2ex}{\large\bf}}
    \def\subsection{\@startsection{subsection}{2}{\z@}{3.25ex plus 1ex minus
 .2ex}{1.5ex plus .2ex}{\bf}}
\makeatother
\def\pf{{\em Proof}.\, }

\def\bysame{\leavevmode\hbox to3em{\hrulefill}\,}
\def\a{\alpha}
\def\b{\beta}
\def\l{\lambda}
\def\e{\epsilon}
\def\g{\gamma}

\def\d{\delta}

\def\af{\mathfrak{a}}
\def\gf{\mathfrak{g}}
\def\qf{\mathfrak{q}}
\def\h{\mathfrak{h}}

\def\kf{\mathfrak{k}}

\def\nf{\mathfrak{n}}

\def\pf{\mathfrak{p}}

\def\s{\sigma}
\def\tf{\mathfrak{t}}
\def\u{\mathfrak{u}}
\def\v{\mathfrak{v}}

\def\Ad{\operatorname{Ad}}

\def\Nc{\mathcal{N}}
\def\diag{\operatorname{diag}}

\def\DD{\mathbb{D}}

\def\O{\mathcal{O}}

\def\Cc{\mathcal{C}}

\def\M{\mathcal{M}}
\def\I{\mathcal{I}}
\def\Aut{\operatorname{Aut}}

\def\H{\mathcal{H}}

\def\J{\mathcal{J}}

\def\L{\mathcal{L}}
\def\NN{\mathcal{N}}
\def\E{\mathcal{E}}

\def\Qc{\mathcal{Q}}
\def\PP{\mathbb{P}}

\def\Sl{\operatorname{Sl}}
\def\Gl{\operatorname{Gl}}
\def\SO{\operatorname{SO}}

\def\W{\mathcal{W}}
\def\Ri{\mathcal{R}}

\renewcommand{\Re}{\mbox{\rm Re}\,}
\renewcommand{\Im}{\mbox{\rm Im}\,}
\makeindex
\begin{document}
\title[Analysis on the crown domain]
{Analysis on the crown domain}
\author{Bernhard Kr\"otz}
\address{Max-Planck-Institut f\"ur Mathematik, Vivatsgasse 7,
D-53111 Bonn, Germany\\ email: kroetz@mpim-bonn.mpg.de}
\author{Eric Opdam}
\address{Korteweg de Vries Institute for Mathematics\\
University of Amsterdam\\ Plantage Muidergracht 24\\ 1018TV
Amsterdam\\ The Netherlands\\ email: opdam@science.uva.nl}
\date{\today}
\thanks{During the preparation of this paper the second
named author was partially supported by a Pionier grant of
the Netherlands Organization for Scientific Research (NWO).
Part of this research was carried out in the fall of 2004, during
which period both authors enjoyed the hospitality of the Research Institute for the Mathematical Sciences in Kyoto,
Japan. It is our pleasure to thank the RIMS for its hospitality and for the stimulating environment it offers.}
\maketitle
\tableofcontents
\section{Introduction}

Our concern is with harmonic analysis on a Riemannian symmetric space
$$X_\R =G_\R/K_\R$$
of the noncompact type. Here $G_\R$ denotes a connected reductive algebraic group and $K_\R$ is a maximal compact
subgroup thereof.

\par Given a $K_\R$-spherical irreducible unitary representation $(\pi, \H_\pi)$ of $G_\R$ with $K_\R$-fixed ray
$\H_\pi^{K_\R-{\rm fix}}=\C v_K$, we obtain an $G_\R$-equivariant continuous map

$$i_\pi: X_\R \to \H_\pi, \ \ gK_\R\mapsto \pi(g) v_K .$$
We assume that $\pi\neq {\bf 1}$ is non-trivial and then $i_\pi$ is injective. The map $i_\pi$ is analytic, hence
admits holomorphic extension to a maximal $G_\R$-neighborhood $\Xi_\pi$ of $X_\R$ in $X_\C=G_\C/K_\C$. It is a
remarkable fact that $\Xi_\pi$ is independent of the choice of $\pi\neq {\bf 1}$ (\cite{KSI},\cite{KSII}, \cite{KOS})
and hence defines a natural domain $\Xi$ in $X_\C$, referred to as the {\it crown domain}. A result in this paper
determines the precise growth rate of $\|i_\pi\|$ when approaching the boundary of $\Xi$.

\par We have to clarify what we understand by the term "approaching the boundary".  The crown domain
admits a natural Shilov-type boundary \cite{GKI}, referred to  as the distinguished boundary $\partial_d\Xi$ of $\Xi$.
In a first step we give a simple description of $\partial_d\Xi$ in terms of the affine Weyl group, hereby extending and
unifying results from \cite{GKI}. At this point it is relevant that the $G_\R$-equivalence classes in $\partial_d\Xi$ are
described by a finite union of Weyl group orbits.

\par Given a distinguished boundary point $z\in \partial_d\Xi$ and $(z_n)\subset \Xi$ a sequence converging radially to
$z$ we are interested in the growth of $\|i_\pi (z_n)\|$ in terms of ${\rm dist}(z_n,z)$. We determine

\begin{itemize}
\item  For fixed $z_n$,  sufficiently close to $z$, optimal lower exponential bounds
for $\|i_\pi(z_n)\|$ in terms of the parameter of $\pi$ and ${\rm dist}(z_n,z)$;
\item For fixed $\pi$, the precise blow up rate of $\|i_\pi(z_n)\|$ for
$z_n\to z$ in terms of ${\rm dist}(z_n,z)$.
\end{itemize}

We use these results  to prove estimates for Maa\ss{} automorphic forms.
For example, a theorem of Langlands asserts
that cuspidal automorphic forms are of rapid decay
\cite{La}, \cite{HC}. An unpublished theorem of J. Bernstein goes beyond and asserts
exponential decay. In this paper we establish precise exponential decay rates.
The basic idea of proof goes back to J. Bernstein and our
contribution lies in a incorporation of geometric methods and
hard estimates.\footnote{JB explained to us
the case of $G=\Sl(2,\R)$, cf. the first half  of Subsection 9.2.}
In particular, we show that the crown domain admits a
natural parameterization by unipotent $G_\R$-orbits which makes Bernstein's
idea work out efficiently.

\par Finally we wish to point out that we make a detailed study
of proper actions of $G_\R$ on $X_\C$ in Section 4.
As a byproduct of these investigations we obtain a complex geometric
classification of the different series of representations of the group $G=\Sl(2,\R)$ (cf. Theorem \ref{th=sl2r}  below).

\medskip

{\it  Acknowledgment:} We express our gratitude to Joseph Bernstein for generously sharing his insights with us.
We thank the number theorists Erez Lapid and Andre Reznikov for pointing out interesting questions and their useful
hints to the literature. We are indebted to Philip Foth who did an elegant matrix computation for us. Finally we thank
Joachim Hilgert for asking a good question during one of our talks.
\section{The complex crown  and its distinguished boundary}

This section is divided into two parts. First we recall the
definition and basic properties of the complex crown $\Xi$
of a Riemannian symmetric space $X$ (see \cite{KSII} for a comprehensive
account). Second we shall
unify and extend results from \cite{GKI} on the
distinguished boundary of $\Xi$.

\subsection{The complex crown}

Let $G$ be a connected, real semisimple, noncompact Lie group.
Write $\gf$ for the Lie algebra of $G$ and
denote by $\gf_\C$ its complexification.
We fix a maximal
compact subalgebra $\kf\subset \gf$ and set $K=\exp(\kf)$.

\par Let us denote by $G_\C$ the universal
complexification of $G$ and by $\iota: G \to G_\C$
the homomorphism sitting over the injection $\gf\hookrightarrow \gf_\C$.
Write $K_\C$ for the analytic subgroup of $G_\C$
corresponding to $\kf_\C$.
\par Our concern
is with the Riemannian symmetric space
$X=G/K$.  The complex symmetric space
$X_\C=G_\C/K_\C$ naturally acts as a complexification of $X$
and the assignment
$gK\mapsto \iota(g)K_\C$ identifies
$X$ as a totally real submanifold of $X_\C$ in a $G$-equivariant way.
We denote the base point $eK_\C$ of $X_\C$ by $x_0$.

\begin{rem}\label{rem=split} Let $\gf=\gf_1+\ldots + \gf_l$ be the
factorization of $\gf$ in simple Lie algebras and let
$\kf=\kf_1+ \ldots+\kf_l$ be the associated splitting
for $\kf$. Denote by $G_j$, $K_j$ the analytic subgroups
of $G$ corresponding to $\gf_j$, $\kf_j$. Then with
$X_j=G_j/ K_j$ there is
the equivariant isomorphism
$$X\simeq X_1\times \ldots \times X_l\, .$$
In a similar manner (and obvious notation)
$$X_\C \simeq X_{1,\C}\times\ldots \times  X_{l,\C}\ .$$
\end{rem}

In the light of the discussion in Remark \ref{rem=split} it is
no loss of generality to assume henceforth that $\gf$ is
simple.

\par Let $\gf=\kf\oplus \pf$ be the Cartan decomposition associated to
the choice of $\kf$, and choose $\af$
a maximal abelian subspace in $\pf$.
The {\it complex crown} $\Xi$ of $X$
is by definition
\begin{equation}
\Xi=G\exp(i\pi \Omega/ 2). x_0\subset
X_\C,
\end{equation}
where $\Omega\subset\af$ is given  by
\begin{equation}
\Omega=\{Y\in \af\mid {\rm spec}({\rm ad} Y)\subset ]-1, 1
[\}.
\end{equation}
According to \cite{AG}, $\Xi$ is a $G$-invariant
open subdomain of $X_\C$ with the $G$-action proper.
Actually $\Xi$ is Stein (see \cite{KSII} and the references
therein).
Let us point out that $\Xi$ is independent of the choice
of the flat $\af$ and therefore naturally attached to $X$.

\par The set $\Omega$ can be described
in terms
of the restricted root system $\Sigma=\Sigma(\gf,\af)$ as follows:
\begin{equation}\label{eq:omega}
\Omega=\{Y\in\af\mid|\a(Y)|<1\,\forall\a\in\Sigma\}.
\end{equation}
In particular we see that $\overline \Omega$ is a compact
$W$-invariant polyhedron. Here $W$, as usual, denotes
the Weyl group of $\Sigma$.

\begin{rem} \label{rem=tan}(Realization in the tangent bundle)
Set $\Omega^K={\pi\over 2}\Ad (K)\Omega$. As $\Omega$ is
an open  $W$-invariant convex subset of $\af$, Kostant's linear
convexity theorem implies that $\Omega^K \subset \pf$ is an
open $K$-invariant convex subset of $\pf$.
Write  $TX=G\times_K \pf$ for the tangent bundle
of $X$. Notice that $G$ acts properly on $TX$ and that
$G\times_K\Omega^K$ is a contractible $G$-equivariant
subset of $TX$ (base and fiber are contractible).
In \cite{AG} it was shown that the map
\begin{equation}\label{eq=tan} G\times_K \Omega^K\to \Xi, \ \ [g,Y]\mapsto
g\exp(iY). x_0
\end{equation}
is a $G$-equivariant diffeomorphism. In particular,
$G$ acts properly on $\Xi$ and $\Xi$ is contractible.
\end{rem}

In the sequel we write $\tf=i\af$ and let $T=\exp \tf$ be the
corresponding torus in $G_\C$. Notice that
$T_\C = A_\C =AT$ with $A=\exp \af$. We will also use the notation
$T_\Omega=\exp(i\pi \Omega/2)$.

\begin{rem}\label{rem=bd} (The boundary of $\Xi$)
\item{(i)} (Semisimple boundary part) The topological boundary
$\partial\Xi$ is a complicated union of $G$-orbits.
This is because not all $G$-orbits in $\partial\Xi$
meet $T. x_0$. Those which do
make up the  semisimple (or elliptic) part of the   boundary
$\partial_s\Xi=
G\exp(i\pi \partial\Omega/ 2). x_0$
of $\Xi$ (see \cite{AG, KSII}).
Equivalently,  $\partial_s\Xi$ describes the closed
$G$-orbits in $\partial \Xi$. One knows
that each $G$-orbit in $\partial\Xi$ has a
a unique semisimple orbit in its closure \cite{FH}, but a
satisfactory general description of $\partial \Xi$
is still missing.

\item{(ii)} (Properness) The polyhedron $\Omega$ is maximal
with regard to proper $G$-action, i.e. there does not exists
a larger connected subset $\tilde \Omega\supset \Omega$
such that $G$ would act properly on $G\exp(i\pi\tilde\Omega/2). x_0$
(cf. \cite{AG}). We  mention that $G$-stabilizers of points
in $\exp(i\pi\partial \Omega/2). x_0$ are noncompact subgroups
\cite{AG}.
\item{(iii)} (Dependence on isogenies)
It follows from (\ref{eq=tan}) that $\Xi$ is homeomorphic to
$\pf\times \Omega^K$. It means in particular that $\Xi$ only
depends on the isogeny class of the connected group $G$.
However, the
situation becomes different once we start to consider the boundary
$\partial \Xi$ of the crown in $X_\C$. It turns out that $\partial
\Xi$ is sensitive with regard to the choice of the connected group
$G$; for instance $\partial \Xi$ is different for $\SO(n,\C)$ and
its simply connected cover $\SO(n,\C)^\sim $. We will comment more
on that when we will discuss collapsing of boundary orbits below
(Example \ref{collapse}).
\end{rem}

\subsection{Distinguished and minuscule boundary of the crown}
\label{sub:distbdy}
The {\it distinguished boundary} $\partial_d\Xi$ of the
crown, introduced in \cite{GKI},
is defined by
$$\partial_d\Xi=G\exp(i\pi\partial_e\Omega/2). x_0\subset \partial\Xi $$
where $\partial_e\Omega$ is the (finite set) of extreme points
of the compact polyhedron $\overline \Omega$. In view of this definition and
the results of the previous subsection we may and will assume
that $\mathfrak{g}$ is simple in this
subsection.
Let us recall that the distinguished boundary plays
the r\^ole of a noncompact
Shilov-type boundary  of $\Xi$; one has
the  following elementary result.

\begin{prop}\label{prop:sup}(\cite{GKI})
Let $f$ be a holomorphic function on $\Xi$ which extends
to a bounded continuous function on $\overline{\Xi}$.
Then
\begin{equation}
\operatorname{sup}_{x\in \Xi}(|f(x)|)=
\operatorname{sup}_{x\in \partial_d(\Xi)}(|f(x)|).
\end{equation}
\end{prop}

In \cite{GKI} a complete characterization of those crowns $\Xi$ was
given which admit symmetric
spaces as components in $\partial_d\Xi$. Cases relevant for
\cite{GKI} are those $\Sigma$ which are not of type
$E_8, G_2$ or $F_4$.

\par The objective of this section is to give
a uniform approach to $\partial_d\Xi$ in the general case.
Our first result is a  description of
$\partial_e\Omega$
in terms of structure theory which is stunningly simple
(cf. Theorem \ref{th=omega} below).
We will define the
{\it minuscule part} of the distinguished boundary and tie
it with the results of \cite{GKI}. After that
we classify the non-symmetric boundary components of $\partial_d\Xi$.
Finally we discuss collapsing of distinguished
boundary orbits.

\par Write $\Sigma^l=\{ \alpha: 2\alpha\not\in \Sigma\}$ for the
irreducible reduced subsystem of unmultipliable roots in
$\Sigma=\Sigma(\mathfrak{g},\mathfrak{a})$. It is clear
that $\Sigma^l$ completely describes $\Omega$, i.e.,
\begin{equation}
\Omega=\{ Y\in \af\mid |\alpha(Y)|<1, \forall \alpha\in \Sigma^l\}
\, .\end{equation}
Fix a basis for $\Sigma^l$, say
$\Pi=\{ \alpha_1,\ldots, \alpha_n\}$, and write $C\subset \af$
for the closure  of the associated Weyl chamber. Let $\beta$ be
the highest root corresponding to $\Pi$ and
$$\beta=k_1 \alpha_1+ \ldots + k_n \alpha_n$$
its expansion in the simple roots (hence $k_i\in \Z_{>0}$).
We record the obvious relation
\begin{equation}\label{eq=aff} \overline \Omega\cap C=\{ Y\in
C : \beta(Y)\leq 1\}\, .
\end{equation}
It means that $ \overline \Omega\cap C$ is a
fundamental domain for the affine Weyl group
$W^{\rm aff}=W\ltimes Q^\vee$ with $Q^\vee={\rm span}_\Z\Sigma^\vee$ the coroot
lattice in $\af$ (observe that
${\rm span}_\Z\Sigma^\vee={\rm span}_\Z(\Sigma^l)^\vee$).

\par Define $\omega_i\in \af$ by $\alpha_j(\omega_i)=\delta_{ij}$.
It is
straightforward from (\ref{eq=aff}) that
\begin {equation}\label{eq=inc}
\partial_e\Omega\cap C\subset\left \{ \omega_1/ k_1, \ldots,
\omega_n/ k_n\right \}\end{equation}
and so $\partial_e \Omega\subset W. \left \{ \omega_1/k_1, \ldots,
\omega_n/ k_n\right\}$ (cf.  \cite{GKI}, Lemma 3.17).

In general the inclusion in (\ref{eq=inc}) is proper and
we have to determine which $\omega_i/ k_i$ actually occur.
The key observation is contained
in Lemma \ref{lem=cone} below.
\par We need some terminology. Let $(V, (\cdot, \cdot))$
be an Euclidean space and $W\subset {\rm O}(V)$ be a Weyl group
of finite type associated to a root system with root basis $F$. We
shall assume that the action is effective, or equivalently
that $V^*={\rm span}_\R F$. For a subset $P\subset F$ let
$W_P<W$ be the corresponding parabolic subgroup. As before
$$C=\{ v\in V: \alpha(v)\geq 0 \forall \alpha\in F\}$$
denotes the closure of the  Weyl chamber.
A closed convex cone $\Gamma\subset V$ will be called non-degenerate if its
edge $E(\Gamma)=\Gamma\cap -\Gamma$ is equal to $\{0\}$. Clearly $\Gamma$
is non-degenerate iff there exists a linear functional
$\omega\in V^*$ such that $\omega|_{\Gamma\backslash\{0\}}>0$.

\begin{lemma}\label{lem=cone} Let $W$ be a Weyl group of finite type
acting effectively on an Euclidean space $V$. Let
$C$ be the closure of the corresponding Weyl chamber.
Then the following statements are equivalent:
\begin{enumerate}\item
$W$ is irreducible.
\item $W_P.C$ is  non-degenerate for all proper subsets
$P\subsetneq F$.
\item $W_P.C$ is  non-degenerate for a maximal proper
subset $P\subsetneq F$.
\end{enumerate}
\end{lemma}
\begin{proof} \par $(i)\Rightarrow (ii)$: If $P=\emptyset$, then
$W_P.C=C$ is non-degenerate. So let us henceforth assume that
$P\neq \emptyset$.
Denote by
$V_{\rm fix}=\{ v\in V\mid (\forall w\in W_P)\ w(v)=v \}$
the space of $W_P$-fixed points. Then
$$V=V_{\rm fix}\oplus V_{\rm eff}$$
with $V_{\rm eff} ={\rm span}_\R P =V_{\rm fix}^\bot$ the effective
part for the $W_P$-action.  We note that
$C\cap V_{\rm fix}\neq \{0\}$ and fix a non-zero element
$u$ in this intersection.
\par  Assume that $W$ is irreducible.
According to \cite{Hump}, Ch. IV, Exc. 8, one has
$(x,y)>0$ for all $x,y\in C\backslash\{0\}$.
In particular if $\omega\in V^*$ is defined by
$\omega(v):=(u,v)$ then $\omega|_{C\backslash\{0\}}>0$.
As $u$ is $W_P$-fixed,
it follows that $\omega|_{W_P. C\backslash\{0\}}>0$
and consequently $W_P. C$ is non-degenerate.
\par $(ii)\Rightarrow(iii)$ is clear, moving on to
$(iii)\Rightarrow(i)$: We argue by contradiction and
assume that $W$ is reducible. Then there exist splittings
$W=W_1\times W_2$, $V=V_1\oplus V_2$,
$F=F_1\amalg F_2$ with $W_1$ irreducible, $F_1\subset P$ and
$V_1, V_2\neq \{0\}$.
But then $V_1\subset W_P. C$ and
$W_P. C$ is degenerate.
\end{proof}
Let us now return to our initial setting with the
irreducible restricted root system $\Sigma=\Sigma(\gf,\af)$
(then the reduced root system $\Sigma^l$ of unmultipliable roots
is irreducible as well).
We write $D$ for the Dynkin diagram associated with the
bases $\Pi$ of $\Sigma^l$, and
$D^*=D(W^{\rm aff})$ for its affine extension. Let
$\Pi_0=\{\alpha_0,\alpha_1,\dots,\alpha_n\}$ denote
the underlying set of affine simple roots.
\begin{thm}\label{th=omega} Let $\gf$ be a simple Lie algebra
and $\Sigma=\Sigma(\gf,\af)$ the associated irreducible root system.
Then for all
$1\leq i\leq n$ the following statements are
equivalent:
\begin{enumerate} \item
$\omega_i/k_i\in \partial_e\Omega$.
\item $D^*-\{\alpha_i\}$ is connected.
\end{enumerate}
\end{thm}
\begin{proof} Fix $1\leq i\leq n$ and denote the stabilizer
of  $\omega_i/k_i$ in $W^{\rm aff}$ by $W^{(i)}$. Notice that
$W^{(i)}\simeq  W(D^*-\{\alpha_i\})$ is a Weyl group of finite type,
and that
$\Pi^{(i)}=\{\alpha_0,\dots,\alpha_{i-1},\alpha_{i+1},\dots,\alpha_n\}$
is a set of simple roots for its root system.
Let us denote by $C^{(i)}$ the associated closed Weyl chamber.

\par Let $U$ denote an open ball around $\omega_i/ k_i$
such that for $w\in W$ one has $w(U)\cap C\not=\emptyset$
iff $w$ fixes $\omega_i/k_i$. The isotropy group
of $\omega_i/k_i$ in $W$ is $W_P$ where $P=\Pi\cap \Pi^{(i)}$.
Observe that $P$ is a maximal proper subset both of
$\Pi$ and of $\Pi^{(i)}$.

Observe that $C \cap \overline \Omega$ is the fundamental
alcove of $W^{\rm aff}$ with respect to $\Pi_0$.
Hence $(C \cap \overline \Omega)\cap U=C^{(i)}\cap U$.
Moreover $\overline \Omega\cap U=W(C \cap \overline \Omega)\cap U
=W_P(C \cap \overline \Omega)\cap U=W_P(C^{(i)}\cap U)=
W^{(i)}_{P}(C^{(i)}\cap U)=W^{(i)}_{P}C^{(i)}\cap U$.
Hence $\omega_i/k_i$ is an extremal point of
$\overline\Omega$ iff the convex cone $W^{(i)}_{P}C^{(i)}$ (with
vertex $\omega_i/k_i$) is non-degenerate.
Apply Lemma \ref{lem=cone}.
\end{proof}

Let us call $\omega_i$ {\it minuscule} if $k_i=1$. Notice
that $\omega_i\in \partial\Omega$. Let us denote the union
of all $W$-orbits through minuscule $\omega_i$ by
$\partial_m \Omega$ and refer to it as the
{\it minuscule part} of $\partial\Omega$. Similarly we
define the {\it minuscule boundary} of
$\Xi$ by
\begin{equation}
\partial_m\Xi=G\exp(i\pi \partial_m \Omega/2). x_0\, .
\end{equation}

\begin{prop}\label{prop:incl} One has
\begin{equation}
\partial_m\Omega\subset \partial_e\Omega\end{equation}
and in particular $\partial_m\Xi\subset\partial_d\Xi$.
\end{prop}
\begin{proof} Let $A=(a_{ij})$ (with $i,j\in\{0,1,\dots,n\}$)
be the generalized Cartan matrix associated with the extended
Dynkin diagram $D^*$.
We consider $A$ as a matrix with respect to
the bases $\{\a_0,\a_1,\dots,\a_n\}$ of affine simple roots.
By elementary theory of generalized Cartan matrices
(cf. \cite[Theorem 4.8]{KacFirstEdition})
the one-dimensional kernel of $A$ is generated by a unique
positive, primitive element $\delta$ in the affine root
lattice, namely $\delta=\a_0+\beta$.
In other words, if we put $k_0=1$
then for each $j$:
$2k_j+\sum_{i\not=j} a_{ji}k_i=0$.
Hence if $\omega_j$ is minuscule (i.e. $k_j=1$) then
either $\a_j$ is an end point in $D^*$
(i.e. has only one neighbor in $D^*$) or else
$\a_j$ has precisely
two neighbors $\a_{i}$, $\a_{l}$
with $k_{i}=k_{l}=1$. But in this
last case $D^*$ must be (by an easy inductive argument)
a circular graph. We conclude in both cases that
$D^*-\{\a_j\}$ is connected as desired.
\end{proof}
For later reference and convenience to the reader
we list $\partial_e\Omega$ and $\partial_m\Omega$.
Theorem \ref{th=omega} and the tables of \cite {B}
yield:
\medskip

{\tiny{
\begin{center}
\begin{tabular}[c]{|c|c|c|}
\hline
\multicolumn{3}{|c|}{
\rule[-2mm]{0mm}{7mm}\textsf{Distinguished and minuscule boundary
of $\Omega$}} \\[2pt]
\hline\hline
$\Sigma$& $\partial_e\Omega\cap C$
& $\partial_m\Omega \cap C$ \\[2pt]\hline\hline
$A_n$  &$\omega_1, \ldots, \omega_n$ &
$\omega_1, \ldots, \omega_n$\\ [2pt]\hline
$B_n$ $(n\geq 3)$ &  $ \omega_1, \omega_n/2$ &
$ \omega_1$\\[2pt]\hline
$C_n, BC_n$  &  $\omega_n$&
$\omega_n$  \\[2pt]\hline
$D_n$ $(n\geq 4)$ & $ \omega_1,  \omega_{n-1}, \omega_n$
&$ \omega_1, \omega_{n-1}, \omega_n$\\ [2pt]\hline
$E_6$  & $\omega_1,\omega_6$
&  $ \omega_1,\omega_6$ \\ [2pt]\hline
$E_7$  &  $\omega_2/2, \omega_7$&
$\omega_7$\\ [2pt]\hline
$E_8$  &  $\omega_1/2, \omega_2/3$  & $\emptyset$ \\ [2pt]\hline
$F_4$  &  $\omega_4/2$  & $\emptyset$ \\ [2pt]\hline
$G_2$  &  $\omega_1/3$  & $\emptyset$ \\[2pt]\hline
\hline
\end{tabular}
\end{center}}}
\centerline{\tt Table 1}

\medskip

\begin{rem}(Correcting literature) The first named author would like to take
the opportunity to point out an error in \cite{GKI}
regarding $\partial_e\Omega$ for the $E_7$-case.
Due to a computational mistake
the $W$-orbit through $\omega_2/2$ was missed.
\end{rem}

For a point $\omega_j/k_j\in \partial\Omega$ set
$$z_j=\exp(i\pi\omega_j/2k_j). x_0\in \partial_d\Xi$$
and denote by $H_j$ the stabilizer of $G$ in $z_j$. We already remarked
earlier that $H_j$ is a noncompact subgroup.  Let us denote by
$\h_j$ its Lie algebra. Our next objective
is to classify the stabilizer algebras $\h_j$ for those
$z_j$ which appear in $\partial_d\Xi$.

\par Write $F=A_\C\cap K_\C =T \cap K$ and notice that
$F$ is a finite two group. We will often identify $A_\C. x_0$ with
$A_\C/F$ and remark that elements $z\in A_\C.x_0$
have well defined squares $z^2 \in A_\C$.
For each $1\leq j\leq n$ let us
define the centralizer subgroup
$$G_j:=Z(z_j^4)=\{ g\in G\mid z_j^4 gz_j^{-4}=g\}$$
and denote by $\gf_j$ its Lie algebra.

\begin{lem} \label{lem=m}Let $1\leq j\leq n$. Then the following assertions
hold:
\begin{enumerate}
\item $\gf_j=\gf$ if and only if $\omega_j$ is minuscule.
\item $\gf_j$ is a $3$-graded reductive Lie algebra
\begin{equation}\gf_j=\gf_{j,-} \oplus \gf_{j,0} \oplus \gf_{j,+}
\end{equation}
where $\gf_{j,\pm}=\{ Y\in \gf\mid
[\omega_j,Y]=\pm k_j Y\}$ and
$\gf_{j,0}=\{ Y\in \gf\mid
[\omega_j,Y]=0\}$.
\end{enumerate}
\end{lem}
\begin{proof}  Associated to $\omega_j$ is the standard grading
\begin{equation}\label{eq=grad} \gf=\sum_{l=-k_j}^{k_j} \gf_l\end{equation}
with $\gf_l=\{ Y\in \gf\mid [\omega_j,Y]=lY\}$. Notice that
$\Ad(z_j^4)$ acts on $\gf_l$ as the scalar $e^{2il\pi/k_j}$. The assertions of the lemma
follow with $\gf_{j,\pm}=\gf_{\pm k_j}$.
\end{proof}

Let us denote by $\theta$ the Cartan involution of $\gf=\kf\oplus \pf$.
Observe that $Y\in \gf$ belongs to $\h_j$ if and only if
$\Ad(z_j^{-1})(Y)\in \kf_\C$, in other words
\begin{equation} \h_j=\{ Y\in \gf\mid \Ad(z_j^2) (\theta Y)=Y\}\,,
\end{equation}
(cf.\ \cite{GKI}, Lemma 3.4).
We reveal the structure of $\h_j$.

\begin{lemma}\label{lem=h} Let $1\leq j\leq n$. Then $\h_j$ is $\theta$-stable
subalgebra of $\gf_j$. Moreover, its Cartan decomposition
is given by
$$\h_j=\gf_{0,j}^\theta \oplus (\gf_{j,-}\oplus \gf_{j,+})^{-\theta}\, .$$
\end{lemma}

\begin{proof} Recall the grading $\gf=\sum_{l=-k_j}^{k_j} \gf_l$
from (\ref{eq=grad}). Then for each $0\leq l\leq k_j$ the operator
$\Ad(z_j^2) \circ \theta $ leaves
$(\gf_l  \oplus \gf_{-l})_\C$ stable; explicitly
$$(Y_l,Y_{-l})\mapsto (e^{il\pi/k_j}\theta (Y_{-l}), e^{-il\pi/k_j}\theta(Y_l))\qquad (Y_l,Y_{-l})\in
(\gf_l  \oplus \gf_{-l})_\C\, .
$$
Hence $\left (\Ad(z_j^2)\circ \theta \right)(\gf_l \oplus \gf_{-l}) \cap \gf\neq\{0\}$ precisely for
$l=0,k_j$. The assertions of the lemma follow.
\end{proof}

\begin{cor} \label{cor=min} Let $1\leq j\leq n$. Then $\dim \h_j\leq \dim \kf$
with equality precisely if $\omega_j$ is minuscule.
\end{cor}

As a consequence of Lemma \ref{lem=m}
and Lemma \ref{lem=h} we can extend \cite{GKI}, Theorem 3.26 (2).

\begin{thm} For a boundary orbit $G.z_j\subset \partial\Xi$ the following statements
are equivalent.
\begin{enumerate}\item $\omega_j$ is minuscule.
\item $\dim \h_j=\dim \kf$.
\item $\Ad(z_j^{-1})(\h_j)_\C=\kf_\C$.
\item $\h_j$ is a symmetric subalgebra of $\gf$.
\item $G.z_j$ is a totally real submanifold of $X_\C$.
\item $G.z_j$ is a totally real submanifold of $X_\C$
of maximal dimension.
\end{enumerate}
\end{thm}
\begin{proof} (i)$\iff$(ii): Corollary \ref{cor=min}.
\par (ii)$\iff$(iii):  $\Ad(z_j^{-1}) \h_j\subset \kf_\C$
holds for all $1\leq j\leq n$ by the definition of $\h_j$.
\par (i)$\Leftarrow$(iv): If $\omega_j$ is minuscule, then
$\gf_j=\gf$ by Lemma \ref{lem=m}(i).  In particular
$\tau_j=\Ad(z_j^2)\circ \theta$ defines an involution and
$\h_j$ being the $\tau_j$-fixed point set is symmetric.
\par (iv)$\Rightarrow$(i): Notice that
$\gf_j$ is a reductive subalgebra properly containing $\h_j$.
Now if $\h_j$ is symmetric, then it is a
maximally reductive proper subalgebra of $\gf$. Thus
$\gf=\gf_j$, i.e. $\omega_j$ is minuscule.
\par (v)$\Rightarrow$(ii):
If $G.z_j$ is totally real, then $\dim_\R G.z_j\leq  \dim_\R X$.
The latter inequality rewrites as  $\dim \h_j \geq \dim \kf$.
Because of $\dim \h_j\leq \dim \kf$ in all cases,
it follows that $\dim \h_j=\dim \kf$.
\par (vi)$\Rightarrow$(v) is clear.
\par (ii) $\Rightarrow$ (vi): (ii) implies that
$\dim_\R G.z_j =\dim_\R X$. It remains to show that
$G.z_j$ is totally real. By $G$-homogeneity, it is sufficient
to show that $T_{z_j} (G.z_j)$ is totally real in $T_{z_j} (X_\C)$.
The assignment $Y\mapsto {d\over dt}\Big|_{t=0} \exp(tY).z_j$
identifies $\gf_\C/ \Ad(z_j) \kf_\C$ with  $T_{z_j} (X_\C)$.
Now observe  that  $\gf_\C/ \Ad(z_j) \kf_\C
=\gf_\C/\h_\C $ by the equivalence of (ii) and (iii).
Thus all we have to show is that $\gf +\h_\C/ \h_\C$
is totally real in $\gf_\C/ \h_\C$, which is apparent.
\end{proof}

\begin{rem} Suppose that $\omega_j$ is minuscule. Then
$\gf=\gf_{j, -}\oplus \gf_{j,0}\oplus \gf_{j,+}$ is a $3$-graduation and
$\tau_j=\Ad(z_j^2)\circ \theta$ is an involution
with fixed point algebra $\h_j$. In other words $(\gf, \h_j)$
is a noncompactly causal (NCC) symmetric pair. Moreover all
(NCC) symmetric pairs arise in this fashion.
For a proof of all this we refer to the paper \cite{Kan2} of
Professor Soji Kaneyuki (specifically Th. 3.1).

\par For the concrete classification of the $\h_j$ in the minuscule case
we allow ourselves to refer alternatively to \cite{GKI}, Th. 3.25.
\end{rem}

We wish to complete the classification of
$\partial_d\Xi$ by listing all the non-minuscule cases.
The most degenerate situation might deserve special attention.

\begin{ex} (The distinguished boundary of $G_2$)
Let us consider the case of $\gf=G_2$.
We use the terminology of \cite{B}. With $\Pi=\{\alpha_1, \alpha_2\}$
the positive roots list as
$$\alpha_1,\alpha_2, \alpha_1+\alpha_2, 2\alpha_1+\alpha_2, 3\alpha_1+\alpha_2, 3\alpha_1+
2\alpha_2\, .$$
We have $\partial_e \Omega= W. {\omega_1/3}$. Hence
\begin{eqnarray*}
\gf_{1,0}&=&\af\oplus  \gf^{\alpha_2} \oplus \gf^{-\alpha_2}\simeq \sl(2,\R)\times \R\\
\gf_{1,1} &=& \gf^{3\alpha_1+\alpha_2} \oplus \gf^{3\alpha_1+2\alpha_2}\simeq \R^2\\
\gf_{1,-1} &=& \gf^{-3\alpha_1-\alpha_2} \oplus \gf^{-3\alpha_1-2\alpha_2}\simeq \R^2
\end{eqnarray*}
and so $\gf_1\simeq \sl(3,\R)$.
Finally Lemma \ref{lem=h} implies
$\h_1\simeq \sl(2,\R)$.
\end{ex}

Let $z_j$ be a non-minuscule boundary points. A glance
at Table 1 above shows that $k_j=2$ except
for $G_2$ and one case in $E_8$. Thus
$\gf=\sum_{j=-k_j}^{k_j} \gf_j$ is a $5$-grading for most
of the cases.
\par Combining Lemma \ref{lem=h} and Lemma \ref{lem=m} with Kaneyuki's
classification of the even part of $5$-graded Lie algebras \cite{Kan}
we arrive at the following two lists.
For the exceptional cases we use \cite{Kan}, Table I, II.

\medskip
{\tiny{
\begin{center}
\begin{tabular}[c]{|c|c|c|c|c|c|}
\hline
\multicolumn{6}{|c|}{
\rule[-2mm]{0mm}{7mm}\textsf{Exceptional non-minuscule cases}} \\[2pt]
\hline\hline
$\gf$ & $\Sigma$ & $j$ &  $\gf_{0,j}$ & $\gf_j$ & $\h_j$ \\ [2pt]\hline\hline
$E_{6(2)}$  &$F_4$& 4& $\so(3,5)\times \R \times i\R$ & $\so(4,6)\times i\R$ &
$\so(1,3)\times \so(1,5)\times i\R$ \\ [2pt]\hline
$E_{6(-14)}$  & $BC_2$& 2&  $\so(1,7)\times \R \times i\R$ & $\so(2,8)\times
i\R$ &
$\so(1,7)\times \R \times i\R$ \\ [2pt]\hline
$E_{7(7)}$ & $E_7$ & 2 &$\sl(7,\R)\times \R $ & $\sl(8,\R)$ &
$\so(1,7)$ \\ [2pt]\hline
$E_{7(-5)}$& $F_4$ & 4& $\so(3,7)\times \su(2)\times \R $ & $\so(4,8)\times \su(2)$ &
$\so(1,3)\times \so(1,7)\times \su(2)$ \\ [2pt]\hline
$E_{8(8)}$ &$E_8$ & 1&  $\so(7,7)\times  \R $ & $\so(8,8)$ &
$\so(1,7)\times \so(1,7)$ \\ [2pt]\hline
$E_{8(8)}$ & $E_8$ & 2&  $\sl(8,\R)\times  \R $ & $\sl(9,\R)$ &
$\so(1,8)$ \\ [2pt]\hline
$E_{8(-24)}$ &$F_4$ & 4& $\so(3,11)\times  \R $ & $\so(4,12)$ &
$\so(1,3)\times \so(1,11)$ \\ [2pt]\hline
$F_{4(4)}$ &$F_4$ &4& $\so(3,4)\times  \R $ & $\so(4,5)$ &
$\so(1,3)\times \so(1,4)$ \\ [2pt]\hline
$F_{4(-20)}$ &$BC_1$ & 1&  $\so(7)\times  \R $ & $\so(1,8)$ &
$\so(1,7)$ \\ [2pt]\hline
$G_2$  & $G_2$ & 1& $\sl(2,\R)\times \R $ & $\sl(3,\R)$ &
$\so(1,2)$ \\ [2pt]\hline
$E_7^\C$  &$E_7$ & 2&  $\sl(7,\C)\times \C $ & $\sl(8,\C)$ &
$\su(1,7)$ \\ [2pt]\hline
$E_8^\C$  & $E_8$ & 1& $\so(14,\C)\times \C $ & $\so(16,\C)$ &
$\so(2,14)$ \\ [2pt]\hline
$E_8^\C$  &$E_8$ & 2& $\sl(8,\C)\times \C$ & $\sl(9,\C)$ &
$\su(1,8)$ \\ [2pt]\hline
$F_4^\C$  & $F_4$& 4& $\so(7,\C)\times  \C $ & $\so(9,\C)$ &
$\so(2,7)$ \\ [2pt]\hline
$G_2^\C$  & $G_2$ & 1& $\sl(2,\C)\times  \C $ & $\sl(3,\C)$ &
$\su(1,2)$ \\ [2pt]\hline
\hline
\end{tabular}
\end{center}}}
\centerline{\tt Table 2}
\medskip

For the classical cases we apply \cite{Kan}, Th. 3.2,
and note that the first two cases
below were already contained in \cite{GKI}, Th. 3.25.
\medskip
{\tiny{
\begin{center}
\begin{tabular}[c]{|c|c|c|c|c|}
\hline
\multicolumn{5}{|c|}{
\rule[-2mm]{0mm}{7mm}\textsf{Classical non-minuscule cases}} \\[2pt]
\hline\hline
$\gf$ & $\Sigma$ & $j$ &  $\gf_j$ & $\h_j$ \\ [2pt]\hline\hline
$\so(p,q)$  $(3\leq p<q)$ & $B_p$& $p$ & $\so(p,p)\times \so(q-p)$
&$\so(p,\C)\times \so(q-p)$\\ [2pt]\hline
$\so(2n+1, \C)$  $(n\geq 3)$ & $B_n$& $n$&  $\so(2n,\C)$ &
$\so^*(2n)$\\ [2pt]\hline
$\su(p,q)$  $(p<q)$ & $BC_p$& $p$& $\su(p,p)\times
\su(q-p)$
&$\sl(p,\C)\times \R \times \su(q-p)$\\ [2pt]\hline
$\sp(p,q)$  $(p<q)$ & $BC_p$& $p$& $\sp(p,p)\times \sp(q-p)$
&$\sl(p,{\Bbb H})\times \R \times \sp(q-p)$\\ [2pt]\hline
$\so^*(2n)$ $(n\geq 5,\  {\rm odd})$ & $BC_{[n/2]}$ & $[n/2]$ &
$\so^*(2n-2)$ &  $\sl((n-1)/2,{\Bbb H})\times \R$  \\ [2pt]\hline
\hline
\end{tabular}
\end{center}}}
\centerline{\tt Table 3}
\medskip

We conclude this section with a discussion of
collapsing of boundary orbits.  Let $z_j, z_l\in \partial_d\Xi$
with $j\neq l$.
If $G_\C$ is simply connected and $G\subset G_\C$,
then $G.z_j\neq G.z_l$, i.e $G.z_j\cap G.z_l=\emptyset$
(cf.\ \cite{GKI}, Th. 3.6). In the general case it might happen
that $G.z_j=G.z_k$ and we say that $\omega_j/k_j$ and $\omega_l/k_l$
{\it collapse} in $\partial_d\Xi$.
Collapsing appears when there exist outer automorphisms.
We refrain from complete results but would like
to mention some important examples.

\begin{ex}\label{collapse} (a) Let $G={\rm PSl}(n,{\Bbb K})$ for ${\Bbb K}=\R, \C,
{\Bbb H}$. Then $\omega_j$ and $\omega_l$ collapse
precisely for $j+l=n$.
\par (b) Let $G={\rm SO}(n,n)$ for $n\geq 4$. Then
$\omega_{n-1}$ and $\omega_n$ collapse.
\end{ex}

\section {New features of $G=\Sl(2,\R)$}\label{sec=nf}

This section is devoted to the crown domain
associated to the basic group $G=\Sl(2,\R)$.
It is divided into two parts. In the first half we
give a description of the full boundary $\partial\Xi$ as
a cone  bundle over the affine symmetric space $G/H=\Sl(2,\R)/
{\rm SO}(1,1)$. In the second part we give a novel
description of the crown as a union of unipotent $G$-orbits.
Later, via  appropriate $\Sl(2,\R)$-reduction, we
will use the material collected there
for our discussion of cusp forms and proper
action.

\subsection{Corner view}

We change perspective. Instead of regarding the crown from
the base point $x_0$  as a thickening of $X$, we may view
$\Xi$ from a corner point $z_j$ as a domain bordered by the
homogeneous space $G/H_j$. The advantage of this perspective
is that it leads to a simple characterization of the
full boundary $\partial \Xi$ of $\Xi$.

\par We will give a detailed discussion of
the boundary of the complex crown when $G=\Sl(2,\R)$.
As $\Xi$ is attached to $X$,  and so independent
of the specific global structure of $G$, we may
replace $\Sl(2,\R)$ by $G={\rm SO}_e(1,2)$.
We regard $K={\rm SO}(2,\R)$ as a maximal
compact subgroup of $G$ under the standard lower right
corner embedding.
\par Let us define a quadratic form $Q$ on $\C^3$ by

$$Q ({\bf z})=z_0^2-z_1^2 -z_2^2, \qquad
{\bf z}=(z_0, z_1, z_2)^T\in \C^3\, .$$
With $Q$ we declare real and complex hyperboloids
by
$$X=\{ {\bf x}=(x_0, x_1, x_2)^T\in \R^3\mid Q({\bf x})=1, x_0>0\}$$
and
$$X_\C=\{ {\bf z} =(z_0,z_1,z_2)^T\in \C^3\mid
Q({\bf z})=1\}\ .$$
We notice that  mapping
$$G_\C/K_\C\to X_\C, \ \ gK_\C\mapsto g.{\bf x}_0 \qquad ({\bf x}_0=(1,0,0))$$
is diffeomorphic and that $X$ is identified
with $G/K$.
\par
At this point it is useful to introduce coordinates
on $\gf=\so(1,2)$. We set

\begin{equation*}
{\bf e_1}=\begin{pmatrix} 0 & 0 & 1 \\ 0 & 0 & 0\\ 1 & 0 & 0\end{pmatrix},
\quad
{\bf e_2}=\begin{pmatrix} 0 & 1 & 0 \\ 1 & 0 & 0\\ 0 & 0 & 0\end{pmatrix},
\quad
{\bf e_3}=
\begin{pmatrix} 0 & 0 & 0 \\ 0 & 0 & 1\\ 0 & -1 & 0\end{pmatrix}\, .
\end{equation*}
We notice that $\kf=\R {\bf e_3}$, $\pf= \R{\bf e_1} \oplus  \R {\bf e_2}$
and make our choice of the flat piece $\af=\R{\bf e_1}$.
Then $\Omega=(-1, 1) {\bf e_1}$,
$\Xi=G\exp (i(-\pi/ 2,\pi/2){\bf e_1}).{\bf x_0}$
and we obtain Gindikin's favorite model of the
crown
$$\Xi=\{ {\bf z}={\bf x} +i{\bf y}\in X_\C\mid
x_0>0, Q({\bf x}) >0\}\, .$$
It follows that the boundary
of $\Xi$ is given by

\begin{equation} \label{be1}
\partial \Xi=\partial_s \Xi \amalg \partial_n \Xi\end{equation}
with semisimple part
\begin{equation} \label{be2}
\partial_s \Xi=\{ i{\bf y}\in i\R^3\mid
Q({\bf y})=  -1\}
\end{equation}
and nilpotent part
\begin{equation} \label{be3}
\partial_n \Xi=\{  {\bf z}={\bf x}+i{\bf y}
\in X_\C\mid x_0>0, Q{(\bf x})=0\}\, .
\end{equation}

Notice that
${\bf z_1}=\exp(i\pi/2 {\bf e_1}).{\bf x_0}=(0,0,i)^T$
and that the stabilizer of ${\bf z}_1$ in $G$ is
the symmetric subgroup $H={\rm SO}_e(1,1)$, sitting inside of $G$
as the upper left corner block.
Hence
\begin{equation}
\partial_s\Xi=\partial_d \Xi =G.{\bf z_1}\simeq G/H
\end{equation}
Write $\tau$ for the involution on $G$
with fixed point set $H$ and let
$\gf=\h \oplus \qf$ the corresponding $\tau$-eigenspace
decomposition. Clearly, $\h=\R {\bf e_2}$ and
$\qf=\af \oplus \kf=\R {\bf e_1} \oplus  \R {\bf e_3}$.
Notice that $\qf$ breaks as an $\h$-module into two pieces
$$\qf=\qf^+ \oplus  \qf^-$$
with
$$\qf^\pm=\{ Y\in \qf\mid [e_2, Y]=\pm Y\} =\R ({\bf e_1}\pm {\bf e_2})
\, .$$
Let us define the $H$-stable pair of half lines
$$\Cc=\R_{\geq 0}({\bf e_1}\oplus  {\bf e_3})\cup  \R_{\geq 0}
({\bf e_1}- {\bf e_3})$$
in $\qf=\qf^+ \oplus \qf^-$.  We remark
that $\Cc$ is the boundary of the $H$-invariant
open cone
$$\W=\Ad(H)(\R_{> 0}{\bf e_1})= \R_{> 0}({\bf e_1}+ {\bf e_3})\oplus  \R_{> 0}
({\bf e_1}- {\bf e_3})\, .$$
Recall that the tangent bundle $T(G/H)$ naturally identifies
with $G\times_H \qf$ and let us mention that $G\times_H \Cc$ is a $G$-invariant
subset thereof.

\begin{theorem}\label{th=bd} For $G={\rm SO}_e(1,2)$, the mapping
$$b: G\times_H \Cc\to \partial \Xi, \ \ [g,Y]\mapsto g\exp(-iY).{\bf z_1}$$
is a $G$-equivariant homeomorphism.
\end{theorem}

\begin{proof} It is  of course clear that
the map is equivariant and continuous. We move to surjectivity.
For $s\in \R$,

\begin{eqnarray*}\label{m-co}
\exp(is({\bf e_1} + {\bf e_3})) &=& \begin{pmatrix}
1 -s^2/2 & s^2/2 & is \\ -s^2/2 & 1+ s^2/2 & is \\
is & -is & 1 \end{pmatrix}, \\
\exp(is({\bf e_1} - {\bf e_3})) &=& \begin{pmatrix}
1 -s^2/2 & -s^2/2 & is \\ s^2/2 & 1+ s^2/2 & -is \\
is & is & 1 \end{pmatrix}\, . \end{eqnarray*}

Therefore,
\begin{equation} \label{be4} b([{\bf 1}, s({\bf e_1}\pm {\bf e_3})])=
\exp(-is({\bf e_1} \pm {\bf e_3})).{\bf z_1}=(s, \pm s, i)^T
\, .\end{equation}
{}From (\ref{be1}) - (\ref{be3}),
\begin{equation} \partial \Xi = G.\{ (s, \pm s, i)^T\mid
s\geq 0\} \end{equation}
and surjectivity is forced by (\ref{be4}) and $G$-equivariance.
\par Next, we prove that $b$ is one-to-one. By
$G$-equivariance, all we have to show is that
\begin{equation}\label{be5}
b([g, s({\bf e_1}\pm {\bf e_3})])=b([{\bf 1}, t({\bf e_1}\pm {\bf e_3})])
\end{equation}
for some $g\in G$ and $s,t\geq 0$, forces $g\in H$ and
$\Ad(g)(s({\bf e_1}\pm {\bf e_3}))=t({\bf e_1}\pm {\bf e_3})$.
We write (\ref{be5}) out and see
\begin{equation}\label{be6}
g.(s, \pm s, i)^T = (t, \pm t, i)^T\, .
\end{equation}
We take imaginary parts of this identity and deduce that
$g(0,0,i)^T=(0,0,i)^T={\bf z_1}$, i.e. $g\in H$.
With this information we go back in (\ref{be6}),
take the real part and get $g(s,\pm s,0)^T= (t,\pm t, 0)^T$.
We observe that the latter means
$\Ad(g)(s({\bf e_1}\pm {\bf e_3}))=t({\bf e_1}\pm {\bf e_3})$
and end the proof of injectivity.

\par Finally we mention that $b$ is an open mapping and
this finishes the proof.
\end{proof}

\begin{cor} For $G={\rm SO}_e(1,2)$ one has
$$\pi_1(\partial\Xi)=\pi_1(G/H)=\Z\, .$$
\end{cor}

\subsection{Unipotent parameterization}

We give now a novel description of the crown
as a union of unipotent $G$-orbits.

\par If not stated otherwise, $G=\Sl(2,\R)$.
The standard choices of coordinates are
$$\af=\R \begin{pmatrix} 1 & 0\\ 0 & -1\end{pmatrix}
\qquad\hbox{and} \qquad \nf=\R \begin{pmatrix} 0 & 1\\ 0 & 0\end{pmatrix}$$
and we observe that
$\Omega=(-1/2, 1/2)\begin{pmatrix} 1 & 0\\ 0 & -1\end{pmatrix}$.

The key observation is contained
in the following lemma.

\begin{lemma}\label{lem=or} Let $G=\Sl(2,\R)$. For all $0\leq |t|<{\pi/4}$, $t\in\R$, one has
the identity
\begin{equation}\label{eq=or} G\begin{pmatrix} 1 &  i \sin 2t \\ 0 & 1
\end{pmatrix}.x_0=G\begin{pmatrix}e^{it} & 0 \\ 0 & e^{-it}
\end{pmatrix}.x_0\, .
\end{equation}
\end{lemma}
\begin{proof} For the proof it is convenient to switch to the
hyperbolic model and replace $G$ by ${\rm SO}_e(1,2)$ (we recall that
$X=G/K$ and $\Xi$ are independent of the globalization $G$ of $\gf$; Remark \ref{rem=bd}(c) ).
As before, we choose $\af=\R {\bf e_1}$.
We come to our choice of $\nf$.
For $z\in \C$ let

$$n_z=\left( \begin{array}{ccc} 1+\frac{1}{2} z^2 & z &
-{1\over 2}   z^2  \\ z & 1 & -z \\
\frac{1}{2} z^2  & z & 1-\frac{1}{2}
z^2 \end{array} \right)
$$
and
$$N_\C=\{ n_z \mid z\in\C\}\, .$$

Further for $t\in \R $ with $|t|<{\pi\over 2}$ we set

$$a_t=\left( \begin{array}{ccc} \cos t & 0 & -i\sin t \\ 0 &
1 & 0 \\ -i \sin t & 0 & \cos t \end{array} \right) \in
\exp (i\Omega)\, .$$
\par The statement of the lemma translates into the
assertion
\begin{equation} Gn_{i\sin t }.{\bf x}_0=G a_t.{\bf x}_0\, .\end{equation}

Clearly, it suffices to prove that
$$a_t.{\bf x}_0=(\cos t, 0, -i\sin t)^T\in Gn_{i \sin t}.{\bf x}_0\, .$$

Now let $k\in   K$ and $b\in A$ be elements
which we write as

$$k= \left( \begin{array}{ccc} 1 & 0 & 0 \\ 0 &
\cos \theta & \sin\theta \\ 0 & -\sin \theta  & \cos \theta
\end{array} \right) \quad \hbox{and}\quad
b=\left( \begin{array}{ccc} \cosh r & 0 & \sinh r  \\ 0 &
1& 0 \\ \sinh r  & 0 & \cosh r \end{array} \right)$$
for real numbers $r,\theta$.
For $y\in \R$, a simple computation yields that

$$kbn_{iy}.{\bf x}_0= \left(\begin{array}{c}
\cosh r (1-{1\over 2} y^2) - {1\over 2}
y^2\sinh r \\
 iy\cos\theta + \sin\theta( \sinh r (1-{1\over 2} y^2) -
{1\over 2}y^2
\cosh r)\\ -iy\sin\theta +
\cos\theta( \sinh r (1-{1\over 2}y^2) -{1\over 2}y^2
\cosh r)
\end{array} \right)\, . $$
Now we make the choice of $\theta={\pi \over 2}$ which gives us
that

$$kb n_{iy}.{\bf x}_0= \left(\begin{array}{c}
\cosh r (1-{1\over 2} y^2) - {1\over 2}
y^2\sinh r \\
 \sinh r (1-{1\over 2}y^2) -{1\over 2}y^2
\cosh r\\ -iy
\end{array} \right)\, . $$
As $y=\sin t $ we only have to verify that we can choose $r$ such that
 $\sinh r (1-{1\over 2}y^2) -{1\over 2}y^2
\cosh r=0$. But this is equivalent to
$$\tanh r =  {{1\over 2} y^2 \over 1 -{1\over 2}y^2}\, .$$
In view of $-1< y=\sin t <1 $, the right hand side is  smaller than one
and we can solve for $r$.
\end{proof}

Let us define a domain in $\nf$
$$\Lambda=\left\{\begin{pmatrix} 0 & x \\ 0 & 0
\end{pmatrix}\in \nf\mid \, x \in \R, |x|<1 \right\}\, .$$
A remarkable consequence of the preceding Lemma is the following
result which we will establish in full generality later on.

\begin{thm}\label{th=sl2} For $G=\Sl(2,\R)$ one has
$$\Xi=G\exp(i\Lambda).x_0\, .$$
\end{thm}

\begin{rem}\label{rem=fb} (a) (Relation to $KNK$) As observed by Kostant, for any semisimple Lie group
$G$ one has $G=KNK$. As one referee pointed out,
a  more careful study of the $KNK$-decomposition of $G$ was undertaken
by H. Lee Michelson in \cite{Mi}. In particular, Prop. 3.1 in \cite{Mi}
applied to $G=\Sl(2,\R)$ states that
\begin{equation}\label{mm}  K \begin{pmatrix} e^t & 0 \\ 0 & e^{-t}\end{pmatrix}K =
K \begin{pmatrix} 1  & 2 \sinh t  \\ 0 & 1 \end{pmatrix}K
\end{equation}
for all $t\in \R$. It is tempting to believe that Lemma \ref{lem=or}
would follow from some sort of analytic continuation of
(\ref{mm}). However, we observe a subtle difference in this matter:
the location of the even prime. This is surprising and we thank this  referee
of having raised the question.
\par\noindent (b) It is not a priori
clear that $G\exp(i\Lambda).x_0$ is open in $X_\C$.
This is because of the fact that the natural
map
$$G\times \nf \to X_\C, \ \ (g, Y)\mapsto g\exp(iY).x_0$$
has singular differential at $(g, 0)$, $g\in G$.
\par\noindent (c) Lemma \ref{lem=or} allows us to give
a characterization of $\Xi$ as a fiber bundle
related to the nilcone. Write $\Nc\subset \gf$ for the
cone of nilpotent elements  in $\gf$ and note that
$\Nc=\Ad(K)\nf$. Define a subset of $\Lambda$ by
$$\Lambda^+=\left\{\begin{pmatrix} 0 & x \\ 0 & 0 \end{pmatrix}\mid
x \in \nf, \, 0\leq x <1 \right\}\, .$$
and put $\Nc^+=\Ad(K)\Lambda^+$. Then
it follows from Lemma \ref{lem=or} that the mapping
$$G\times_K \Nc^+\to \Xi, \ \ [g,Y]\mapsto g\exp(iY).x_0$$
is a homeomorphism.
\end{rem}

While it is not possible to enlarge $\Xi$ to a larger
domain in hyperbolic directions, i.e. beyond $\Omega$,
the situation is quite different for unipotent elements.

\begin{lemma}\label{lem=np} The differential of the mapping
$$G\times \nf \to X_\C, \ \ (g, Y)\mapsto g\exp(iY).x_0$$
is invertible at all points $(g,Y)\in G\times \nf$
with $Y\neq 0$.
\end{lemma}
\begin{proof} By $G$-equivariance of the map, it will
be the sufficient to show that the map is submersive
at all points $({\bf 1},Y)$ with $Y\neq 0$.
This assertion in turn translates into the identity
$$e^{-i {\rm ad} Y}\gf +i \nf + \kf_\C =\gf_\C$$
which is satisfied whenever $Y\neq 0$.
\end{proof}

For $a<b$ we define an open subset of $\nf$ by
$$\Lambda_{a,b}=\left\{\begin{pmatrix} 0 & x \\ 0 & 0
\end{pmatrix}\in \nf\mid x \in \R, a < x<b \right\}\, $$
and declare $G$-invariant connected subsets of $X_\C$
by
$$\Xi_{a,b}=G \exp (i\Lambda_{a,b}).x_0\, .$$
Of further interest for us is the limiting object
for $a\to -\infty , b \to \infty$,

$$\Xi_N= G\exp(i\nf).x_0=GN_\C.x_0\, .$$

\begin{lemma} For all $a<b$, the sets $\Xi_{a,b}$ are open.
In particular $\Xi_N$ is open.
\end{lemma}

\begin{proof} If $0\not \in (a,b)$, then the
assertion follows from Lemma \ref{lem=np}.
Thus we may assume that $0\in (a,b)$.
Suppose first that $|a|, |b|\leq 1$. Then by
Lemma \ref{lem=or} there exists a symmetric, i.e. $W=\Z_2$-invariant,
interval $\Omega_{a,b}\subset \Omega$ such that
$\Xi_{a,b}=G\exp(i\Omega_{a,b}).x_0$. The latter set is open
by Remark \ref{rem=tan}. Finally assume that $b>1$ or $a<-1$.
Then
$$\Xi_{a,b} =\Xi_{\max\{a, -1\},  \min\{1,b\}}
\cup  G\exp( i \Lambda_{a,b} \setminus\{0\}).x_0$$
is the union of two open sets (use Lemma \ref{lem=np} for the second
term) and we conclude the proof of the lemma.
\end{proof}

We exhibit the structure of the domain $\Xi_N$.
For that it is useful to move to the hyperboloid picture
with $G={\rm SO}_e(1,2)$. Define the
horocycle space of $X$ by

$${\rm Hor} (X)=\{ \xi\in \R^3\mid Q (\xi)=0,  \xi_0>0\}$$
and notice that the map
$$G/N\to {\rm Hor}(X), \ \ gN\mapsto g.\xi_0$$
with $\xi_0=(1, 0, 1)^T$ is a diffeomorphism.
Let us denote by
$${\bf z}\cdot {\bf w}=z_0 w_0-  z_1 w_1 - z_2 w_2 $$
the complex bilinear form obtained from polarizing
$Q$.

\begin{prop}\label{prop=xic} For the domain $\Xi_N$ the following
assertions hold:
\begin{enumerate}
\item $\Xi_N=\{ {\bf z}\in X_\C\mid {\bf z}\cdot \xi =1\  \hbox{for some}
\ \xi \in {\rm Hor}(X)\}$.
\item $\Xi_N= X_\C - \partial_d\Xi - \{{\bf z}\in X_\C \mid Q({\bf x})>0,
x_0<0\} $.
\item For all $z\in \Xi_N$, the $G$-stabilizer
$$G_z=\{ g\in G\mid g.z=z\}$$
is a compact subgroup of $G$.
\end{enumerate}
\end{prop}
\begin{proof} (i) We only have to notice that
$$N_\C.x_0=\{{\bf z}\in X_\C\mid {\bf z}\cdot \xi_0=1\}\, .$$
\par \noindent(ii) We use the characterization
of $\Xi_N$ from (i). We have to show that

\begin{equation}\label{eq=xin}
\Xi_N=\{{\bf  z}={\bf x}+i {\bf y}\in X_\C\mid {\bf x}\neq 0\}\, .
\end{equation}
For elements $z\in X_\C$ we will distinguish
three cases: $Q ({\bf x})>0$, $Q ({\bf x})<0$
and  $Q ({\bf x})=0$.
Before we do our case by case analysis let us
mention the fact that  elements ${\bf z}={\bf x}+i{\bf y} \in \C^3$
belong to $X_\C$ precisely when
$$Q({\bf x})- Q({\bf y})=1 \qquad\hbox{and}
\qquad {\bf x}\cdot {\bf y}=0\, .$$

\par\noindent Case 1: $Q ({\bf x})>0$ and $x_0<0$. We claim that
the $G$-orbit through ${\bf z}$ has a representative of the
type ${\bf z}=(x_0, iy_1, 0)^T$. In fact, as $Q({\bf x})>0$,
the $G$-orbit through ${\bf x}$ has a representative
$(x_0, 0, 0)^T$ with $x_0<0$.
{}From $ {\bf x}\cdot {\bf y}=0$ we then conclude that
${\bf y}=(0, y_1, y_2)^T$. Further we may alter
${\bf y}$ by the stabilizer $K={\rm SO}(2, \R)$ of ${\bf x}$.
Thus we may assume that ${\bf y}=(0, y_1, 0)^T$.
But then ${\bf z}\cdot \xi=1$ for $\xi=(1/x_0, 0, 1/x_0)^T\in {\rm Hor}(X)$.
In particular
$$\{z\in X_\C\mid Q({\bf x})>0\}\subset \Xi_N\, .$$

\par\noindent Case 2: $Q ({\bf x})<0$. We claim that
the $G$-orbit through ${\bf z}$ has a representative of the
type ${\bf z}=(0, iy_1, x_2)^T$. Indeed, as $Q({\bf x})<0$,
we may assume that ${\bf x}=(0, 0, x_2)^T$ with $x_2>0$.
Orthogonality $ {\bf x}\cdot {\bf y}=0$ then implies that
${\bf y}=(y_0, y_1, 0)^T$. Notice that
$$Q({\bf y})=y_0^2 -y_1^2 = -1 - x_2^2 <0\, .$$
It is allowed to change
${\bf y}$ by displacements of $H={\rm SO}(1,1)$, the stabilizer
of ${\bf x}$. As $H$ acts transitively on all connected
component of the level sets of  $y_0^2 -y_1^2$, it is
no loss of generality to assume that
${\bf y}=(0,y_1, 0)^T$. But then
${\bf z}\cdot \xi =1$ for $\xi=(1/x_2, 0, -1/x_2)^T\in {\rm Hor}(X)$
and we conclude that
$$\{z\in X_\C\mid Q({\bf x})<0\}\subset \Xi_N\, .$$
\par\noindent Case 3:   $Q ({\bf x})=0$ and
${\bf x}\neq 0$.  We assert that the $G$-orbit
through ${\bf z}$ has a representative of the type
${\bf z}=(1+iy_0,1+iy_0,\pm 1)^T$. Namely, as
$Q({\bf x})=0$ and ${\bf x}\neq 0$, the $G$-orbit
through ${\bf x}$ contains the element ${\bf x}=(1,1,0)^T$.
Then ${\bf x}\cdot {\bf y}=0$ and $Q({\bf y})=-1$ force
$y=(y_0, y_0, \pm 1)$. Hence we can choose ${\bf z}$ of the
asserted form. But then $\xi=( {y_0^2+1\over 2} , {y_0^2- 1\over 2}, \mp y_0)^T
\in {\rm Hor }(X)$ with ${\bf z}\cdot \xi=1$.

\par Finally, we observe that elements of the type ${\bf z}= i{\bf y}$
cannot belong to $\Xi_N$ by (i).

\par\noindent (iii) Notice that $g.{\bf z}={\bf z}$ means that
$g.{\bf x}={\bf x}$ and $g.{\bf y}={\bf y}$.
We analyze the three cases in (ii).
If $Q({\bf x})>0$, then $G.{\bf x}\simeq X$ and
$G_{\bf x}$ is compact.
If  $Q({\bf x})<0$, then we may assume that
${\bf z}=(0, iy_1, x_2)^T$. Hence $g.{\bf x}={\bf x}$ forces
$g\in H$ and then $g.{\bf y}={\bf y}$ yields
$g={\bf 1}$.
Finally if $Q({\bf x})=0$ and ${\bf x}\neq 0$, then our choice
of ${\bf z}$ can be ${\bf z}=(1+iy_0,1+iy_0,\pm 1)^T$.
Then $g.{\bf x}={\bf x}$ implies that $g$ is unipotent while
$g.{\bf y}={\bf y}$ forces $g$ to be hyperbolic. Hence
$g={\bf 1}$ in this case also.
\end{proof}

The statement in Proposition \ref{prop=xic} (iii) suggest
that the $G$-action on $\Xi_{a,b}$ should be proper. However, this
is not always the case as our next result shows.
For  that

\begin{prop} \label{prop=ab}Suppose that $(a,b)\cap (-1,1)\neq \emptyset$.
Then the  following assertions hold:
\begin{enumerate}
\item If $\max\{ |a|, |b|\}\leq 1$, then the $G$-action on $\Xi_{a,b}$ is proper.
\item If $\min\{ |a|, |b|\} >1$, then the $G$-action on
$\Xi_{a,b}$ is not proper.
\end{enumerate}
\end{prop}

\begin{proof} If $|a|, |b|\leq 1$, then $\Xi_{a,b}\subset \Xi$ and as the
$G$-action is proper on $\Xi$, the same holds for $\Xi_{a,b}$.
We move to (ii). Assume now that $|a|>1$ and  $|b|>1$. Then
$\Xi_{a,b}$ contains both elements
${\bf w}_+=n_i.{\bf x}_0$ and ${\bf w}_-=n_{-i}.{\bf x}_0$.
We note that
$${\bf w}_+:=\begin{pmatrix} 1/2 \\ i\\ -1/2\end{pmatrix}=
n_i.{\bf x}_0=
\left( \begin{array}{ccc} 1/2  & i & 1/2     \\ i & 1 & -i \\
-1/2   & i & 1/2
\end{array} \right)\cdot  \begin{pmatrix} 1 \\ 0\\ 0\end{pmatrix}\,.
$$
For $n\in \N$ we define elements ${\bf z}_n\in X_\C$ by
$${\bf z}_n =(1/2, 0, -1/2 + e^{-n})^T + i( 0, \sqrt{{3/4} + (e^{-n} -1/2)^2}, 0)
^T\, .$$
Notice that $\lim_{n\to \infty} {\bf z}_n={\bf w}_+$.
Hence there exists an $n_0\in \N$ such that ${\bf z}_n\in \Xi_{a,b}$ for
all $n\geq n_0$.
Now set

$$b_n=\left( \begin{array}{ccc} \cosh n & 0 & \sinh n  \\ 0 &
1& 0 \\ \sinh n  & 0 & \cosh n \end{array} \right)\in A\, .$$
Note that eigenvectors of $b_n$ are
$${\bf f}_1=\begin{pmatrix}  1 \\ 0 \\1\end{pmatrix},
\quad {\bf f}_2=\begin{pmatrix}  0 \\ 1 \\0\end{pmatrix},
\quad {\bf f}_3=\begin{pmatrix}  1 \\ 0 \\-1\end{pmatrix}$$
with eigenvalues $e^n$, $1$ and $e^{-n}$ respectively.
Thus
\begin{align*} b_n.{\bf z}_n = & b_n.{\bf x_n} + i {\bf y}_n \\
= & e^n\cdot {x_{0,n} +x_{2,n}\over 2}\cdot {\bf f}_1 +
e^{-n} \cdot {x_{0,n} -x_{2,n}\over 2}\cdot {\bf f}_2 + \\
& + i \sqrt{{3/4}
 + (e^{-n} -1/2)^2}\cdot {\bf f}_3 \\
=&  1/2 \cdot{\bf f}_1 + e^{-n} (1-e^{-n}/2)\cdot {\bf f}_2 +
i \sqrt{{3/4} + (e^{-n} -1/2)^2}\cdot {\bf f}_3
\end{align*}
and thus $\lim_{n\to \infty} b_n.{\bf z}_n ={\bf w}_-\in \Xi_{a,b}$.
Hence $(b_n.{\bf z}_n)_{n\geq n_0}$ stays in a compact subset
of $\Xi_{a,b}$ but with $(b_n)_{n\geq n_0}$ an unbounded
sequence. Thus the action of $G$ on $\Xi_{a,b}$ is
not proper.
\end{proof}

We conclude this section with a final result
for proper $G$-action.

\begin{prop} \label{prop=sl2} Let $G=\Sl(2,\R)$ and let $D\subset X_\C$
be a $G$-invariant domain with $X\subset D$. If the action of $G$ on $D$ is proper,
then:
\begin{enumerate}
\item $\partial_s\Xi\cap D=\emptyset$.
\item $\partial_n\Xi\not\subseteq D$.
\end{enumerate}
In particular, if $\partial_n\Xi\cap D =\emptyset$, then
$D\subseteq \Xi$.
\end{prop}

\begin{proof} Let $X\subset D\subset X_\C$ be an open $G$-invariant
domain with proper $G$-action. Suppose that $D\cap \partial \Xi\neq \emptyset$ and
let $z$ be a point thereof. Then $z\not \in \partial_d\Xi=G/H$ as $H$ is noncompact
and the $G$-action on $D$ is proper. Hence $z\in \partial_n\Xi$. It follows from
Theorem \ref{th=bd} that
$$\partial_n\Xi= Gn_i.{\bf x}_0 \amalg Gn_{-i}.{\bf x}_0\, .$$
Thus $D\cap \Xi_{a,b}\neq \emptyset $ for some $a,b$ with  $\max\{|a|, |b|\}>1$
-- the assertion now follows from the previous proposition.
\end{proof}

\begin{rem}\label{remidemmi} There exist larger $G$-domains $D\supsetneq\Xi$ with
the $G$-action proper. We provide the recipe for their
construction in case of $G=\Sl(2,\R)$. Recall
that $X$ identifies with the upper halfplane
and henceforth we view $X$ in the projective space $\PP^1(\C)$. Notice that
$G_\C$ acts on $\PP^1(\C)$ by fractional linear transformation.
Denote by $\overline X$ the lower half plane and notice that
$\Xi$ is $G$-isomorphic to $X\times \overline X$. In this
realization $X$ sits in $\Xi=X\times \overline X$ via
$z\mapsto (z,\overline z) $.
We view $\Xi\in \PP^1(\C)\times \PP^1(\C)$ and note that
$$X_\C=\{ (z,w)\in \PP^1(\C)\times \PP^1(\C)\mid z\neq w\}\, .$$
Furthermore
$$\partial_s\Xi=\{(x,y)\in \PP^1(\R)\times \PP^1(\R)\mid
x\neq y\}\, $$
and
$$\partial_n\Xi= X\times \PP^1(\R)\amalg \PP^1(\R)\times \overline X\, .$$
In particular we see that
$$D=(X\times \PP^1(\C))\cap X_\C$$
provides a $G$-domain in $X_\C$ such that
\begin{itemize}
\item $\partial_s\Xi\cap D=\emptyset$,
\item $\partial_n\Xi\not\subset D$,
\item $G$ acts properly on $D$.
\end{itemize}
With this picture of $\Xi$ one can easily sharpen
Proposition \ref{prop=ab} to: $G$ acts properly
on $\Xi_{a,b}$ if and only if $\min\{|a|,|b|\}\leq 1$.
\end{rem}

\section{Properness and maximality of holomorphic extension}

The first part of this section is valid for general $G$; the subsection after
for $G=\Sl(2,\R)$ only.

\par As we mentioned earlier in Remark \ref{rem=bd}(ii), it was proved in \cite{AG}, that $\Omega$ is maximal with respect
to proper $G$-action. We will refine this result in Theorem \ref{th=p} below. This new geometric fact translates into a
maximality assertion for holomorphic extension of representations.

\par \begin{thm}\label{th=p} Let $X\subset D$ be a $G$-domain in $X_\C$
with the $G$-action proper. Then the following assertions hold:
\begin{enumerate}
\item $\partial_s\Xi\cap D=\emptyset$.
\item $\partial_n\Xi\not\subset D$. In particular, if
$\partial_n\Xi \cap D=\emptyset$, then $D\subset \Xi$.
\end{enumerate}
\end{thm}

\begin{proof} (i) It was shown in \cite{AG} that
$G$-stabilizers on $\partial_s\Xi$ are noncompact.
Hence the assertion.
\par\noindent (ii) Suppose that $\partial_n\Xi\subset D$
and let $z$ be a point of $\partial_n\Xi$.
As $D$ is open we may assume that $z$
is generic in the sense of \cite{FH}, Section 4.2. It follows from
\cite{FH}, Th. 4.3.5., that there is a subgroup $G_0\subset G$ which is
locally isomorphic to $\Sl(2,\R)$ such that
the crown $\Xi_0$ associated to $G_0$ embeds $G_0$-equivariantly
into $\Xi$ with $z\in \partial_n \Xi_0$ in addition.
As $\partial_n\Xi_0\subset \partial_n\Xi$
we obtain a contradiction to
Proposition \ref{prop=sl2}.
\end{proof}

We turn to applications in representation theory.
For that it is convenient to look at the preimage
$$\tilde \Xi = G\exp(i\pi/2\Omega)K_\C$$
of $\Xi$ in $G_\C$.

We let $(\pi, \H)$ be a unitary
irreducible representation of $G$ and write
$\H_K$ for the associated Harish-Chandra modul
of $K$-finite vectors.
Then, for $v\in \H_K$,  it was shown in \cite{KSI} that
the orbit map
$$F_v: G\to \H, \ \ g\mapsto \pi(g) v$$
extends to a $G$-equivariant holomorphic map
$\tilde \Xi\to \H$, also denoted by $F_v$ in the sequel.
We wish to show that $\tilde \Xi$ is maximal
and want to relate this to the properness of the
action of $G$ on $\Xi$. The link is
established through the following fact.

\begin{lemma} \label{lem=hm}Let $(\pi, \H)$ be a unitary representation
of a reductive group $G$ which does not contain
the trivial representation. Then $G$ acts properly
on $\H -\{0\}$.
\end{lemma}

\begin{proof} Let $C\subset \H -\{0\}$ be a compact subset
and $C_G=\{ g\in G\mid \pi(g)C\cap C\neq \emptyset\}$.
Suppose that $C_G$ is not compact. Then there exists
a sequence $(g_n)_{n\in \N}$ in $C_G$ and a
sequence $(v_n)_{n\in \N}$ in $C$ such that
$\pi(g_n)v_n \in C$ and $\lim_{n\to \infty} g_n =\infty$.
As $C$ is compact we may assume that $\lim_{n\to\infty} v_n =v$ and
$\lim_{n\to\infty} \pi(g_n)v_n =w$ with $v,w\in C$.
We claim that
\begin{equation} \label{hm}
\lim_{n\to\infty} \langle \pi(g_n)v, w\rangle \neq 0\, . \end{equation}

In fact $\|\pi(g_n) v_n -\pi(g_n)v\|=\| v_n -v\|\to 0$
and thus $\pi(g_n)v\to w$ as well. As $w\in C$, it follows
that $w\neq 0$ and our claim is established.
\par Finally we observe that (\ref{hm}) contradicts
the Riemann-Lebesgue lemma for representations which asserts
that the matrix coefficient vanishes at infinity.
\end{proof}

{}From Lemma \ref{lem=hm} we deduce the following
result.

\begin{thm}\label{t1} Let $(\pi, \H)$ be an irreducible unitary
representation of $G$ which is not trivial. Let $v\in \H_K$,
$v\neq 0$,
be a $K$-finite vector. Let $\tilde D$ be a maximal
$G\times K_\C$-invariant domain in $G_\C$ with respect to
the property that the
orbit map $F_v: G\to \H, \ \ g\mapsto \pi(g)v$
extends to a $G$-equivariant holomorphic map
$\tilde \Xi\to \H$. Then $G$ acts properly on
$\tilde D/ K_\C \subset X_\C$.
\end{thm}

\begin{proof} We argue by contradiction and assume that
$G$ does not act properly on $D=\tilde D/K_\C$.
We obtain sequences $(z_n')_{n\in \N}\subset D$ and
$(g_n)_{n\in \N} \subset G$ such that
$\lim_{n\to\infty} z_n' =z'\in D$, $\lim_{n\to\infty} g_n z_n' =w'\in D$ and
$\lim_{n\to\infty} g_n=\infty$. We select preimages $z_n$, $z$ and $w$ of
$z_n'$, $z'$ and $w'$ in $\tilde D$. We may assume
that $\lim_{n\to\infty} z_n=z$ and find a sequence $(k_n)_{n\in \N}$ in $K_\C$
such that $\lim_{n\to\infty} g_n z_n k_n = w$.
\par Before we continue we claim that
\begin{equation}\label{eq=nz} (\forall z\in \tilde D)\qquad \pi(z)v\neq 0
\end{equation}
In fact assume $\pi(z)v=0$ for some $z\in \tilde D$. Then
$\pi(g)\pi(z)v=0$ for all $g\in G$. In particular
the map $G\to \H, \ \ g\mapsto \pi(g)v$ is constantly
zero. However this map extends to a holomorphic map
to a $G$-invariant neighborhood in $G_\C$. By the
identity theorem for holomorphic functions this
map has to be zero as well. We obtain a contradiction
to $v\neq 0$ and our claim is established.

\par Write $V= {\rm span} \{\pi(K)v\}$ for the finite dimensional
space spanned by the $K$-translates of $v$.
In our next step we claim that

\begin{equation}\label{eq=bou} (\exists c_1, c_2>0) \qquad  c_1 < \|\pi(k_n)v\|< c_2\, .
\end{equation}
In fact from

$$\lim_{n\to\infty} \pi(g_nz_nk_n)v = \pi(w)v\quad \hbox{and}\quad
\|\pi(g_n z_n k_n)v\|=\|\pi (z_n) \pi(k_n)v\|$$
we conclude with
(\ref{eq=nz}) that
there are positive constants $c_1',c_2'>0$ such that
$c_1'<\|\pi (z_n) \pi(k_n)v\|<c_2'$ for all $n$.
We use that $\lim_{n\to\infty} z_n =z\in \tilde D$ to obtain
$\pi(z_n)|_V-\pi(z)|_V\to 0$ and our claim
follows.

\par We define $C$ to be the closure of the sequences
$(\pi(z_nk_n)v)_{n\in \N}$ and $(\pi(g_nz_nk_n)v)_{n\in \N}$
in $\H$.
With our previous claims (\ref{eq=nz}) and (\ref{eq=bou}) we obtain
that $C\subset \H -\{0\}$ is a compact subset. But
$C_G=\{g\in G\mid \pi(g) C\cap C\neq \emptyset\}$ contains
the unbounded sequence $(g_n)_{n\in \N}$ and hence is
not compact - a contradiction to Lemma \ref{lem=hm}.
\end{proof}

\begin{rem} Let $(\pi, \H)$ and $v\in \H$ be as in
the theorem. Then we might ask whether the stronger statement
$$\lim_{z\to \partial \tilde\Xi} \|\pi(z)v\| =\infty $$
holds true. For the special case
of $v=v_K\in \H_K$ a $K$-fixed vector this
was established in \cite{KSII}, Th. 2.4.
\end{rem}

\subsection{Domains of holomorphy for the unitary dual
of $G=\Sl(2,\R)$}
Let now $G=\Sl(2,\R)$. With the coordinates of Remark \ref{remidemmi}
we have
$$X_\C=\left (\PP^1(\C)\times \PP^1(\C)\right \backslash \diag(\PP^1(\C)), $$
$$\Xi=X\times \overline X, $$
where $X$ denotes the upper and $\overline X$ the lower halfplane. Then there are two interesting
$G$-domains in $X_\C$
which contain $\Xi$. These are:
\begin{itemize}
\item $S^+=(\PP^1(\C)\times\overline X)\cap X_\C$,
\item $S^-=(X\times \PP^1(\C))\cap X_\C$.
\end{itemize}

\begin{prop} The following assertions hold:
\begin{enumerate}
\item $S^+=G\exp(i\Lambda_{(-1,\infty)}).x_0$,
\item $S^-=G\exp(i\Lambda_{(-\infty,1)}).x_0$.
\end{enumerate}
In particular, $S^\pm$ are maximal $G$-domains in $X_\C$ on which $G$ acts properly.
\end{prop}

\begin{proof} The first two assertions come down to a very elementary
computation; the last one  follows from Proposition \ref{prop=sl2}.
\end{proof}

\begin{rem} As $\gf$ is of Hermitian type,
the $\kf_\C$-module $\pf_\C$ splits into two
inequivalent subspaces $\pf_\C=\pf^+ \oplus \pf^-$
with
$$\pf^\pm =\C \cdot \begin{pmatrix} 1 & \pm i \\ \pm i & 1 \end{pmatrix}\, .$$
Set $P^\pm =\exp (\pf^\pm)$. Then the preimages $\tilde S^\pm$ of
$S^\pm$ in $G_\C$ are given by

 \begin{itemize}
\item $\tilde S^+=G K_\C P^+ $,
\item $\tilde S^-= G K_\C P^-$.
\end{itemize}
\end{rem}

We obtain the following result.

\begin{theorem} \label{th=sl2r} Let $(\pi, \H)$ be an irreducible unitary representation
of $G=\Sl(2,\R)$. Let $v\in \H_K$ be a non-zero vector
and $f_v: G\to \H, \ g\mapsto \pi(g)v$ the corresponding
orbit map. Then the domains of holomorphy
of $f_v$ are given by:

\begin{enumerate}
\item $G_\C$,  if $\pi$ is trivial.
\item $\tilde S^+$,  if $\pi$ is a non-trivial highest weight representation.
\item $\tilde S^-$,  if $\pi$ is a non-trivial
lowest weight representation.
\item $\tilde \Xi$,  if $\pi$ is none of the above, i.e. a unitarizable
principal series.
\end{enumerate}
\end{theorem}

\begin{proof}(i) is clear.
\par\noindent(ii) If $\H_K$ is a highest weight module, then all
its vectors are $K_\C P^+$-finite. Hence $\tilde S^+=GK_\C P^+$ lies in the domain of holomorphy of $f_v$. By the
preceding Proposition $\tilde S^+$ is maximal for proper action and the assertion follows from Theorem \ref{t1}.
\par\noindent (iii) Analogous to (ii).
\par\noindent(iv) For $\tilde \Xi$ to be contained
in the domain of holomorphy we refer to the general result of \cite{KSII}, Th. 1.1. If $\pi$ is $K$-spherical, then
$\tilde \Xi$ is indeed maximal as it follows from \cite{KOS}, Th. 5.1 and Remark (\ref{rem=corr}) below. Finally, the
case of non-spherical principal series is similar to the spherical case (the same proof as in \cite{KOS} applies).
\end{proof}

\begin{rem}\label{rem=corr} (Correcting literature) In the proof of Th. 5.1 in \cite{KOS} there
is an inaccuracy which we wish to correct here. Actually we have to address  the proof of the key result Th. 5.4 in
\cite{KOS}: it asserts for $G=\Sl(2,\R)$ that a spherical function with imaginary parameter blows up at the boundary of
$\Xi$. Now $\partial\Xi=\partial_s\Xi\amalg \partial_n\Xi$. The arguments given for the blow-up at the semisimple
boundary $\partial_s\Xi$ are fine; the ones for the blow-up at $\partial_n\Xi$ are not correct and should be modified.
With the notation of \cite{KOS} we have for $a_r=\begin{pmatrix} r & 0 \\ 0 & {1\over r}\end{pmatrix}\in A$, $r>0$, and
$-1< t< 1$  that
$$P\left(a_r\begin{pmatrix}  1& it \\ 0 & 1\end{pmatrix}.x_0\right)=r^2 +{1\over r^2} - t^2 r^2\, .$$
In particular, if $|t|>1$, then there would exist a sequence $r_n\to r_0$
such that $P\left(a_{r_t}\begin{pmatrix}  1& it \\
0 & 1\end{pmatrix}\right)\to -2^+$. Now we can use the argument given in the proof of Th. 5.4 in \cite{KOS}.
\par Secondly, let us mention that Th. 5.1 and Th. 5.4 are true for all positive definite
spherical functions $\neq {\bf 1}$, not only for those with imaginary parameters as stated -- the argument is literally
the same. \end{rem}

\begin{rem} Theorem \ref{th=sl2r} says that there are different domains of holomorphy for different series of representations. We
expect an analogous result for arbitrary semisimple Lie groups and intend to return to this topic elsewhere.\footnote{Added in proof:
This is now established, see \cite{K}.} Let us
mention that the crown is maximal for the class of non-trivial spherical unitary representations  of $G$ by \cite{KOS},
Th. 5.4 and the remark above. \end{rem}

\section{Holomorphic extension of spherical functions}

This section is a short essay on  spherical functions
on $X$ which highlights their natural holomorphic extension to the crown
$\Xi$. Here $G$ is arbitrary semisimple (within our self-imposed
restrictions).

\par As always some notations upfront.
For $\alpha\in \Sigma$ write $\gf^\alpha$ for the corresponding
root space. Choose a positive system $\Sigma^+$
and define $\nf=\sum_{\alpha\in \Sigma^+} \gf^\alpha$. Set
$N=\exp \nf$.
The Iwasawa decomposition $G=NAK$ yields
the analytic diffeomorphism
\begin{equation}\label{eq=iwa}
N\times A\mapright{\simeq} X,\  \ (n,a)\to na.x_0\, .\end{equation}
In particular, every $x\in X$ can be uniquely
written as $x=n(x) a(x).x_0$ with $n(x)\in N$ and $a(x)\in A$ both
depending analytically on $x$.

\par Let $N_\C =\exp\nf_\C$. If we complexify the Iwasawa decomposition of $X$ we obtain
a Zariski open  subset $N_\C A_\C .x_0\subsetneq X_\C$
which contains the crown, i.e.
\begin{equation}\label{eq=cin} \Xi\subset N_\C A_\C .x_0 \end{equation}
(see \cite{KSI} for classical groups  and \cite{Hu}\footnote{We would like to caution the reader that the proof in
\cite{Hu} is severely wrong; a correct proof -- unfortunately not emphasized --  appeared later in \cite{FH}.}
as well as \cite{Ma} in general).
Let us mention that $\Omega\subset \af$ is a maximal domain
for the inclusion (\ref{eq=cin}) to hold, i.e.
$G\exp(i\pi \tilde \Omega/2)\not\subset N_\C A_\C.x_0$ for
any domain $\tilde \Omega\subset \af$ strictly containing $\Omega$
(cf.\ \cite{Ba} and \cite{KSII}, Th. 2.4 with proof).
\par

Define the finite $2$-group $F=A_\C\cap K_\C =T\cap K$ and record that
the map
$$N_\C \times A_\C/ F\mapright{\simeq} N_\C A_\C.x_0, \ \
(n,aF)\mapsto na.x_0$$
is biholomorphic. It follows that
each element $z\in N_\C A_\C.x_0$ can be uniquely
expressed as $z=n_\C(z) a_\C (z).x_0$ with
$n_\C(z)\in N_\C$ and $a_\C(z)\in A_\C/F$ both holomorphic
in $z$. One obtains an $N$-invariant
holomorphic assignment
\begin{equation} a_\C: \Xi \to A_\C/ F
\end{equation}
We have already remarked that $\Xi$ is contractible and this yields
that $a_\C$ lifts  to a holomorphic map $\Xi \to A_\C$,
as well  denoted by $a_\C$, such that $a_\C (x_0)=e$. Likewise
there is a holomorphic logarithm $\log a_\C : \Xi \to \af_\C$
extending $\log a: X\to \af$.  In particular,
for all $\lambda\in \af_\C^*$ we can define
the holomorphic $\lambda$-power of $a_\C$  by

$$a_\C(z)^\lambda=e^{\lambda (\log a_\C (z))} \qquad (z\in \Xi)$$

\par We would like to mention the complex convexity theorem
(\cite{GKII}, \cite {KOt}) which states
that
\begin{equation} \label{cc}
\Im \log a_\C (G\exp(iY).x_0)={\rm co} (W.Y) \qquad (Y\in \pi\Omega/ 2) \end{equation}
with ${\rm co}(\cdot)$ denoting the convex hull of $(\cdot)$.
As a consequence we obtain a refinement of the inclusion
(\ref{eq=cin}):
\begin{equation} \Xi \subset N_\C AT_\Omega.x_0
\end{equation}

\par For $\alpha\in\Sigma$ let us define  $m_\alpha=\dim \gf^\alpha$ and note
that the multiplicity assignment $\alpha\mapsto m_\alpha$ is $W$-invariant.
As usual we set $\rho={1\over 2}\sum_{\alpha\in \Sigma^+} m_\alpha \alpha$.
Motivated by our previous discussion we define the spherical function
with parameter $\lambda\in \af_\C^*$ ab initio as a holomorphic function on
$\Xi$:

\begin{equation} \phi_\lambda(z)=\int_K a_\C (kz)^{\rho +i\lambda} \, dk \qquad (z\in \Xi)
\end{equation}

Of later relevance for us
will be the doubling formula for
spherical functions (\cite{KSI}, Th. 4.2).
For the convenience of the reader we briefly recall
the short argument.
We translate the inclusion (\ref{eq=cin}) into representation theory:
Using the compact realization of a spherical minimal principal
series module $(\pi_\lambda,\H_\lambda)$ one
shows that the orbit map of a spherical vector $v_{\lambda}\in\H_\lambda$
\begin{equation}\label{or1}
F: X\to \H_\lambda,\ \  x\to\pi_\lambda(x)v_\lambda
\end{equation}
extends to a holomorphic map
\begin{equation} \label{or2}
F: \Xi\to \H_\lambda\ \ z\to\pi_\lambda(z)v_\lambda
\end{equation}
see \cite{KSI}, Prop. 4.1. This allows us to express the
spherical function $\phi_\lambda$ as a holomorphic matrix
coefficient
$\phi_\lambda(z)=\langle v_\lambda,\pi_\lambda(z)v_\lambda \rangle$ for
$z\in \Xi$ (where we have adopted the physicist's
convention
that sesquilinear pairings are linear on the right hand side,
and anti-linear on the left hand side).

Now let $z\in AT_\Omega$ such that $z^2\in AT_\Omega$ and
observe that
$$
\phi_\lambda(z^2)=\langle \pi_\lambda^*(\overline{z^{-1}})
v_\lambda,\pi_\lambda(z)v_\lambda\rangle
$$
with $\pi_\lambda^*$ the conjugate contragredient representation.
It follows that $\phi_\lambda|_{AT_\Omega}$ extends to a
holomorphic function on $AT_\Omega^2$.
In particular we see from this formula in the case of the unitary
spherical minimal principal series $\lambda\in \af^*$ that the
function $\phi_{\l}$ is positive on $T_\Omega^2$,
and that for all $x=gt.x_0\in \Xi$ (recall the notation
$v_{\l}^x:=\pi_\lambda(x)v_{\l}$):
\begin{align}\label{eq:l2norm}
\langle v_{\l}^x,v_{\l}^x\rangle&=\phi_{\l}(t^2)\\
\nonumber &=\int_K |a_\C(kt)^{2(\rho+i\lambda)}|dk.
\end{align}
\section{Sharp uniform lower bound for holomorphically extended orbit maps
of spherical representations}
Given a non-trivial unitary spherical representation
$(\pi, \H)$ of $G$ with normalized $K$-spherical
vector $v_K$ we wish to control the norm of the
holomorphically extended orbit map
$$F_\pi: \Xi\to \C, \ \ z\mapsto \pi(z)v_K$$
in two aspects:
\begin{itemize}
\item For $z\in \Xi$ sufficiently close to $\partial_d\Xi$ we
are aiming to give optimal lower bounds for $\|F_\pi(z)\|$ uniform
in the representation parameter $\lambda(\pi)\in \af_\C^*$;
\item For fixed $\pi$ we  are looking for optimal upper
bounds of $\|F_\pi(z)\|$ for $z$ approaching the distinguished
boundary.
\end{itemize}
In view of the fundamental identity (\ref{eq:l2norm}) we can translate
the problems above into growth behavior of analytically continued
spherical functions. In this section and the next we will address
these two aspects. We begin with the uniform lower bounds.

Fix a distinguished boundary point $t=z_j=\exp(i\pi\omega_j/2k_j).x_0$
of the crown domain. For $0<\e<1$ set
$$t_\e =\exp(i(1-\e)\pi\omega_j/2k_j).x_0\, .$$
The objective of this section is to provide sharp lower estimates
for $\phi_\lambda(t_\e^2)$
which are uniform
in $\e$ and $\lambda\in \af^*$.
Our approach is based on the doubling identity
(\ref{eq:l2norm})
which implies
that
\begin{equation}\label{ble} \phi_\lambda(t_\e^2)\geq \int_U
|a_\C (kt_\e)^{2(\rho +i\lambda)}|
\ dk \end{equation}
where $U$ is any neighborhood of $e\in K$.
It turns out that the desired  estimate will depend on the nature of
the distinguished boundary point $t=z_j$, in particular
whether $z_j$ is minuscule or not.
We will treat the minuscule case first and later reduce
the general case to the minuscule situation.

\par If $t=z_j$ is minuscule,  then $G=Z(t^4)$ by Lemma \ref{lem=m}, i.e.
$t^4$ is central.
The following lemma, especially seen in the context of
(\ref{cc}), is quite remarkable.

\begin{lemma} \label{lem=key}Let $t=z_j$ be a minuscule boundary point
and $U$ a connected and simply connected compact
neighborhood of $e\in K$ such that $Ut\subset N_\C A_\C .x_0$.
Then for all $k\in U$ the middle projection
$a_\C(kt)\in A_\C$ is well defined and
we have $a_\C(kt)=r(kt)t$ with $r(kt)\in A$ continuously
depending on $k$.
\end{lemma}

\begin{proof} The assertion of the lemma is local and thus it is
no loss of generality to assume that $G\subset G_\C$ with
$G_\C$ simply connected. In particular the Cartan
involution $\theta:G\to G$ extends to a holomorphic involution on $G_\C$,
again denoted by $\theta$. We notice that $G_\C^\theta=K_\C$.
Likewise $G_\C$ admits a complex conjugation $g\mapsto \overline g$
with respect to $G$
\par Fix $k\in K$. Then $kt=nak'$ for some $n\in N_\C$, $a\in A_\C$ and
$k'\in K_\C$. Define $x:=kt\theta(kt)^{-1}$ and note that

\begin{equation} \label{*}
x=kt^2k^{-1}=na^2\theta(n)^{-1}\, .\end{equation}
On the other hand, as $t^4$ is central,
\begin{equation}\label{**}
t^{-4}x=kt^{-2}k^{-1}=\overline{kt^2k^{-1}}
=\overline{x}=\overline {n}\overline{a}^2\theta(\overline n)^{-1}\, .
\end{equation}
Combining the information of (\ref{*}) and (\ref{**}) yields
$$t^4 \overline {n} \overline {a}^2 \theta(\overline {n})^{-1}=
x=na^2 \theta(n)^{-1}\, .$$
Once more  we use the fact $t^4$ is central and obtain
$$t^4 \overline{a}^2 = \underbrace{(\overline {n}^{-1} n)}_{\in N_\C}
a^2
\underbrace{\theta(n^{-1}\overline {n})}_{\in \theta(N_\C)}\, .
$$
Bruhat implies that $\overline{n}=n$. Consequently
$t^4=a^2\overline{a}^{-2}$, and this forces $a=r(tk) t$ for some
$r(tk)\in A$.
\end{proof}

Choose $0<\e_0<1$ small enough such that
$Ut_\e\subset N_\C A_\C .x_0$ for all $\e\in (0,\e_0)$.
In particular $a_\C(kt_\e)$ is well defined
for all $k\in U$ and $\e\in (0, \e_0)$.
As $a_\C (kt)\in A t$ for all $k\in U$ by the lemma, linear
Taylor approximation yields that
there are balls  $B_r, B_{r'}$ in $\af$ centered at $0$ with radii $r, r'>0$
such that

\begin{equation} a_\C(kt_\e)\in  t \exp(B_r) \exp(i\e B_{r'})
\end{equation}
for all $k\in U$ and $\e\in (0,\e_0)$.
Thus it follows that there exists a constant $c>0$ such that
for all $\lambda\in \af^*$, $k\in U$ and $\e\in (0,  \e_0)$
the estimate
\begin{equation}\label{zztop}
|a_\C(kt_\e)^{2(\rho+ i\lambda)}|\geq  c e^{\lambda(\pi\omega_j) -
r''\e |\lambda|}
\end{equation}
holds for some $r''\geq r'$.

\begin{prop}\label{prop=before} Let $t=z_j=\exp(i\pi \omega_j/2).x_0$ be a minuscule
boundary boundary point of $T_\Omega$.
Then there exist constants
$\e_0\in(0,1)$ and $R>0$, $C>0$
such that
\begin{equation}
\phi_\lambda(t_\e^2)\geq C\max_{w\in W}e^{ \pi\l(w\omega_j)(1 - R\e)},
\end{equation}
for all $\l\in\af^*$ and $\e\in(0,\e_0)$.
\end{prop}
\begin{proof} According to Harish-Chandra one has
$\phi_\lambda=\phi_{w\lambda}$ for all
$\lambda$. Thus it is no loss of generality to  assume that
$\lambda(\omega_j)=\max_{w\in W} \lambda(w\omega_j)$.
We notice that $||\lambda||:=\max_{w\in W} \lambda(w\omega_j)$
defines a norm on $\af^*$. Hence, by the equivalence of norms
on Euclidean spaces, there exist a constant $d>0$ such that
$|\cdot |\leq d||\cdot ||$.
Now the the assertion
follows from (\ref{zztop}) and the
basic lower estimate (\ref{ble}).
\end{proof}

Let us now turn to the general case where
$t=z_j=\exp(i\pi\omega_j/ 2k_j).x_0$ is an arbitrary
extremal boundary point of $T_\Omega$. We recall
the groups $G_j=Z(t^4)$ with Lie algebra $\gf_j$. The
main result of this section is:

\begin{thm}\label{lowerth}
Let $t=\exp(i\pi  \omega_j/2k_j)$ be an extremal boundary point of
$T_\Omega$. Then there exist constants $\e_0\in(0,1)$
and $R>0$, $C>0$ such that
\begin{equation}
\phi_\lambda(t_\e^2)\geq
C \e^{(\operatorname{dim}\gf-\operatorname{dim}
\gf_j)/4}\max_{w\in W} e^{\pi\l(w\omega_j)(1-R\e)}
\end{equation}
for all $\l\in\af$, and for all $\e\in(0,\e_0)$.
\end{thm}
\begin{proof} First, for $\omega_j$ minuscule one
has $\gf_j=\gf$ and the assertion follows from
Proposition \ref{prop=before} above.
The general case will be reduced to this situation.
\par We begin with some remarks on the reductive Lie algebra
$\gf_j$. Recall that $\gf_j$ is $\theta$-stable and hence
$\gf_j=\kf_j\oplus \pf_j$ with $\kf_j=\kf\cap \gf_j$ and
$\pf_j=\pf\cap \gf_j$.
By definition
$\af\subset\gf_j$ and hence $\af$ is maximal abelian
in $\pf_j$. Let $\Sigma_j=\Sigma(\gf_j,\af)$ be the corresponding
reduced root system and $\Omega_j\subset \af$ the associated
polyhedron.
A quick look at our classification of the $\gf_j$'s shows that
$\gf_j$ is simple modulo a compact ideal. Hence $\Sigma_j$ is
irreducible and $\omega_j/k_j$ becomes a minuscule boundary
point of $\Omega_j$.

\par Write $\kf_j^\perp$ for the orthogonal complement to $\kf_j$
in $\kf$ with respect to the Cartan-Killing form of $\gf$.
Let $V_j,V_j'$ be small balls around $0$ in $\kf_j$,$\kf_j^\perp$
such that the map
$$V_j\times V_j'\to K,\ \  (v,v')\mapsto \exp(v)\exp(v')$$
is a diffeomorphism. Set $U_j=\exp(V_j)$, $U_j'= \exp(V_j')$
and define $U=U_j U_j'$. Then $U$ is a connected and simply
connected neighborhood of $e$ in $K$. We assume that
$Ut\subset N_\C A_\C.x_0$ and choose $\e_0>0$
such that $Ut_\e\subset N_\C A_\C.x_0$ for
$\e\in (0,\e_0)$ holds in addition.

\par Our previous discussion combined with Lemma \ref{lem=key}
implies that
\begin{equation}\label{eq=1} a_\C (kt)\in r(kt) t \qquad \forall k\in U_j'
\end{equation}
and $r(kt)\in A$ depending continuously on $k$.
Next consider the map $\psi: U\to A_\C , k\mapsto a_\C(kt)$.
We claim that
\begin{equation} \label{eq=2}
d\psi(e)= 0\,  .
\end{equation}
In fact, this is well known, and follows from
$\operatorname{pr}_{\af_\C}(\operatorname{Ad}(t^{-1})\kf)=\{0\}$ with
$\operatorname{pr}_{\af_\C}: \gf_\C \to \af_\C$ the linear projection
along $\nf_\C \oplus \kf_\C$.
\par Using the information of (\ref{eq=1}) and
(\ref{eq=2}), linear Taylor approximation
yields constants $r,r'>0$ such that
\begin{equation}\label{eq=3}
a_\C(\exp(v)\exp(v')t_\e)\in t \exp(B_r)\exp(i(\e +\|v'\|^2)  B_{r'})
\end{equation}
for all $(v,v')\in V_j\times V_j'$ and $\e\in(0,\e_0)$.
It follows from equation (\ref{ble}) that there exists
a constant $c>0$ such that
\begin{equation} \label{ble'}\phi_\lambda(t_\e^2)\geq c
\int_{V_j} \int_{V_j'}
|a_\C (\exp(v)\exp(v')t_\e)^{2(\rho +i\lambda)}|
\ dv dv' \end{equation}
for all $\lambda\in \af^*$ and $\e\in (0,\e_0)$.
Thus if we choose $V_j'$ to be ball
of radius $\sqrt{\e}$, then (\ref{eq=3}) and (\ref{ble'})
yield  constants  $r'', c'>0$ such that
 $$\phi_\lambda(t_\e^2) \geq c' \e^{ (\operatorname{dim} \kf_j^\perp)/2}
e^{\lambda(\pi\omega_j/k_j) -r''\e|\lambda|}\, $$
for all $\lambda\in \af^*$ and $\e\in (0,\e_0)$ (note that the $\e$
-dependence of $V_j'$ is incorporated in the factor
$\e^{ (\operatorname{dim} \kf_j^\perp)/2}$).
We observe that $\operatorname{dim} \gf-\operatorname{dim}
\gf_j= 2 \operatorname{dim} \kf_j^\perp$ and finish the proof
with the same argument for Proposition \ref{prop=before} before.
\end{proof}

\begin{rem} We consider the lower estimate in
Theorem \ref{lowerth} is optimal. It is for
the following reason: the crucial point in the above argumentation
was the fact that $d\psi(e)=0$, to be very
precise it was the fact $d (\Im \log \psi)(e)|_{\kf_j^\perp}=0$
which entered.
This is actually the best one can hope for as
the second derivative $d^2(\Im \log \psi)(e)$ is already
non-degenerate on $\kf_j^\perp\times \kf_j^\perp$. In fact,
fix $\lambda \in \af^*$ regular, and
set $F_\lambda =\lambda \circ \Im\log \psi$.
Then \cite{DKV}, pp. 343--346, implies
$$d^2 F_\lambda (e)(Z,W)=-{1\over 2}\sum_{\alpha\in \Sigma^+}
\langle \alpha, \lambda\rangle \Im (1-t^{-2\alpha}) \langle  Z_\alpha, W_\alpha
\rangle\qquad (Z,W\in \kf)$$
where $Z_\alpha$, resp. $W_\alpha$
is the orthogonal projection of $Z$, resp. $W$,  onto $\kf\cap (\gf^\alpha+
\gf^{-\alpha})$.  In particular if $\alpha\in \Sigma\backslash \Sigma_j$,
then $\Im (1-t^{-2\alpha})\neq 0$. It follows that $d^2 F_\lambda(e)$
and hence $d^2(\Im \log \psi)(e)$
is non-degenerate on $\kf_j^\perp\times\kf_j^\perp$.
\end{rem}
\section{Sharp upper bound for holomorphically extended orbit maps
of spherical representations}
In this section we consider the problem to give an upper estimate
for the square norm of the holomorphic extension of the orbit map
$\Xi\ni x\to v^x\in\mathcal{H}$.
Recall from (\ref{eq:l2norm}) that
\begin{equation}\label{eq:orbnorm}
(v^{x}_{i\l},v^{x}_{i\l})=\phi^X_{i\l}(t^2.x_0)
\end{equation}
if $x=gt.x_0\in\Xi$ with $g\in G$ and $t\in T_\Omega$.
Here $\phi^X_\mu$ denotes the elementary spherical function on
$X$ with spectral parameter $\mu\in\af_\C^*$.
Therefore we concentrate on the question of estimating the
singular behavior of the holomorphic extension of the restriction
of the elementary spherical function $\phi^X_{\mu}|_{A}$ to
$A_\C\supset A$ when we approach $t^2.x_0\in T.x_0\subset A_\C.x_0$
where $t=t(\eta)=\exp(i\pi\eta/2)$ with $\eta=\omega_j/k_j\in\af_\C$
an extremal boundary point of $\Omega$ in $\Omega\cap C$.
Thus we are interested in the singular expansion in
$\epsilon$ of the pull-back of (the holomorphic continuation of)
spherical functions $\phi_\mu$ via the embedding
$\mathbb{D}^\times\ni\epsilon\to A_\mathbb{C}.x_0$
of a small punctured disk $\mathbb{D}_r^\times
=\{\e\in\C\mid 0<|\e|< r\}$ given by
$\epsilon\to t_\e^2.x_0$. For $\mu$ fixed
the restriction of $\phi^X_\mu$ has a convergent logarithmic
singular expansion at $\epsilon=0$.
This means that there exists a finite set
$S\subset\C\times\mathbb{Z}_{\geq0}$
such that if $(s,l),(s^\prime,l^\prime)\in S$
then $s-s^\prime\not\in\mathbb{Z}\backslash\{0\}$
and such that we have a unique decomposition
(for $\epsilon$ varying in any sector
$S_{r,\theta_1,\theta_2}=\{\e\in\mathbb{D}_r^\times\mid \theta_1<
\operatorname{arg}(\e)<\theta_2\}$
of $\mathbb{D}_r^\times$) of the form
\begin{equation}\label{eq:leadexp}
\phi^X_{\mu}(t^2_\e.x_0)=\sum_{(s,l)\in S}
\e^s\log^l\e f_{s,l}(\e)
\end{equation}
where each $f_{s,l}$ is holomorphic on $\mathbb{D}_r$
and such that $f_{s,l}(0)\not=0$.

The projection of the set $S$ on the first factor
$\C$ is called the set of exponents of the pull
back of $\phi^X_\mu$ to $\mathbb{D}^\times_r$. In our case
this set will always belong to $\mathbb{R}$. The minimum
of this set is denoted by $s^X_{\eta,\mu}$ and is called
the \emph{leading exponent} of the
singular expansion of the pull back of $\phi^X_\mu$ to
$\mathbb{D}^\times_r$. We call the largest $l\in\mathbb{Z}_{\geq0}$
such that the $(s^X_{\eta,\mu},l)\in S$ the
\emph{logarithmic degeneracy} of the leading exponent.

So our problem boils down to the determination of
the leading exponent $s^X_{\eta,\mu}$ at $\e=0$ of the pull back of
$\phi^X_\mu$ on the distinguished embedded punctured disk given above,
and its logarithmic degeneracy. In the Appendix \ref{app:exp}
we define an appropriate
notion of the exponent of a regular holonomic system of
differential equations and using the basic properties of these
exponents we compute the exponents of $\phi_\mu$ at the extremal
boundary points $\eta$ of $\Omega\cap C$ for $\phi_\mu$ a solution
of a more general system of differential equations, namely the system
of hypergeometric equations associated with the root system $\Sigma$.
This system of
equations is a parameter deformation of the system of equations
for the restriction for the elementary spherical functions
$\phi^X_\mu$ to $A_\C.x_0$. This deformation is an essential
ingredient for the computations of the exponents.
These results imply the following:
\begin{thm}\label{thm:leadexp}
We use the notations as introduced above.
Consider the functions
$s_\eta(m)$ and
$d_\eta(m)$ of the multiplicity parameters $m=(m_\a)$ as
listed in the table in  Theorem \ref{thm:table}.
Suppose that the Riemannian
symmetric space $X$ has root system $\Sigma$ and root
multiplicity parameters $m^X=(m^X_\a)$, then we put
$s^X_\eta:=s_\eta(m^X)$ and $d_\eta^X=d_\eta(m^X)$.

For all $\mu\in\af_\C^*$ we have
$s^X_{\eta,\mu}\geq s^X_\eta$, and if
$s^X_{\eta,\mu}= s^X_\eta$ then $d^X_\eta$ is an upper bound for
the logarithmic degeneracy of $s^X_{\eta,\mu}$.
\end{thm}
We postpone the proof of this theorem in the general case
to the Appendix \ref{app:exp}. For the complex cases (i.e.
when $X$ is a Riemannian symmetric space of type IV) the
proof will be given below.

As an immediate consequence of theorem \ref{thm:leadexp}
we have:
\begin{thm}\label{thm:esteta}
We use the notations as introduced above.
Given an extremal boundary point $\eta=\omega_j/k_j$ of
$\Omega$ we consider $t_\epsilon=\exp(i\pi\eta/2)\in A_\C$.
Fix $-\pi<\theta_1<\theta_2<\pi$.
Let $\mu\in\af_\C^*$, then there exist constants
$r>0$, $C>0$ such that for all $\e\in S_{r,\theta_1,\theta_2}$:
\begin{equation}\label{eq:mainest}
|\phi^X_\mu(t^2_\epsilon.x_0)|\leq C\e^{s^X_\eta}|\log(\e)|^{d^X_\eta}
\end{equation}
where $s_\eta^X=s_\eta(m^X)$ and $d_\eta^X=d_\eta(m^X)$ for the
functions $s_\eta$ and $d_\eta$ listed in
Theorem \ref{thm:table}.
\end{thm}
For later applications it is useful to have a slightly weaker
but more handy version of the estimate above. Let us define

\begin{equation}\label{def=ds} s^X:=\max_\eta s_\eta^X\quad\hbox{and}
\quad
d^X:= \max_{\eta:s_\eta^X=s^X} d_\eta^X\, .
\end{equation}
The theorem above combined with the maximum principle of
holomorphic functions then yields:

\begin{thm}\label{thm=se} For each $\mu\in \af_\C^*$ there
exists a constant $C=C(\mu)>0$ such that
for all $Y\in \partial\Omega$ and $0<\e<1$
\begin{equation}
|\phi^X_\mu(\exp(i(1-\e) \pi Y).x_0)|\leq C\e^{s^X}|\log(\e)|^{d^X}\, .
\end{equation}
\end{thm}

Since $\Xi$ only depends on the isogeny class of $G$ (Remark
\ref{rem=bd}(c)) it suffices to do the analysis in the situation
where $G_\C$ is simply connected. In addition we assume the
restricted root system $\Sigma$ of $X$ to be irreducible. Recall
that twice the character lattice of $A_\C$ is equal to the weight
lattice of the restricted root system $\Sigma^l$. The categorical
quotient $W\backslash A_\C$ (as well as $W\backslash A_\C/F$) is
affine space.

There are two special cases which can be treated by direct methods,
the real rank one case and the complex case. It is both instructive and
useful to consider these cases first before going to the
general case which is treated in the Appendix \ref{app:exp}.
\subsubsection{The real rank one case} This case was treated in
detail in \cite[Theorem 5.1]{KSI} but it is useful to discuss
the difference between the approach used in \cite{KSI} and the
approach in the present paper.

Let $G$ be a real semisimple group with real rank one. Then
$\Sigma=\Sigma(\mathfrak{g},\mathfrak{a})$ is of the form
$\Sigma=\{\pm\alpha\}$ (reduced case) or
$\Sigma=\{\pm\alpha/2,\pm\alpha\}$ (non-reduced case).
Let $\eta\in\mathfrak{a}$ be such that $\alpha(\eta)=1$, so that
$\eta$ is at the boundary of $\Omega$.
We put $q:=\dim\mathfrak{g}^\alpha$ and
$p:=\dim\mathfrak{g}^{\alpha/2}$ (so that $q\geq 1$, with
$p=0$ (reduced case) or $p\geq q$ (non-reduced case);
this is convenient but admittedly a bit unconventional
in the reduced case).

In \cite{KSI} the Harish-Chandra integral representation for the
spherical function $\phi_\lambda$ is analyzed directly to obtain
the precise asymptotic behaviour of the holomorphic extension of
$\phi_\lambda(\exp(i\pi(1-\epsilon)\eta))$ as $\epsilon\to 0$ if
$\lambda\in i\mathfrak{a}^*$. The result of \cite{KSI} says that
for $\lambda(\alpha^\vee)\in i\mathbb{R}$ one has
\begin{equation}\label{eq:asexplog}
\phi_\lambda(\exp(i\pi(1-\epsilon)\eta))\asymp |\log\epsilon|
\end{equation}
if $q=1$, and
\begin{equation}\label{eq:asexp}
\phi_\lambda(\exp(i\pi(1-\epsilon)\eta))\asymp \epsilon^{1-q}
\end{equation}
if $q>1$.

The method of the present paper is not based on the analysis of
the integral representation of the spherical function
$\phi_\lambda$ but rather on the analysis of the radial system
of eigenfunction equations for $\phi_\lambda$ with respect to the commutative
algebra of $G_\mathbb{C}$-invariant differential operators on $\Xi$.
In the case at hand, this amounts to the well known fact that
\begin{equation}
\phi_\lambda(\exp(i\pi(1-\epsilon)\eta))=F(a,b,c;z)
\end{equation}
with
\begin{align*}
a&=\lambda(\alpha^\vee)+p/4+q/2\\
b&=-\lambda(\alpha^\vee)+p/4+q/2\\
c&=1/2+p/2+q/2\\
z&=(1-\cos(\pi(1-\epsilon)))/2=1-\pi^2\epsilon^2/4+\dots
\end{align*}
where $F(a,b,c;z)$ (with $a,b,c\in\mathbb{C}, |z|<1$)
denotes the Gauss hypergeometric function. This function
is the unique holomorphic solution of
the hypergeometric differential equation
\begin{equation}\label{eq:Ghyp}
z(1-z)\frac{d^2F}{d^2z}+(c-(a+b+1)z)\frac{dF}{dz}-abF=0
\end{equation}
on the unit disc, normalized by $F(a,b,c;0)=1$.
The hypergeometric differential equation is a second order
Fuchsian equation with regular singularities at $z=0, z=1$ and
$z=\infty$. The exponents of (\ref{eq:Ghyp}) at $z=1$ are easily
seen to be equal to
$0$ and $c-a-b=1/2(1-q)$. Hence the problem of finding the singular
expansion of $\phi_\lambda(\exp(i\pi(1-\epsilon)\eta))$ at
$\epsilon=0$ is now transformed to a simple exercise on
the analytic continuation of the solutions of (\ref{eq:Ghyp}).
If $\phi_\lambda(\exp(i\pi(1-\epsilon)\eta))=F(a,b,c;z)$
is holomorphic at $z=1$ (i.e. at $\epsilon=0$) then $F(a,b,c;z)$
extends holomorphically to $\mathbb{C}$.
Since $\infty$ is a regular singular point of (\ref{eq:Ghyp})
this implies
that $F(a,b,c;z)$ is a polynomial in $z$ in this case (essentially a
so-called Jacobi polynomial). The hypergeometric series
terminates iff $a$ or $b$ is a non-positive integer (since
$c$ is not a negative integer in our situation), i.e. iff
$\lambda(\alpha^\vee)=\pm(p/4+q/2+n)$ with
$n\in\mathbb{Z}_{\geq 0}$. Hence for these values of $\lambda$
the function $\phi_\lambda(\exp(i\pi(1-\epsilon)\eta))$ is asymptotic
to a constant (for $\epsilon\to 0$), and for all other values of
$\lambda$ we find $(\ref{eq:asexplog})$ or $(\ref{eq:asexp})$
(whichever is relevant).

In the general case we need to deal with Harish-Chandra's radial system of
eigenfunction equations for the $G$-invariant differential operators, which is a
holonomic system of differential equations on $W\backslash A_\mathbb{C}/F$ with
regular singularities along the discriminant locus.
This gives rise to various questions and problems which are addressed in
the Appendix \ref{app:exp}.
First of all it is not clear in general how to define a set
of singular exponents at a point of the singular locus of a regular
singular holonomic system. In Subsection \ref{sub:expAR}
we give an appropriate definition of the notion of exponents in the special
case of a holonomic system of differential equations with regular singularities
along an affine hyperplane arrangement. We apply this definition to the pull back
to $\mathfrak{a}_\mathbb{C}$ of the Harish-Chandra system of equations.
The next problem is the actual computation of the set of exponents at $\eta$.
We do not know how to do this directly in an algebraic way in general;
instead we use the relation between the exponents and monodromy of the solutions
of the system of equations and, in a crucial way, the parameter deformation of
the Harish-Chandra system provided by the so-called hypergeometric system
of differential equations for root systems. Via the monodromy representation
affine Hecke algebras play a role in dealing with these computations.
A final problem in the higher rank
case is the selection of those exponents among the full set of exponents
at $\eta$ which may be involved in the expansion of $\phi_\lambda$ at $\eta$,
and among those, the leading exponent.
Again some basic representation theory of the affine Hecke algebra plays
an important role at this stage. The outcome is remarkable;
the leading exponent is attached to a specific irreducible character of
the isotropy group $W^a_\eta$ of $\eta$ in the affine Weyl group $W^a$.
We call this irreducible character the \emph{leading character} at $\eta$.
It turns out that the leading character (for generic $\lambda$) depends on
the local geometry of $\Omega$ at $\eta$, but not on the root multiplicities.
\subsubsection{The complex case}
We consider the complex case $X=G/U$ where $G$ is the
connected simply connected complex simple group and $U$ its
maximal compact subgroup. As is well known by the work of Harish-Chandra
\cite{HC0} the spherical functions are of an elementary nature
in this case, allowing us to analyze the desired exponents directly.
The restricted root system $\Sigma$ of $X$ is equal to twice the root
system of $G$, and all
root multiplicities are equal to $2$. We introduce the Weyl
denominator $\delta$ on $A_\C/F$ by
\begin{equation}
\delta(a)=\prod_{\a\in\Sigma_+}(\a(a)-\a(a)^{-1})
\end{equation}
and we denote by $A_\C^{\operatorname{reg}}/F$ the complement of
the set $\d=0$ in $A_\C/F$. The algebra
$\Ri_X$ of radial parts of invariant differential operators on $X$
consists of the differential operators on
$A_\C^{\operatorname{reg}}/F$ of the form
\begin{equation}
\Ri_X=\{\delta^{-1}\circ\partial(p)\circ\delta\mid p\in\C[\af^*]^W\}
\end{equation}
The Harish-Chandra isomorphism $\gamma_X:\Ri_X\to\C[\af^*]^W$ is
given by
$\gamma_X(\delta^{-1}\circ\partial(p)\circ\delta)=\partial(p)$.
Now $\phi^X_\mu|_{A_\C}$ satisfies the following system of
eigenfunction equations
\begin{equation}\label{eq:cpx}
D\phi=\gamma_X(D)(\mu)\phi,\ \forall D\in\Ri_X
\end{equation}
This is a $W$-equivariant system of differential equations on
$A_\C^{\operatorname{reg}}/F$. It is also equivariant for the
action of the $2$-group $F=A_\C\cap U$. Therefore we can view
(\ref{eq:cpx}) as a system of differential equations on
$W\backslash A_\C/F-\{d=0\}$, where $d=\delta^2$ is the
discriminant of $W$, viewed as a polynomial on $W\backslash
A_\C/F$. We call $\L_X$ the sheaf of local
solutions of (\ref{eq:cpx}).

A general local solution to this set of equations is of the form
\begin{equation}\label{eq:factor}
\phi=\d^{-1}\psi
\end{equation}
where $\psi$ is a local solution of the constant coefficient
system
\begin{equation}\label{eq:stein}
\partial(p)\psi=p(\mu)\psi,\ \forall p\in\C[\af^*]^W
\end{equation}
Let $\exp:\af_\C\to A_\C$ denote the exponential
map such that $\exp(2\pi iX)=1$ iff $X\in
Q(\Sigma^\vee)$. Consider the covering map
\begin{equation}\label{eq:cov}
\pi:\af_\C^{\operatorname{reg}}\to
W\backslash A_\C^{\operatorname{reg}}/F
\end{equation}
where $\af_\C^{\operatorname{reg}}$ is the complement in $\af_\C$
of the set of affine root hyperplanes, the zero sets of the affine
roots $a=\a-n$ (with $\a\in\Sigma$ and $n\in\mathbb{Z}$), and
$\pi$ is given by $\pi(X)=W\exp(\pi i X)F$. The following
proposition is well known.
\begin{prop}\label{prop:dim}
The space of solutions of (\ref{eq:cpx}) on a nonempty open ball
$U\subset W\backslash A_\C^{\operatorname{reg}}/F$ consists of holomorphic
functions and has dimension $|W|$ (independent of $\mu$). Let
$V\subset\pi^{-1}(U)$ be a connected component. The pull
back of a local solution of (\ref{eq:cpx}) on $U$ via $\pi|_V$
extends to a global holomorphic function on
$\af_\C^{\operatorname{reg}}$.
\end{prop}
\begin{proof}
We use the general form (\ref{eq:factor}) of the local solutions.
By a well known result of Steinberg \cite{S} the global solution
space of (\ref{eq:stein}) on
$\af_\C$ has dimension equal to $|W|$ (independent of $\mu$) and
consists of entire functions. On the other hand the left
ideal in the ring of differential operators with holomorphic
coefficients on $\af_\C$ generated by the operators
$\partial(p)-p(\lambda)$ (with $p\in\C[\af^*]^W$) is cofinite, a
complement being generated by constant coefficient operators
$\partial(q)$ where $q$ is running over a set of polynomials
representing a basis of the coinvariant algebra. Hence the local
solution space is at most of dimension $|W|$. The proposition
follows.
\end{proof}
\begin{cor}
The system (\ref{eq:cpx}) is holonomic of rank $|W|$. Upon
choosing a base point $p\in\af_\C^{\operatorname{reg}}$
we may view the monodromy representation as a representation of
the group of deck transformations of the covering map $\pi$, which is
the affine Weyl group $W^a=W\ltimes Q(\Sigma^\vee)$ acting on
$\af_\C$. The restriction of the monodromy representation to $W$ is
equivalent to the regular representation.
\end{cor}
\begin{proof} All is clear except for the last assertion.
By equation (\ref{eq:factor}) it is enough to know this for the
space of solutions of (\ref{eq:stein}). This is well known, and
follows from the case $\mu\in\af_\C^{\operatorname{reg}}$ by rigidity
of characters of a finite group.
\end{proof}
Notice that the center of the group ring of $W^a$ is
$\C[Q(\Sigma^\vee)]^W$. Thus the central characters of irreducible
representations of $W^a$ correspond to $W$-orbits of points of the
complex algebraic torus
\begin{equation}
T^L=\af_\C^*/P(\Sigma)
\end{equation}
whose exponential map we will denote by $\exp^L$ (i.e.
($\exp^L(\mu)=1$ iff $\mu\in P(\Sigma)$). Then $T^L$
is the dual torus of $A_\C$.
\begin{prop} The monodromy representation of
(\ref{eq:cpx}) has central character
$W\operatorname{exp}^L(\mu)\in W\backslash T^L$ (observe that this
central character is unitary iff $\mu\in\af^*$ is
real). The monodromy representation is irreducible iff
$\exp^L(\mu)$ has trivial isotropy for the action of $W$ on $T^L$.
\end{prop}
\begin{proof}
The monodromy representation clearly has central character
$W\exp^L{\mu}$, and dimension $|W|$ by Proposition \ref{prop:dim}.
The last assertion follows easily from the Mackey induction
procedure.
\end{proof}
The spherical function $\phi^X_\mu|_{A_\C}$ is a special solution of
(\ref{eq:cpx}), which can be characterized by saying that it is a
nonzero $W$-fixed vector in the monodromy representation. By the
above corollary the space of $W$-fixed vectors is
one-dimensional. Looking at (\ref{eq:factor}) we see that
$\phi^X_\mu|_{A_\C}$ is of the form $\d^{-1}\psi$ where $\psi$ is a
$W$-skew solution (say on $\af_\C$) of (\ref{eq:stein}). Then
$\psi$ is divisible by $\prod_{\a\in\Sigma_+}\a$, and thus
$\phi^X_\mu|_{A_\C}$ extends to a holomorphic solution on
$AT_\Omega^2$.

We continue the discussion by considering $\phi^X_\mu|_{A_\C}$ via
$\pi$ as a holomorphic function on $i\af+\Omega$ of the form
$\d^{-1}\psi$ with $\psi$ a $W$-skew solution of (\ref{eq:stein}).
Let $\eta=\omega_i/k_i$ be as before.
We denote by $W_{\eta}$ the isotropy group of ${\eta}$ in $W$,
and by $W_{\eta}^a$ the isotropy group of ${\eta}$ in the affine Weyl
group $W^a$.
\begin{lem}
The natural homomorphism
$W^a\to W$ restricts to an isomorphism from $W^a_{\eta}$ onto the isotropy
group $W_{t(\eta)^2F}\subset W$ of $t(\eta)^2F=\exp(i\pi\eta)F$
for the action of $W$ on $A_\C/F$. In particular
$W_{t(\eta)^2F}=W$ if $\eta\in\partial(\Omega)$ is minuscule.
\end{lem}
\begin{proof}
The injectivity is clear. The isomorphism of complex algebraic tori
$A_\C/F\approx A_\C$ given by $aF\to a^2$ is $W$-equivariant.
Therefore $W_{t(\eta)^2F}$
is equal to the isotropy group of $W_{t(\eta)^4}$ for the action
of $W$ on $A_\C$. Since $A_\C$ has the weight lattice $P(\Sigma)$
as its character lattice the group $W_{t(\eta)^4}$ is generated
by reflections (by a well known result of Steinberg \cite{S2}).
Now suppose that $s=s_\a\in W_{t(\eta)^4}$ is a reflection.
But $s\in W_{t(\eta)^4}\Longleftrightarrow
s(\eta)-\eta\in Q(\Sigma^\vee)\Longleftrightarrow
\a(\eta)=n\in\Z$. Hence the affine reflection
$s_{\a-n}$ satisfies $s_{\a-n}\in W^a_\eta$ and
is mapped to $s$, proving the surjectivity.
\end{proof}
 Let
$\Sigma_{+,{\eta}}^a$ be the set of positive affine roots which
vanish on ${\eta}$, and let
$\Sigma_{+,{\eta}}=\Sigma_{+,{\eta}}^a\cap\Sigma$. Then the Dynkin
diagram of $\Sigma_{+,\eta}^a$ is $D^*-\{\a_\eta\}$ where
$\a_\eta$ is the unique simple root of $\Sigma$ such that
$\a_\eta(\eta)\not=0$, and
$\Sigma_{\eta}\subset\Sigma_{\eta}^a$ is the maximal standard
parabolic subsystem in which we delete $\a_0$ from the set of
simple roots of $\Sigma_{\eta}^a$.

In a small neighborhood of the extremal boundary point $\eta$ of
$\Omega$ we consider the Taylor expansion of $\psi$. The lowest
homogeneous term $h_{\eta,\mu}$ of $\psi$ at $\eta$ is a
$W$-harmonic polynomial which is $W_{\eta}$ skew, and
$\phi^X_\mu|_{A_\C}\circ\pi$ can be uniquely expressed in the form
\begin{equation}
\phi^X_\mu|_{A_\C}\circ\pi=
(\prod_{a\in\Sigma_{\eta,+}^a}a)^{-1}(h_{\eta,\mu}+
\mathrm{higher\ order\ terms\ at\ } \eta)
\end{equation}
Since $h_{\eta,\mu}$ is divisible by
$\prod_{\a\in\Sigma_{\eta,+}}\a$ we see that in this case
the leading exponent $s_{\eta,\mu}^X$ satisfies
\begin{equation}
s_{\eta,\mu}^X\geq  s_\eta^X:=-|\Sigma_{\eta,+}^a-\Sigma_{\eta,+}|
\end{equation}
For $\mu$ in an open, dense subset of $\af_\C^*$ this bound is sharp.
The bound is sharp if $\mu=i\lambda$ with $\lambda\in\af^*$. It is
not so easy to describe the function $\mu\to s_{\eta,\mu}^X$
exactly. We observe that this function is upper semi-continuous.

The above analysis can not be used directly in general, since the
spherical functions do not have a simple factorization formula
like (\ref{eq:factor}) in general. For our later use it
is helpful to describe the above result in terms of the monodromy
representation. By the above we see that $h_{\eta,\mu}$ belongs
to space of $W$-harmonic polynomials which transform by the sign
representation under the action of $W_{\eta}$. This means that
$h_\eta$ is  a $W$-harmonic polynomial in the direct sum of the
isotypical components of the irreducible characters of
$W^a_{\eta}$ which are induced from the sign representation
$\operatorname{det}_{\eta}$ of $W_{\eta}$. Therefore the homogeneous
degree of $h_\eta$ at $\eta$ is at least equal to the harmonic birthday
of the irreducible character (\emph{the leading character})
$\tilde\sigma_\eta\in\operatorname{Irr}(W_{\eta}^a)$ given by
the \emph{truncated induction}
\begin{equation}
\tilde\sigma_\eta\in\operatorname{Irr}(W_{\eta}^a)=
j_{W_{\eta}}^{W_{\eta}^a}(\operatorname{det_{\eta}})
\end{equation}
It follows by truncated induction (see \cite[Section
11.2]{Ca}) that the harmonic birthday of this irreducible
character is equal to $|\Sigma_{\eta,+}|$, and that this
representation has multiplicity $1$ in this degree. Moreover,
the same is true in the space of $W$-harmonic polynomials.
\begin{prop}\label{cpx} Assume we are in the complex case, so
$X$ is a Riemannian symmetric space of type IV with restricted root
system $\Sigma$. Let $\eta\in\partial_e{\Omega}\cap\C$ be an extremal
boundary point of $\Omega$ as before, and define
$-s_\eta^X:=|\Sigma_{\eta,+}^a-\Sigma_{\eta,+}|$
as the number of roots $\a$ in $\Sigma_+$ with $\a(\eta)=1$.
Let $s^X_{\eta,\mu}$ be the leading exponent at $z=t(\eta)^2.x_0$
in the sense of (\ref{eq:leadexp}).
For all $\mu\in\af_\C^*$ the
logarithmic degeneracy of the leading exponent $s^X_{\eta,\mu}$ is $0$.
For all $\mu\in\af_\C^*$ we have $s_{\eta,\mu}^X\geq s_\eta^X$.
For $\mu$ in a dense open set of $\af_\C^*$ containing $i\af^*$
this inequality is an equality (cf.  Theorem \ref{thm:table}).
\end{prop}
\begin{thm}\label{thm:table}
In Table 4 below we have used the numbering of the
extremal boundary points $\eta_j=\omega_j/k_j\in\partial\Omega$
corresponding to the distinguished boundary orbits as in Table 1.
Table 4  displays lower bounds for the leading exponents
of the holomorphically extended elementary spherical functions at
the extremal points $\eta_j$
in the sense of (\ref{eq:leadexp}). In the case where
this lower bound is attained the table displays the
corresponding logarithmic degeneracy and a leading
character (leading characters are explained
in Appendix \ref{app:exp}).

The convention for the root multiplicities is as follows.
We use $m_1\geq 1$ for the root multiplicity of a long root $\a$
(or simply $m$ if $\Sigma$ is reduced and simply laced).
The multiplicity of
half a long root in $\Sigma$ is denoted by $m_{1/2}\geq 0$
(i.e. we view $C_n$ as the special case of $BC_n$ where $m_{1/2}=0$).
The multiplicity of unmultipliable roots $\b\in\Sigma$
(i.e. $2\b\not\in\Sigma$) which are not long roots is denoted
by $m_2\geq 1$ (if such roots exist).

\medskip
{\tiny{
\begin{center}
\begin{tabular}[c]{|c|c|c|c|c|c|}
\hline
\multicolumn{6}{|c|}
{\rule[-2mm]{0mm}{7mm}\textsf{The leading character
$\sigma_\eta\in\operatorname{Irr}{W_\eta^a}$, the leading exponent $s_\eta$
and its degeneracy $d_\eta$}} \\[2pt]
\hline\hline
$\Sigma$& $\eta$&$\Sigma_\eta^a$
&$\sigma_\eta$&$s_\eta$&$d_\eta=1$ iff\\[2pt]\hline\hline
$A_{2r}(r\geq1)$&$\omega_j$($j\leq r$)&
$A_{2r}$&$(2r-j+1,j)$&
$j(1-(2r+2-j)m/2)$&
\\[2pt]\hline
$A_{2r-1}(r\geq2)$&$\omega_j$($j\leq r$)&
$A_{2r-1}$&$(2r-j,j)$&
$j(1-(2r+1-j)m/2)$
&$m=1\&j=r$
\\[2pt]\hline\hline
$B_l(l\geq 3)$&$\omega_1$&
$B_l$&$(l-1,1)$&
$1-(l-1)m_1-m_2$&
\\[2pt]\hline
$B_3$&$\omega_3/2$&
$A_3$&$(3,1)$&
$1-2m_1$&
\\[2pt]\hline
$B_{2r}(r\geq 2)$&$\omega_{2r}/2$&
$D_{2r}$&$(r,r)$&
$r(1-r m_1)$&
\\[2pt]\hline
$B_{2r+1}(r\geq 2)$&$\omega_{2r+1}/2$&
$D_{2r+1}$&$(r+1,r)$&
$r(1-(r+1)m_1)$&
\\[2pt]\hline\hline
$BC_{1}$&$\omega_{1}$&
$A_{1}$&$1^2$&$1-m_1$
&$m_1=1$
\\[2pt]\hline
$BC_{2r}(r\geq1)$&$\omega_{2r}$&
$C_{2r}$&$(r,r)$&$r(1-rm_2-m_1)$&
\\[2pt]\hline
$BC_{2r+1}(r\geq1)$&$\omega_{2r+1}$&
$C_{2r+1}$&$(r,r+1)$&$(r+1)(1-rm_2-m_1)$
&$m_1=1$
\\[2pt]\hline\hline
$D_l(l\geq 4)$ & $\omega_1$&
$D_l$&$((l-1,1),-)$&$2-lm$&$m=1$
\\[2pt]\hline
$D_{2r}(r\geq2)$&$\omega_{2r}$
&$D_{2r}$&$(r,r)^\prime$
&$r(1-rm)$&
\\[2pt]\hline
$D_{2r+1}(r\geq2)$&$\omega_{2r+1}$&$D_{2r+1}$&$(r+1,r)$
&$r(1-(r+1)m)$&
\\[2pt]\hline\hline
$E_6$  & $\omega_1$
&$E_6$&$\phi_{20,2}$&$2-9m$& \\[2pt]\hline\hline
$E_7$  &  $\omega_7$
&$E_7$&$\phi_{21,3}$&$3-15m$&$m=1$\\[2pt]\hline
&$\omega_2/2$&$A_7$
&$(7,1)$&$1-4m$&\\[2pt]\hline\hline
$E_8$&$\omega_1/2$&$D_8$&
$((7,1),-)$&$2-8m$&$m=1$\\[2pt]\hline
&$\omega_2/3$&$A_8$&
$(8,1)$&$1-9/2m$&\\[2pt]\hline\hline
$F_4$&$\omega_4/2$&$B_4$&
$(3,1)$&$1-3m_1-m_2$&\\[2pt]\hline\hline
$G_2$&$\omega_1/3$&$A_2$&
$(2,1)$&$1-3/2m_1$&
\\[2pt]\hline
\hline
\end{tabular}
\end{center}}}
\centerline{\tt Table 4}
\end{thm}
The proof Table 4 and
these facts is given in the Appendix
Section \ref{app:exp}.
In the parameter family of hypergeometric
functions $\phi_{\mu,m}$ (cf. Appendix \ref{app:exp})
with real multiplicities $m$ the indicated lower bounds are
sharp generically in $m$, provided that $m$ satisfies the
inequalities $1\leq m_1\leq m_2$.
The leading character is independent of $m$
in this cone (hence only depends on the geometry of $\Omega$
at the extremal point $\eta$).
\section{Unipotent model for the crown domain}
In this section we give a new geometrical characterization of the
crown by unipotent $G$-orbits. Fot the beginning there
is no restriction on $G$ and we define a connected $G$-subset
of $X_\C$ by

$$\Xi_N=G\exp(i\nf).x_0=GN_\C.x_0\, .$$
For $G=\Sl(2,\R)$ we have shown that $\Xi_N$ is an open
subset, but in the general case this is not clear to us.
The next lemma contains the crucial information.

\begin{lemma} $\Xi\subset\Xi_N$.
\end{lemma}

\begin{proof} Let $ Y\in \pi\Omega/ 2$.
 We recall the complex convexity theorem (\ref{cc})
$$\Im \log a_\C (K\exp(iY).x_0)={\rm co} (W.Y) \, .$$
In particular, there exists a $k\in K$ such that
$\Im \log a_\C (k\exp(iY).x_0)=0$, or, in other words,
$$k\exp(iY) \in N_\C A K_\C = A N_\C K_\C\, .$$
We conclude that $G\exp(iY)\subset GN_\C K_\C$ and
then $G\exp(i\Omega)\subset GN_\C K_\C$, i.e.
$\Xi\subset \Xi_N$.
\end{proof}

Let us define a domain $\Lambda\subset \nf$  by

$$\Lambda=\{ Y\in \nf\mid \exp(iY).x_0\in \Xi\}_0$$
where $\{\cdot \}_0$ stands for the connected component of $\{\cdot \}$
containing $0$. As $\Xi$ is open, it is clear that $\Lambda$ is
open as well.

\begin{lemma} \label{lem=lb}Suppose that $\Omega_c\subset \Omega$ is a compact subset.
Then the set
$$\Lambda_c=\{ Y\in \nf\mid \exp(iY).x_0\in G\exp(i\Omega_c).x_0\}$$
is compact in $\nf$.
\end{lemma}

\begin{proof} First we observe that $\Lambda_c$ is closed
as $G\exp(i\Omega_c).x_0$ is closed in $X_\C$.

\par Let $Y\in \Lambda_c$. Then
$\exp(iY)=g\exp(iZ).x_0$ for some $g\in G$ and $Z\in \Omega_c$.
With  $g=n^{-1} a^{-1} k$ for
$n\in N$, $a\in A$ and $k\in K$ we obtain that
\begin{equation} \label{eq=ii}k\exp(iZ).x_0= an\exp(iY).x_0\, .\end{equation}

We recall that $\Xi\subset A_\C N_\C .x_0$ and that
there is a well defined holomorphic projection
$$\tilde n: A_\C N_\C.x_0\to N_\C\, .$$
Further we note that the map
$$N\times \nf\to N_\C , \ \ (n, Y)\mapsto n\exp(iY)$$
is a diffeomorphism. In particular $N\backslash N_\C \simeq \nf$ under
a homeomorphic map $\psi$.
Consider the continuous map
$$f=\psi\circ \tilde n: A_\C N_\C.x_0\to \nf$$
and note that (\ref{eq=ii}) shows that
$f(k\exp(iZ).x_0)= Y$. Therefore
$$ \Lambda_c\subset f(K\exp(i\Omega_c).x_0)\, .$$
Since  $K\exp(i\Omega_c).x_0$ is compact and $f$ is continuous, the
assertion of the lemma follows.
\end{proof}

We arrive at the main result of this section.

\begin{thm} $\Xi=G\exp(i\Lambda).x_0$.
\end{thm}

\begin{proof} We argue by contradiction. Suppose that the assertion
is false. Then there exists $Z\in    \Omega$ such that $\exp(iZ).x_0\not\in G\exp(i\Lambda).x_0$
and sequences $(Z_n)_{n\in \N}\subset \Omega$, $(g_n)_{n\in \N}\subset G$ and
$(Y_n)_{n\in \N}\subset \Lambda$ such that $Z_n\to Z$ and
$\exp(iZ_n).x_0=g_n \exp(iY_n).x_0$. Let $\Omega_c\subset \Omega$ be a
compact subset of $\Omega$ with $(Z_n)_{n\in \N}\subset \Omega_c $.
Then $\exp(iY_n).x_0\in G\exp(i\Omega_c).x_0$ and we conclude
from Lemma \ref{lem=lb} that $(Y_n)_{n\in \N}$ is a bounded sequence in
$\nf$. W.l.o.g. we may assume that $Y_n\to Y\in \nf$.
As $\Omega_c$ is compact, the set $G\exp(i\Omega_c).x_0$ is closed in $X_\C$
(cf. Remark \ref{rem=tan}) and thus $\exp(iY).x_0\in G\exp(i\Omega_c).x_0
\subset \Xi$. Hence $Y\in \Lambda$. Because $G$ acts properly on $\Xi$ we conclude
that $(g_n)_{n\in \N}$ is bounded and it is no loss of generality
to assume that $\lim_{n\to\infty} g_n =g$. But then $\exp(iZ).x_0=g\exp(iY).x_0\in
G\exp(i\Lambda).x_0$, a contradiction.
\end{proof}

\begin{rem} The determination of the precise shape of $\Lambda$ is a difficult problem,
especially for higher rank groups. Generally one might ask:
Is $\Lambda$ always bounded? Is $\Lambda$ convex?

\par Recall the fact that $G\exp(iZ).x_0=G\exp(iZ').x_0$ for
$Z,Z'\in \pi\Omega/ 2$ means that $W.Z=W.Z'$.
Thus we obtain a well defined map
$$p: \Lambda\to \Omega/ W$$
via $G\exp(iY).x_0=G\exp(i\pi p(Y)/2).x_0$ for $Y\in \Lambda$.
The following would be interesting to know:
What are the fibers of the map $p$? What are  the preimages
of the extreme points? Is there an expressable relationship
between $Y$ and $p(Y)$?
\end{rem}

\subsection{The case of real rank one}

In this subsection we will determine the precise
shape of $\Lambda$ for groups $G$ with real rank one.
We begin with a criterion which will allow us
explicit computations.

\begin{lemma}\label{lem=r1} Suppose that $G$ has real rank one. Then
$$\Lambda=\{ Y\in \nf\mid N\exp(iY).x_0\subset \overline N_\C A_\C.x_0\}_0$$
with $\{ \cdot  \}_0$ denoting the connected component of $\{ \cdot \}$
containing $0$. \end{lemma}

\begin{proof} Set
$$\Lambda_1=\{ Y\in \nf\mid N\exp(iY).x_0\subset \overline N_\C A_\C.x_0\}_0$$
and note that
\begin{equation}\label{eq=l1}
\Lambda_1=\{ Y\in \nf\mid \exp(iY).x_0\subset\bigcap_{n\in N} n \overline N_\C A_\C.x_0\}_0\, .\end{equation}
We recall the fundamental fact on (general) complex crowns that
$$\Xi=\left[\bigcap_{g\in G} g\overline N_\C A_\C. x_0\right]_0$$
with $[\cdot]_0$ denoting the connected component of $[, ]$ containing
$x_0$. We are now going to use the fact that $G$ has real rank one.
In particular $W=\{ 1, w\}=\Z_2$ and the Bruhat decomposition
of $G$ reads $G=N MA\overline N \cup w MA\overline N$. Hence
$$\Xi=\left[\bigcap_{n\in N} n\overline N_\C A_\C. x_0\cap N_\C A_\C.x_0\right]_0\,.$$
As a result (\ref{eq=l1}) translates into
$$\Lambda_1=\{ Y\in \nf \mid \exp(iY)\in \Xi\}_0 =\Lambda\, $$
and the proof is complete.
\end{proof}

We introduce coordinates on $\nf=\gf^\alpha + \gf^{2\alpha}$.
As usual we write $p=\dim \gf^\alpha$ and $q=\dim \gf^{2\alpha}$ and
let $c={1\over 4(p+4q)}$. We endow $\nf$ with the inner
product $\langle Y_1, Y_2\rangle = - \kappa (Y_1, \theta(Y_2))$ where $\kappa$
denotes the Cartan-Killing form of $\gf$.
For $z\in \Xi$ we write in the sequel $a_\C (z)$ for the $A_\C$-part
of $z$ in the Iwasawa decomposition $\overline N_\C A_\C.x_0$.
For $Y\in \gf^\alpha, Z\in \gf^{2\alpha}$ we recall
the formula

\begin{equation}\label{eq=e1} a_\C (\exp(Y+Z).x_0)^\rho =\left[(1+c\|Y\|^2)^2 + 4c \|Z\|^2\right]^{p+2q\over 4}\,.
\end{equation}
The complex linear extension of $\langle\cdot,\cdot \rangle$ to $\nf_\C$ shall
be denoted by the same symbol.
We obtain the following criterion for
$\Lambda$.

\begin{lemma}\label{lem=exc1} If $G$ has real rank one, then
\begin{align*}\Lambda  =&\{ (Y,Z)\in \gf^\alpha \oplus\gf^{2\alpha}\mid
(\forall (Y',Z')\in \nf) \quad (1+c\langle Y'+iY, Y'+iY\rangle )^2 \\
&\quad  +4c \langle Z'+iZ +i{1/2} [Y',Y],
 Z'+iZ +i{1/2} [Y',Y]\rangle \neq 0\}_0\, .\end{align*}
\end{lemma}

\begin{proof} A standard argument (see \cite{KSI}, Lemma 1.6) combined
with
Lemma \ref{lem=r1} yields that
$$\Lambda=\{ Y\in \nf\mid (\forall n\in N)\  a_\C (n\exp(iY)) \hbox{ is defined}\}_0\, .
$$
Now for $n=\exp(Y'+Z')\in N$ with $(Y',Z')\in\nf$ and $(Y,Z)\in \Lambda$ one has
$$n\exp(i(Y+Z))=\exp(Y'+iY + Z'+i(Z +1/2 [Y',Y]))$$
and the assertion follows in view of the explicit formula (\ref{eq=e1}).
\end{proof}

We use the criterion in Lemma \ref{lem=exc1} to determine
$\Lambda$ explicitly. However, this is not so easy
as it looks in the beginning.
We shall begin with two important special cases and start
with the Lorentz groups $G={\rm SO}_e(1,p+1)$ where $q=0$.

\begin{lemma}\label{lem=lo} Assume that $G$ is locally ${\rm SO}_e(1,p+1)$. Then
$c={1\over 4p}$ and
$$\Lambda=\{ Y\in \nf=\R^p\mid  c
\|Y\|^2 < 1\}\, .$$
\end{lemma}

\begin{proof} In view of the previous lemma we have to look at the
connected component of those $Y\in \nf$ such that
$$  1+c\langle Y'+iY, Y'+iY\rangle= 1 +c ( \|Y'\|^2 -\|Y\|^2 - 2i \langle Y, Y'\rangle) \neq 0 $$
for all $Y'\in \nf$. The assertion follows.
\end{proof}

Next we consider the case of the group $G={\rm SU}(2,1)$.
Here $p=2$ and $q=1$ and so $c={1\over 24}$. Define matrix elements

$$X_\alpha=\begin{pmatrix} 0 & 1 & 0\\ -1 & 0 & 1\\ 0 & 1 & 0\end{pmatrix},
\quad Y_\alpha=\begin{pmatrix} 0 & i & 0\\ i & 0 & -i\\ 0 & i & 0\end{pmatrix}
, \quad X_{2\alpha}={1\over 2}
\begin{pmatrix} i & 0 & -i\\ 0 & 0 & 0\\ i & 0 & -i\end{pmatrix}\,$$
and note that
$$\gf^\alpha= \R X_\alpha \oplus \R Y_\alpha\quad \hbox{and}\qquad \gf^{2\alpha}=\R X_{2\alpha}\, .$$
We record the commutator relation $[X_\alpha, Y_\alpha]=4 X_{2\alpha}$ and
the orthogonality relation $\langle X_\alpha, Y_\alpha\rangle =0$. Finally
we need that $\|X_\alpha\|^2=\|Y_\alpha\|^2={1\over c}$ and $\|X_{2\alpha}\|^2 ={1\over 4c}$.

\begin{lemma}\label{lem=su} For $G$ locally ${\rm SU}(2,1)$ one has
\begin{align*}\Lambda & =\{ xX_\alpha+ yY_\alpha +zX_{2\alpha}\in \nf \mid  2(x^2+y^2) + |z|< 1\}\\
&=\{(Y,Z)\in \nf\mid 2c \|Y\|^2 + 2\sqrt{c} \|Z\|<1\}\, .
\end{align*}
\end{lemma}

\begin{proof} We want to determine those $Y\in \nf$ which belong to $\Lambda$.
By the $M$-invariance of $\Lambda$ we may restrict our attention to elements
of the form $Y=xX_\alpha+ zX_{2\alpha}\in \Lambda$.
We have to find the connected component of those $x,z$ such that
\begin{align*}& (1+c\langle(u+ix)X_\alpha + vY_\alpha, (u+ix)X_\alpha + vY_\alpha\rangle)^2\\
& \quad +4c \langle (w+iz) X_{2\alpha} + {1\over 2} ixv[Y_\alpha, X_\alpha],
(w+iz) X_{2\alpha} + {1\over 2} ixv[Y_\alpha, X_\alpha]\rangle =0\end{align*}
 has no solution for $u,v,w\in \R$. We employ the precedingly  collected material
on commutators, orthogonality and norms and obtain the equivalent version

$$ (1+(u+ix)^2 + v^2)^2 + (w+i(z+ 2xv))^2 =0$$
for $u,v,w\in \R$. However, this is equivalent to
$$ 1+(u+ix)^2 + v^2 = \pm i(w+i(z+2xv))= \pm i w  \mp (z+2xv)\,. $$
Comparing real and imaginary part yields the system of equations

$$ (1 -x^2 \pm 2z) + u^2 +v^2 = \mp 2xv$$
$$ 2ux =\pm w\, .$$
We can always choose $w$ so that the second equation is satisfied.
Hence we look for $x,z$ such that
the quadratic equation in $v$

$$ v^2 -2xv + (1 -x^2  \pm  z + u^2)=0\,.$$
has no solution for all $u$. Clearly we can take $u=0$ and
assume $\pm z= -|z|$. We are left with analyzing the
discriminant
$$ 4x^2 - 4( 1-x^2 -|z|)<0\,.$$
This inequality translates into $2x^2+ |z|<1$ and
concludes the proof of the lemma.
\end{proof}

\begin{rem} Consider the domain
$$\tilde\Lambda=\{Y\in \nf\mid \exp(iY).x_0\in \overline N_\C A_\C.x_0\}_0\,.$$
It is clear that $\Lambda\subset\tilde\Lambda$ and it is easy
to determine $\tilde \Lambda$ explicitly:
\begin{align*} \tilde\Lambda& =\{ (Y,Z)\in \nf\mid (1 - c\|Y\|^2)^2 -4c\|Z\|^2 >0\}\\
&=\{ (Y,Z)\in \nf\mid c\|Y\|^2 +2\sqrt{c}\|Z\| <1\}\, .
\end{align*}
Now for $q=0$ we have seen that $\tilde\Lambda=\Lambda$. However, as
our previous analysis of the ${\rm SU}(2,1)$-case shows, one has
$\Lambda\neq \tilde\Lambda$ in general.
\end{rem}

Before we come to the determination of  $\Lambda$ for all rank one
cases some remarks concerning the nature of the constant
$c={1\over 4(p +4q)}$ are appropriate.

\begin{rem}
\label{rem=norm}Let $\gf$ be of real rank one. Let $E\in
\gf^{\alpha}$, resp. $E\in \gf^{2\alpha}$ and set $F=\theta(E)$,
$H=[E,F]$. If we assume that $\{H, E, F\}$ is an $\sl(2)$-triple,
i.e. $[H,E]=2E$ and $[H,F]=-2F$, then  elementary
$\sl(2)$-representation theory gives $\|E\|^2={1\over c}$ if $E\in
\gf^\alpha$ and $\|E\|^2={1\over 4c}$ if $E\in \gf^{2\alpha}$.
\end{rem}

\begin{thm}\label{th=r1} Let $G$ be a simple Lie group of real rank one.
\begin{enumerate} \item If $q=0$, then
$$\Lambda=\left\{ Y\in \nf\mid  c
\|Y\|^2 < 1\right\}\, .$$
\item If $q>0$, then
$$\Lambda=\{ (Y,Z)\in \nf\mid  2c\|Y\|^2 + 2\sqrt{c}\|Z\|< 1\} \, .$$
\end{enumerate}
\end{thm}

\begin{proof} (i) Lemma \ref{lem=lo}.
\par \noindent (ii) Set
$$\tilde \Lambda=\{ (Y,Z)\in \nf\mid  2c\|Y\|^2 + 2\sqrt{c}\|Z\|<
1\} \, .$$ We first show that $\tilde \Lambda \subset \Lambda$.
Let $(Y,Z)\in \tilde \Lambda$, $Y\neq 0$ and $Z\neq 0$. Consider
the Lie algebra $\gf_0$ generated by $Y, \theta(Y), Z, \theta
(Z)$. Standard structure theory says that  $\gf_0\simeq \su(2,1)$, see
\cite{Hel}, Ch. IX, \S 3.
Choose $E_\alpha\in \R Y $ and $E_{2\alpha}\in \R Z$ such that $\{
[E_\alpha, \theta (E_\alpha)], E_\alpha, \theta(E_\alpha)\} $ as
well as  $\{ [E_{2\alpha}, \theta (E_{2\alpha})], E_{2\alpha},
\theta(E_{2\alpha})\} $ are $\sl(2)$-triples. Let $y,z\in \R$ such
that $Y=yE_\alpha$ and $Z=zE_{2\alpha}$. In view of the previous
Remark \ref{rem=norm} the condition $ 2c\|Y\|^2 + 2\sqrt{c}\|Z\|<
1$ is equivalent to $2y^2  + |z|< 1$. But this is just the
condition for $\exp(i(Y+Z)).x_0$ to be contained in the crown
domain $\Xi_0$ for the group $G_0=\langle \exp \gf_0\rangle < G$
(see Lemma \ref{lem=su}). Now, with the obvious notation, we have
$\Omega=\Omega_0$ and so $\Xi_0\subset \Xi$. This concludes the
proof of $\tilde \Lambda \subset \Lambda$.
\par It remains to verify that $\Lambda\subset \tilde \Lambda$. For that it is
sufficient to show the following: if $(Y,Z)\in \partial \tilde \Lambda$, then
$(Y,Z)\not\in \Lambda$. Let also $(Y,Z)\in \partial \tilde \Lambda$ with
$(Y,Z)\in \Lambda$
and let
$\gf_0$ as before. Note that $\exp(i(Y+Z)).x_0\in \partial\Xi_0$.
As $\Xi_0\subset \Xi$ is closed and $\exp(i(Y+Z)).x_0\in \Xi$,
there would exist $G_0$-domain $\Xi_0'\subset X_{0,\C}$, properly containing
$\Xi_0$, and on which $G_0$ acts properly. But this contradicts
Theorem \ref{th=p}.
\end{proof}

\subsection{The case where $\Xi$ is a Hermitian symmetric space}
It can happen that $\Xi$ allows additional symmetries, i.e. the
group of holomorphic automorphisms is strictly larger then $G$.
For example, when $X=G/K$ is a Hermitian symmetric space, then
$\Xi$ is biholomorphic to $X\times \overline X$ where $\overline
X$ denotes $X$ endowed with the opposite complex structure (see
\cite{KSII}, Th. 7.7). In this example $\Xi=X\times \overline X$
is again a Hermitian symmetric space and $\Aut(\Xi)= G\times G$ is
twice the size of $G=\Aut(X)$. In \cite{KSII}, Th. 7.8, one  can
find a classification of all those cases where $\Xi$ is a
Hermitian symmetric space for a larger group $S$. For all these
cases it turns out that there is an interesting subset
$\Lambda^+\subset \Lambda$ such that $\Xi=G\exp(i\Lambda^+).x_0$,
and moreover, it is possible to give a precise relation between
unipotent $G$-orbits through $\exp(i\Lambda^+).x_0$ and   the
elliptic $G$-orbits through $\exp(i\Omega).x_0$. As the case of
the symplectic group is of special interest, in particular for
later applications to automorphic forms, and as it is always good
to have a illustrating example, we shall begin with a discussion
for this group.

\subsubsection{The symplectic group}

In this section $G={\rm Sp}(n,\R)$ for $n\geq 1$.
Let us denote by ${\rm Sym}(n,\R)$, resp.  $M_+(n,\R)$, the symmetric, resp. strictly upper triangular,
matrices in $M(n,\R)$. Our choices
of $\af$ and $\nf$ shall be

$$\af=\left\{ {\rm diag}(t_1, \ldots, t_n, -t_1, \ldots, -t_n)\mid t_i\in \R\right\}$$
and
$$\nf=\left\{ \begin{pmatrix} Y & Z\\ 0 & -Y^T\end{pmatrix}\mid Z\in {\rm Sym}(n,\R), Y\in M_+(n,\R)\right\}\, .$$
Of special interest is an  abelian subalgebra of $\nf$
$$\nf^+=\left \{ \begin{pmatrix} 0& Z\\ 0 & 0\end{pmatrix}\mid Z\in {\rm Sym}(n,\R)\right\}, .$$
We recall that the maximal compact subgroup $K<G$ is isomorphic to $U(n)$ and that
$X=G/K$ admits a natural realization as a Siegel upper halfplane:

$$X={\rm  Sym}(n,\R)+ i {\rm  Sym}^+(n,\R)\subset {\rm Sym}(n,\C)$$
where ${\rm  Sym}^+(n,\R)$ denotes the positive definite symmetric matrices.
The action of $G$ on $X$ is given by generalized fractional
transformations: if $g=\begin{pmatrix} A& B\\ C & D\end{pmatrix}\in G$ and $Z\in X$, then
$$g(Z)=(AZ +B) (CZ+D)^{-1}\, .$$
Notice that the base point $x_0$ with stabilizer $K$ becomes $x_0=i I_n$ with $I_n$ the
identity matrix.
The natural realization of $\overline X$ is the lower half plane
$$\overline X =  {\rm  Sym}(n,\R) - i {\rm  Sym}^+(n,\R)\, .$$
In the sequel we view $X$ inside of $X\times \overline X$
as a totally real submanifold via the embedding
$$X\hookrightarrow X\times \overline X, \ \ Z\mapsto (Z, \overline Z)$$
where $\overline Z$ denotes the complex conjugation in ${\rm Sym}(n,\C)$ with
respect to the real form ${\rm Sym}(n,\R)$.
As we remarked earlier, $\Xi$ is naturally biholomorphic to $X\times \overline X$.
Let us now consider a domain in $\nf^+$

$$\Lambda^+=\{ Y\in \nf^+\mid \exp(iY).x_0\in \Xi\}_0\, .$$

\begin{thm} \label{th=sp}For $G={\rm Sp}(n,\R)$ the following assertions hold:
\begin{enumerate}
\item If $\|\cdot \|$ denotes the operator norm on $M(n,\R)$, then
$$\Lambda^+=\left\{ \begin{pmatrix} 0& Z\\ 0 & 0\end{pmatrix}\in \nf^+ \mid \|Z\|< 1\right\}\, . $$
\item With $\Lambda^{++}=\Lambda\cap {\rm  diag}(n,\R)$ one has
\begin{enumerate} \item $\Lambda^{++}=\{ {\rm diag}(t_1, \ldots, t_n)\in \nf^+\mid |t_i|<1\}$
\item$\Lambda^+=\Ad K_0(\Lambda^{++})$ with $K_0={\rm SO}(n,\R)<K $.
\end{enumerate}
\item $\Xi=G\exp(i\Lambda^+).x_0=G\exp(i\Lambda^{++}).x_0$.
\end{enumerate}
\end{thm}

\begin{proof} (i) Let $Z\in {\rm Sym}(n,\R)$ and
$\tilde Z=  \begin{pmatrix} 0& Z\\ 0 & 0\end{pmatrix}$ the corresponding element in
$\nf^+$. Then
$$\exp(i\tilde Z)=\begin{pmatrix} I_n& iZ\\ 0 & I_n\end{pmatrix}$$
and accordingly
$$\exp(i\tilde Z).x_0=\exp(i\tilde Z) (iI_n, -iI_n)=
(i(I_n+Z), -i(I_n-Z))\, .$$
Therefore $\exp(i\tilde Z).x_0\in  X\times \overline X$ if and only if
$I_n+Z \in {\rm Sym}^+(n,\R)$ and $I_n-Z \in {\rm Sym}^+(n,\R)$.
Clearly, this is equivalent to $\|Z\|<1$ and the proof of (i)
is finished.
\par\noindent (ii) This is immediate from (i).
\par\noindent (iii) It is enough to show that $\Xi=G\exp(i\Lambda^{++}).x_0$ and
for that it suffices to verify that $\exp(i\Omega).x_0\subset G\exp(i\Lambda^{++}).x_0$.
Now, as $\Sigma$ is of type $C_n$, the domain $\Omega$ is a cube
$$\Omega =\{{\rm diag}(t_1, \ldots, t_n, -t_1, \ldots, -t_n)\mid |t_i|<{\pi \over 4}\}\, .$$
Let us write  $E_{ij}$ for the elementary matrices in $M(2n,\R)$. Then for
each $1\leq j\leq n$ we define an $\sl(2)$-subalgebra $\gf_j$ by
$$\gf_j={\rm span}_\R\{ E_{jj}-E_{j+n, j+n}, E_{j, j+n}, E_{n+j, j}\}\, .$$
We note that the $\gf_j$ pairwise commute and so $\gf_0=\gf_1\oplus\ldots \oplus
\gf_n$ is a subalgebra of $\gf$ which contains $\af$. Now we are in the situation
to use $\sl(2)$-reduction and the assertion  becomes a consequence of
Lemma \ref{lem=or}.
\end{proof}

For later reference we wish to make the last part of the above theorem more precise.
For ${\bf z}=(z_1, \ldots, z_n) \in \C^n$ let us define a matrix
in $N_\C^+$

$$n_{\bf z}=\begin{pmatrix} I_n & {\rm diag}({\bf z})\\  0& I_n\end{pmatrix}\, .$$
Moreover if ${\bf z}\in (\C^*)^n$, then we set
$$ a_{{\bf z}}={\rm diag}(z_1,\ldots, z_n, z_1^{-1}, \ldots, z_n^{-1})\in A_\C \, .$$
In the course of the proof of Theorem \ref{th=sp} (iii) we have
shown the following:

\begin{lemma}\label{lem=orn} Let ${\bf t}=(t_1,\ldots, t_n)\in \R^n$ with
$|t_i|<{\pi\over 4}$. Set ${\bf e}^{i\bf t}=(e^{it_1}, \ldots, e^{it_n})$ and
${\bf sin}(2 {\bf t})=(\sin 2t_1, \ldots, \sin 2t_n)$. Then
$$G n_{{\bf  sin}(2{\bf t})}.x_0= G a_{{\bf e}^{i{\bf t}}}.x_0\, .$$
\end{lemma}

\subsubsection{The general case of Hermitian $\Xi$}\label{ssh}

In this subsection we will assume that $\Xi$ is a Hermitian
symmetric space for an overgroup $S\supset G$. From a technical
point of view it is however better to work with
an alternative characterization, namely  (cf.  \cite{KSII}, Th. 7.8)
$\Sigma$ is of type $C_n$ or $BC_n$ for $n\geq 2$ or
$\gf=\so(1,k)$ with $k\geq 2$ for the rank one cases.

\par If $\Sigma$ is of type $C_n$ or $BC_n$, then

$$\Sigma=\{ \pm \gamma_i\pm \gamma_j\mid 1\leq i,j\leq n\}\backslash \{0\}
\cup\{\pm {1\over 2}\gamma_i: 1\leq i\leq n\}$$
with the second set on the right  to be considered not present in the
$C_n$-case.
As a positive system of $\Sigma$ we choose
$$\Sigma^+=\{ \gamma_i\pm \gamma_j\mid 1\leq i\leq j\leq n\}\backslash \{0\}
\cup\{{1\over 2}\gamma_i: 1\leq i\leq n\}\, .$$
Further we consider the $A_{n-1}$-subsystem
$$\Sigma_0=\{ \pm \gamma_i \mp \gamma_j\mid 1\leq i\neq  j\leq n\}\, $$
and set
$$\Sigma^{++}=\Sigma^+\cap (\Sigma^+\backslash\Sigma_0) \quad \hbox{and}\quad
\Sigma^{--}=-\Sigma^{++}\, .$$
Next we define subalgebras of $\gf$ by

$$\nf^+=\bigoplus_{\alpha\in \Sigma^{++}} \gf^\alpha, \quad
\nf^-=\bigoplus_{\alpha\in \Sigma^{--}} \gf^\alpha \quad\hbox{and}
\quad \gf(0)=\af \oplus\mathfrak{m} \oplus \bigoplus_{\alpha\in \Sigma_0} \gf^\alpha\, .$$
We note that $\nf^+$ is a subalgebra of $\nf$ and that
$$\gf=\nf^- \oplus \gf(0) \oplus \nf^+$$
is a direct decomposition with $[\gf(0),\nf^{\pm}]\subset \nf^\pm$.
Define elements  $T_j\in \af$ by the requirement $\gamma_i(T_j)=\delta_{ij}$ and note
that
$$\Omega=\bigoplus_{j=1}^n \left ]-{\pi\over 4}, {\pi\over 4}\right [ T_j \, .$$

Now for an element $Y_j\in \gf^{2\gamma_j}$ we find $E_j\in \R Y_j$, unique up to
sign,  such that
$\{ T_j, E_j, \theta(E_j)\}$ form an $\sl(2)$-triple.  Define
$y_j\in \R$ by  $Y_j= y_j E_j$. With that
we can define an open ball in $\bigoplus_{j=1}^n \gf^{2\gamma_j}$ by

$$\Lambda^{++}=\left\{ Y=\sum_{j=1}^n Y_j \in \bigoplus_{j=1}^n \gf^{2\gamma_j}\mid |y_j|< 1\right \} $$
Further write $\kf(0)=\gf(0)\cap \kf$ and set $K(0)=\exp \kf(0)$.
Finally we define a subset of $\nf^+$ by
$$\Lambda^+=\Ad K_0 (\Lambda^{++})\, .$$

\begin{rem} Define an abelian subspace of $\nf^+$ by
$\nf^{++}=\bigoplus_{\alpha\in \Sigma^{++}\cap C_n} \gf^\alpha$.
Then  $\Lambda^+$ is a bounded convex domain in
$\nf^{++}$.
\end{rem}

\begin{thm}\label{ts} If $\Xi$ is a Hermitian symmetric space for an overgroup
$S\supset G$, then the following assertions hold.
\begin{enumerate}
\item If $T =\sum_{j=1}^n t_j T_j\in \Omega$, and $\{ T_j, E_j, \theta(E_j)\}$ is
any $\sl(2)$-triple with $E_j\in \gf^{2\gamma_j}$, then
$$G\exp\left(i\sum_{j=1}^n \sin (2t_j)E_j\right).x_0= G \exp\left(i\sum_{j=1}^n t_j T_j\right).x_0\, .$$
\item $\Xi =G\exp(i\Lambda^+).x_0=G\exp(i\Lambda^{++}).x_0$.
\end{enumerate}
\end{thm}

\begin{proof} (i) We define subalgebras of $\gf$ which are isomorphic
to $\sl(2,\R)$ by
$\gf_j=\R T_j \oplus \R E_j \oplus \R \theta(E_j)$. The $\gf_j$'s commute in $\gf$
and so $\gf_0=\gf_1\oplus\ldots\oplus \gf_n$ defines a subalgebra of $\gf$. In view of
Lemma \ref{lem=or}, the assertion now
follows by $\sl(2)$-reduction.
\par\noindent (ii) This is a consequence of (i).
\end{proof}

\begin{rem} If $\gf$ is Hermitian and of tube type, then
$\Lambda^+$ is a bounded open convex set in $\nf^+$ and
$$\Lambda^+=\{ Z\in \nf^+\mid  \exp(iZ).x_0\in \Xi\}_0\, .$$
This can be proved as in the ${\rm Sp}(n)$-case by employing the
machinery of Jordan algebras.
\end{rem}

\subsection{Some partial results for the special linear groups}\label{sssl}

In this subsection we exclusively deal with $G=\Sl(n,\R)$.
To determine the exact shape of $\Lambda$ for $n\geq 3$ seems
to be very challenging; already the case of $n=3$ appears to be
very intricate. Instead we will exhibit a fairly
large cube-domain inside of $\Lambda$; further
we will estimate the corresponding
hyperbolic parameterization.

\par In order to perform reasonably efficient computations
we use the matrix model for $X_\C$. Let us denote by
${\rm Sym} (n,\C)_{{\rm det}=1}$ the affine variety of complex
symmetric matrices with unit determinant.
The map
$$X_\C=\Sl(n,\C)/{\rm SO}(n,\C)\to {\rm Sym} (n,\C)_{{\rm det}=1}, \ \ gK_\C \mapsto
gg^t$$
is an isomorphism. Within this model for $X_\C$, the Riemmannian
symmetric space $X$ identifies with ${\rm Sym} (n,\R)_{{\rm det}=1}^+$,
the determinant one section in the cone of
positive definite symmetric matrices.
Now the crown domain $\Xi$ contains the determinant
one cut $\Xi_0$ of the tube domain, i.e.

$$\Xi_0=\{ Z\in X_\C \mid  \Re Z \gg 0\}\, .$$

\par As usual we write $E_{ij}=(\delta_{ki}\delta_{ij})_{lj}$ for the elementary matrices.
We choose $N$ to be the group of unipotent upper triangular
matrices and consider the mapping
$$m: \C^{n-1}\to N_\C , \ \ (z_1,\ldots, z_{n-1})
\mapsto \exp(z_1 E_{12})\cdot \ldots \cdot
\exp(z_{n-1}E_{n-1\, n})\, .$$
In matrix notation $m$ is given by

$$m(z_1,\ldots, z_{n-1})=\begin{pmatrix} 1 & z_1 &        &   &\\
                                           & 1   & z_2    &   &  \\
                                           &     & \ddots &\ddots   &\\
                                           &     &        & 1 & z_{n-1} \\
                                           &     &        &   & 1 \end{pmatrix}\, .$$

We define a subset of $N_\C$ by

$$\NN^+=m\left (i\prod_{j=1}^{n-1} (-1,1)\right)$$
and claim:

\begin{prop} $\NN^+\cdot x_0\subset \Xi_0$. In particular $G\NN^+\cdot x_0\subset \Xi_0\subset \Xi$.
\end{prop}

\begin{proof} This is an elementary matrix computation.
Let $(t_1,\ldots, t_{n-1})\in \R^{n-1}$, $|t_i|<1$ and set
$$n=\begin{pmatrix} 1 & it_1 &        &   &\\
                                           & 1   & it_2    &   &  \\
                                           &     & \ddots &\ddots   &\\
                                           &     &        & 1 & it_{n-1} \\
                                           &     &        &   & 1 \end{pmatrix}\, .$$

One has to verify that $\Re (nn^t)\gg 0$. A straightforward calculation
yields

$$nn^t=\begin{pmatrix} 1 & it_1 &        &   &\\
                       it_1                    & 1-t_1^2   & it_2    &   &  \\
                                           & it_2    & \ddots &\ddots   &\\
                                           &     & \ddots       & 1-t_{n-2}^2 & it_{n-1} \\
                                           &     &        & it_{n-1}  & 1-t_{n-1}^2 \end{pmatrix}\, .$$
Therefore

$$\Re (nn^t)=\begin{pmatrix} 1 &           &        &\\
                               & 1-t_1^2   &        & \\
                               &           &\ddots  & \\
                               &           &        &  1-t_{n-1}^2\end{pmatrix}\gg 0 $$

\end{proof}

Next we discuss hyperbolic parameterization for elements in $\NN^+$. Here
our results are somewhat partial but perhaps still interesting.
For what follows we are indebted to Philip Foth.
For $t\in \R$ with $|t|<1$ we consider the element
$$z(t)= m(i(t, \ldots, t))= \begin{pmatrix} 1 & it &        &   &\\
                                           & 1   & it     &   &  \\
                                           &     & \ddots &\ddots   &\\
                                           &     &        & 1 & it  \\
                                           &     &        &   & 1 \end{pmatrix}\, .$$

We wish to estimate the element $a(t)\in \exp(i\pi \Omega/2)$ for which
$$Gz(t)\cdot x_0=Ga(t)\cdot x_0$$
holds. The result is as follows.

\begin{prop}\label{prost} Let $G=\Sl(n,\R)$. Fix $t\in\R$, $|t|<1$ and set $z(t)=m(i(t, \ldots, t))$.
Then $Gz(t)\cdot x_0=Ga(t)\cdot x_0$ with
$$a(t)=\diag (e^{i\phi_1(t)}, \ldots,
e^{i\phi_n(t)}),\qquad  \diag(\phi_1(t), \ldots, \phi_n(t))\in \pi\Omega/ 2$$ and
\begin{equation}\label{es1}
|\phi_j(t)|\leq \left|{1\over 2} \tan^{-1} \left({2t\over 1-t^2}\right)\right| \qquad (1\leq j\leq n)\, .
\end{equation}
\end{prop}

\begin{proof} We proceed indirectly and use the complex convexity theorem
(\ref{cc}). For $k\in K$ we have to show that
the components of $\Im \log a_\C (kz(t))$ satisfy the
estimate (\ref{es1}). To compute $a_\C (kz(t))$ we write
the corresponding matrix identity out:

$$
\underbrace{\begin{pmatrix} 1 & * & \ldots       & \ldots  &*\\
                                           & 1   & *  &\ldots & * \\
                                           &     & \ddots &\ddots   &\vdots\\
                                           &     &        & 1 & *  \\
                                          &     &        &   & 1 \end{pmatrix}}_{\in N_\C}
\cdot
\underbrace{\begin{pmatrix} * &\ldots          &*\\
                * &\ldots         &*\\
                \vdots &\vdots    &\vdots\\
                 * &\ldots  &*\\
                k_1 & \ldots             & k_n\end{pmatrix}}_{\in K}
\cdot \underbrace{\begin{pmatrix} 1 & it &        &   &\\
                                           & 1   & it     &   &  \\
                                           &     & \ddots &\ddots   &\\
                                           &     &        & 1 & it  \\
                                           &     &        &   & 1 \end{pmatrix}}_{=z(t)} =$$

$$=
\underbrace{\begin{pmatrix} a_{\C,1}(t)&  &               &        &   \\
                                          & a_{\C,2}(t)   &        &    \\
                                          &               & \ddots &     \\
                                          &               &        &  a_{\C,n}(t)
\end{pmatrix}}_{\in A\exp(i\pi\Omega/2)}
\cdot
\underbrace{\begin{pmatrix} *     & \ldots       & *          \\
                            *     & \ldots       & *          \\
                         \vdots   & \vdots       & \vdots     \\
                           k_1' &  \ldots        &  k_n'
\end{pmatrix}}_{\in K_\C}\, .$$

We match the bottom rows and arrive at:

$$(k_1, it k_1 +k_2, itk_2+k_3, \ldots, it k_{n_1} +k_n)= a_{\C,n}(t)(k_1', \ldots, k_n')\, .$$
We square the entries and sum them up:
$$1-t^2(k_1^2+\ldots+ k_{n-1}^2) +2it(k_1 k_2+\ldots+ k_{n-1}k_n)= a_{\C,n}(t)^2\, .$$
The result is

$$|\phi_n(t)| = \left| {1\over 2} \tan^{-1}\left ({ 2t(k_1k_2+\ldots+ k_{n-1} k_n)\over
 1-t^2(k_1^2+\ldots +k_{n-1}^2)}\right)\right|\, .$$
Finally we use the estimates
$k_1^2+\ldots +k_{n-1}^2\leq 1$ and
$k_1k_2+\ldots+ k_{n-1} k_n\leq 1$ and obtain that

$$|\phi_n(t)| = \left| {1\over 2} \tan^{-1}\left({ 2t\over
 1-t^2}\right)\right|\, .$$
This proves (\ref{es1}) for the last entry. The general
case follows by Weyl group invariance.
\end {proof}

\begin{rem} One can show that $z(t)\cdot x_0\in \Xi$ precisely for
$|t|<1$.
\end{rem}

\section{Exponential decay of Maa\ss{} cusp forms I: The example of $G=\Sl(2,\R)$}

It is a result obtained by Langlands that cuspidal automorphic forms
are of rapid decay. But actually more is true and the decay is of exponential
type. The purpose of this section is to give an introduction to this circle
of problems with a  solid discussion of the case of $G=\Sl(2,\R)$.
We will restrict our attention to Maa\ss{} cusp forms and to the
modular group $\Gamma=\Sl(2,\Z)$ in order to keep the exposition basic.
It is possible to verify the exponential decay by our explicit
knowledge of the Whittaker functions in this case. This
will be presented first. In general however, concrete
knowledge of the Whittaker functions is not available and
an alternative approach is needed. It was Joseph
Bernstein who came up with the idea to use
analytic continuation to obtain exponential decay.
We shall present his ideas in the geometric framework which we
developed in the preceding section.

\subsection{Concrete approach}

For the rest of this section we let $G=\Sl(2,\R)$ and keep
our choices including notation from Section \ref{sec=nf}.
In the sequel we will identify $X=G/K$ with the upper
half plane, i.e. $X=\{ z\in \C \mid \Im z >0\}$ and
$G$ acting by fractional linear transformations:

$$g(z)={az+ b\over cz+d} \qquad \hbox{for}\quad g=\begin{pmatrix} a & b\\ c & d\end{pmatrix}\in G, \ z\in X\, .$$
In these coordinates the base point $x_0$ is the imaginary unit
$x_0=i$ and the Iwasawa decomposition states that
the map

$$N\times A \to X, \ \ (n_x, a_y)\mapsto n_xa_y(i)=x+iy\, , $$
where

$$n_x=\begin{pmatrix} 1 & x \\ 0 & 1\end{pmatrix} \in N \quad \hbox{and}\quad
a_y=\begin{pmatrix} \sqrt{y} & 0 \\ 0 & {1\over \sqrt{y}}\end{pmatrix} \in A $$
with $x\in \R, y>0$,
is a diffeomorphism.
The Laplace-Beltrami operator of $X$ is given by $\Delta=-y^2(\partial_x^2 + \partial_y^2)$
and we note that $\DD(X)=\C[\Delta]$.
We make a simpler choice for a lattice $\Gamma<G$, namely
$\Gamma=\Sl(2,\Z)$ the modular group.
Then by a {\it Maa\ss{} automorphic form} we understand an analytic function
$\phi: X\to\C$ such that
\begin{itemize}
\item $\phi$ is $\Gamma$-invariant
\item $\phi$ is an eigenfunction for $\DD(X)$, i.e. there exists $\lambda\in \C$ such that
$\Delta \phi =\lambda(1-\lambda)\phi$.
\item $\phi$ is of moderate growth, i.e. there exists $\alpha\in\R$ such that
$$|\phi(x+iy)|\ll  y^\alpha \qquad (y>1)\, .$$
\end{itemize}

Moreover, a Maa\ss{} automorphic form is called a {\it cusp form} if
$$\int_{N\cap \Gamma\backslash N  } \phi(nz) \ dn=0 \qquad \hbox {for all $z\in X$}\, .$$
Note that $N\cap \Gamma=\begin{pmatrix}  1 & \Z\\ 0 & 1\end{pmatrix}$ so that
$N\cap \Gamma\backslash N \simeq \Z\backslash  \R$ is a circle.
In our special case the results of
Langlands reads as follows.

\begin{thm} Maa\ss{} cusp forms $\phi$ are of rapid decay, i.e.
$$|\phi(x+iy)| \ll y ^\alpha \qquad (y>1)$$
for any $\alpha\in \R$.
\end{thm}

However, more is true and we can state:

\begin{thm}\label{th=mcf}  Maa\ss{} cusp forms $\phi$ are of exponential decay, i.e.
there is a constant $C>0$ such that
$$|\phi(x+iy)| \leq C  e^{-2\pi y}\qquad (y>1) \, .$$
\end{thm}

Before we prove this theorem, we will recall the Whittaker
expansion of a Maa\ss{} cusp form: If
$\phi$ is a Maa\ss{} cusp form with $\Delta \phi=\lambda(1-\lambda)\phi$, then

\begin{equation} \label{eq=FB}
\phi(x+iy)=\sum_{n\in \Z^\times} a_n \sqrt{y} K_\nu(2\pi |n| y) e^{2\pi i nx}\,
\end{equation}
where $K_\nu$ is the McDonald Bessel function

$$K_\nu(y)= {1\over 2} \int_0^\infty  e^{-y(t+{1\over t})/2} t^\nu {dt\over t} \qquad (y>0)\, $$
with parameter $\lambda={1\over 2} + \nu $
and the $a_n$ are complex numbers satisfying the Hecke bound

\begin{equation} \label{eq=hb} |a_n| \ll |n|^{1\over 2} \, .\end{equation}

As a final piece of information we need the asymptotic expansion of
the Bessel function

\begin{equation}\label{eq=Kas} K_\nu(y)\sim \left({\pi\over 2y}\right)^{1\over 2}
e^{-y} \cos (\nu \pi)\, . \end{equation}
We can now prove  Theorem \ref{th=mcf}.

\begin{proof} We plug the estimates (\ref{eq=hb}) and (\ref{eq=Kas}) in the
Fourier expansion (\ref{eq=FB}) and use the convention that
$C$ denotes a positive constant whose actual value may change from line
to line: for $y$ large we obtain

\begin{align*} |\phi(x+iy)| &\leq \sum_{n\neq 0} |a_n| \sqrt{y} \cdot |K_\nu(2\pi |n|y)|  \\
&\leq C \sum_{n\neq 0} |n|^{1\over2} \sqrt{y} \left({\pi\over 2\pi  |n|y}\right)^{1\over 2}
e^{-2\pi|n|y} \\
&\leq C \sum_{n\neq 0} e^{-2\pi|n|y} \\
&\leq C {e^{-2\pi y}\over 1- e^{-2\pi y }}\\
&\leq C e^{-2\pi y} \,.\end{align*}
\end{proof}

\subsection{The method of analytic continuation}

We now present an alternative approach to Theorem \ref{th=mcf},
essentially due to J. Bernstein, which uses the method of analytic continuation.
The final result is slightly weaker than the optimal estimate
in Theorem \ref{th=mcf}, but this will be balanced by the
conceptionality of the approach.

\par Let $\phi$ be a Maa\ss{} cusp form. Let us fix $y>0$ and consider
the $1$-periodic function

$$F_y: \R \to \C, \ \ u\mapsto \phi(n_ua_y(i))=\phi(u+iy)\, .$$
This function being smooth and periodic admits a Fourier expansion

$$F_y(u)=\sum_{n\neq 0} A_n(y) e^{2\pi i n x}\, .$$
Here,  $A_n(y)$ are complex numbers depending on $y$.
Now observe that
$$n_ua_y=a_y a_y^{-1}n_u a_y= a_y n_{u/y}$$
and so
$$F_y(u)=\phi(a_y n_{u/ y}.x_0)\, .$$
As $\phi$ is a $\DD(X)$-eigenfunction, it admits a
holomorphic continuation to $\Xi$ and thus
it follows from
Lemma \ref{lem=or} and Theorem \ref{th=sl2} that
$F_y$ admits a holomorphic continuation to the
strip domain
$$S_y=\{ w=u+iv \in \C\mid |v|< y\}\, .$$
Let now $\e>0$, $\e$ small. Then, for $n>0$, we
proceed with Cauchy

\begin{align*} A_n(y) & =\int_0^1 F_y( u -i(1-\e)y) e^{-2\pi i n (u-i(1-\e)y)}\ du \\
&= e^{-2\pi n (1-\e)y} \int_0^1 F_y( u-i(1-\e)y) e^{-2\pi i n u} \ du \\
&= e^{-2\pi n (1-\e)y} \int_0^1 \phi ( a_y n_{u/y} n _{-i(1-\e)}.x_0) e^{-2\pi i n u} \ du \, .
\end{align*}

Thus we get, for all $\e>0$ and $n\neq 0$ the inequality

\begin{equation}\label{eq=ineq}
|A_n(y)|\leq e^{-2\pi |n| y(1-\e)} \sup_{\Gamma g\in \Gamma\backslash G} |\phi(\Gamma g n_{\pm i(1-\e)}.x_0)|
\end{equation}

We need an estimate.

\begin{lem} \label{lem=esti}Let $\phi$ be a Maa\ss{} cusp form. Then there exists a
constant $C$ only depending on $\lambda$ such that
for all $0< \e < 1$
$$  \sup_{\Gamma g\in \Gamma\backslash G} |\phi(\Gamma g n_{i(1-\e)}.x_0)|\leq  C |\log \e|^{1\over 2}$$
\end{lem}

\begin{proof} Let $-\pi/4< t_\e <\pi/4$ be such that $\pm (1-\e)=\sin 2t_\e$.
Then, by Lemma \ref{lem=or}  we have
$G n_{\pm i(1-\e)}.x_0 =G a_{\e}.x_0$ with
$a_\e=\begin{pmatrix} e^{it_\e} & 0 \\ 0& e^{-it_\e}\end{pmatrix}$.
Now note that $t_\e \approx \pi/4 -\sqrt{2\e}$ and thus
\cite{KSI}, Th. 5.1 and Th. 6.17 , give that
$$\sup_{\Gamma g\in \Gamma\backslash G} |\phi(ga_\e.x_0)|\leq C |\log \e|^{1\over 2}\, .$$
This concludes the proof of the lemma.\end{proof}

We use the estimates in Lemma \ref{lem=esti} in (\ref{eq=ineq}) and
get

\begin{equation}\label{eq=ineq1}
|A_n(y)|\leq C e^{-2\pi |n| y(1-\e)} |\log \e|^{1\over 2}\, ,
\end{equation}
and specializing to $\e=1/y$ gives that
\begin{equation}\label{eq=ineq2}
|A_n(y)|\leq C e^{-2\pi |n| (y-1)} (\log y)^{1\over 2} \, .
\end{equation}
This in turn yields for $y>2$ that

\begin{align*} |\phi(iy)| & = |F_y(0)|\leq \sum_{n\neq 0} |A_n(y)|\\
&\leq C (\log y)^{1\over 2}\sum_{n\neq 0} e^{-2\pi |n| (y-1)}\\
&\leq C (\log y)^{1\over 2} \cdot e^{-2\pi y}\end{align*}
It is clear, that we can replace $F_y$ by $F_y(\cdot +x)$ for any $x\in \R$
without altering the estimate. Thus
we have proved:

\begin{thm} Let $\phi$ be a Maa\ss{} cusp form. Then there
exists a constant $C>0$, only depending on $\lambda$,  such that
$$|\phi(x+iy)|\leq C (\log y)^{1\over 2} \cdot e^{-2\pi y} \qquad (y>2)\, .$$
\end{thm}

\begin{rem} It is not too hard to make the constant in the theorem precise. We will do this
in the next section when we give a general discussion of the
rank one cases. \end{rem}

\section{Exponential decay of Maa\ss{} cusp forms II: the rank one cases}

The example of $G=\Sl(2,\R)$ admits a straightforward generalization
to all rank one cases and this will be outlined below.
Throughout this section we let $G$ be of real rank one, i.e. $\dim \af=1$.
We fix a noncocompact lattice $\Gamma<G$ and call a parabolic
subgroup $MAN$ {\it cuspidal} for $\Gamma$ if $\Gamma\cap N$ is a
lattice in $N\cap \Gamma$. Notice that this implies that
$\Gamma\cap Z(N)$ is a lattice in $Z(N)$ where $Z(N)$ is the center
of $N$.
Recall the constant $c={1\over 4 (p+4q)}$ and let
$d={1\over \sqrt c}$ if $q=0$ and $d={1\over 2\sqrt{c}}$ otherwise.
We define the {\it period } $r_\Gamma$ of $\Gamma$ to be the
positive number

$$r_\Gamma={1\over d} \min\{ \|\log \gamma\|: \gamma\in Z(N)\cap \Gamma, \gamma\neq {\bf 1},
N \ \hbox{cuspidal}\}\, .$$

We fix now $MAN$ and $E'\in \log (Z(N)\cap \Gamma)$, $E'\neq 0$, such that
$\|E'\|$ is minimal for all possible choices of $N$. Then
$\|E'\|= d r_\Gamma$.
\par Next let $E\in \R^+ E'$ be such that for $F=\theta(E)$ and $H=[E,F]$ the set
$\{ H, E, F\}$ forms an $\sl(2)$-triple. Recall from
Remark \ref{rem=norm} that $\|E\|=d$ so that
\begin{equation} \label{eq=period} E' =r_\Gamma E\, . \end{equation}
For $y>1$ we set $a_y=\exp(\log y\cdot  H/2)\in A$.

\par We fix a Maa\ss{} cusp form  $\phi$ for $\Gamma$, fix $n\in N$ and $y>1$ and consider
the function

$$F_{n,y}: \R\to \C, \ \ u\mapsto \phi(\exp(uE)na_y.x_0)\, .$$
{}From the relation (\ref{eq=period}), it follows that
$F_{n,y}$ is periodic with period $r_\Gamma$. Thus
$F_{n,y}$ admits a Fourier expansion

$$F_{n,y}(u)=\sum_{k\in \Z^\times} A_k(n,y) e^{{2\pi i k\over  r_\Gamma} u}\, .$$
As $\exp (uE)\in Z(N)$ we notice next that
$$\exp(uE)na_y=na_y \exp(u/y E)\, $$
and we conclude with Theorem \ref{th=r1} that $F_{n,y}$ extends a
holomorphic function on the strip domain
$S_y=\{ u+iv\in \C\mid |v|<y\}$. For $0<\e<1$ we obtain, as in
the previous section, the coefficient estimate

\begin{equation}\label{esti=c1}
|A_k(n,y)| \leq e^{-{2\pi y(1-\e)\over r_\Gamma}} \sup_{\Gamma g\in \Gamma\backslash G}
|\phi(\Gamma g \exp(i(1-\e)E).x_0)|\, .
\end{equation}

\par From now on we make the slightly restrictive assumption that
$\phi$ corresponds to a spherical principal series representation
$\pi_\lambda$ with $\lambda\in i\af^*$. Often we will identify
$i\af^*$ with $i\R$ via $\lambda=\lambda\cdot \rho$.

\begin{lem}\label{lem=esti1} Let $\phi$ be a Maa\ss{} cusp form associated
to $\pi_\lambda$. Then for all $0<\e<1$ the following estimate
holds
\begin{equation}
\sup_{\Gamma g\in \Gamma\backslash G} |\phi (\Gamma g\exp(i(1-\e)E).x_0)|
\leq C(\lambda) \begin{cases} |\log \e|^{1\over 2} & \text{if $p=1$ and $q=0$} \\
\e^{{1-p\over 4}} & \text{if $p>1$ and $q=0$}\\
 |\log \e|^{1\over 2} & \text{if $q=1$} \\
\e^{ {1-q\over 4}} & \text{if $q>1$}\end{cases}
\end{equation}
where
\begin{equation} C(\lambda)=C \cdot e^{{\pi\over 2}|\lambda|} (|\lambda|+1)^{1+[\dim X/2]}
\end{equation}
and $C>0$ a constant independent of $\lambda$.
\end{lem}

\begin{proof} In first order approximation we
have $G\exp(i(1-\e)E).x_0=G\exp(i (\pi/4 -2\sqrt\e)H).x_0$ as in the proof
of Lemma \ref{lem=esti}. Now for fixed $\lambda$, the assertion
follows from \cite{KSI}, Th. 5.1  (or alternatively from our table
in Theorem \ref{thm:table}) and Th. 6.17.
The precise shape of the constant $C(\lambda)$ is found by tracing
the proofs in \cite{KSI}.\end{proof}

Finally, specializing to $y={1\over \e}$ in (\ref{esti=c1}) we obtain
from Lemma (\ref{lem=esti1}) the following result:

\begin{thm} Let $\phi$ be a Maa\ss{} cusp form associated
to $\pi_\lambda$ with $\lambda\in i\af^*$. Then, for all
$n\in N$ and $y>2$ the following estimate holds:
\begin{equation}
\sup_{n\in N} |\phi (na_y.x_0)|
\leq C(\lambda) e^{-{2\pi y\over r_\Gamma}} \begin{cases} (\log y)^{1\over 2} & \text{if $p=1$ and $q=0$} \\
y^{{p-1\over 4}} & \text{if $p>1$ and $q=0$}\\
 (\log y)^{1\over 2} & \text{if $q=1$} \\
y^{{q-1\over 4}} & \text{if $q>1$}\end{cases}\, .
\end{equation}
\end{thm}

\section{Exponential decay of Maa\ss{} cusp forms III: the higher
rank cases}

Throughout this section we denote by $G$ a simple Lie group
and by $\Gamma<G$ a noncocompact lattice.
We say that a parabolic subgroup $P=MAN$ is {\it cuspidal}
if $\Gamma_N=\Gamma\cap N$ is a lattice in $N$.

\begin{dfn} A $\Gamma$-invariant
smooth ${\mathbb D}(X)$-eigenfunction on $X$ is called a {\it weak
Maa\ss{} automorphic form}. A weekly automorphic Maa\ss{} form
is called a {\it Maa\ss{} cusp form} if it is of moderate
growth
and
$$\int_{\Gamma_N\backslash N} f(\Gamma_N n g)\ d(\Gamma_Nn)=0$$
for all $g\in G$ and all proper cuspidal parabolic
subgroups $P=MAN$.  \end{dfn}

\begin{rem} The crucial fact for us is that
all ${\mathbb D}(X)$-eigenfunctions on $X$ extend holomorphically
to $\Xi$ (cf.\ \cite{KSII}, Th. 1.1). Hence all
weak Maa\ss{} automorphic forms extend to holomorphic
functions on $\Xi$. If moreover $\Gamma$ is torsion free, then
$\Gamma$ acts properly on $\Xi$ (as the $G$-action is proper) and
we can form the quotient $\Gamma\backslash \Xi$ in the category of
complex manifolds. Thus Maa\ss{} forms have
$\Gamma\backslash \Xi$ as their natural domain of definition.
\end{rem}

For the rest of this section we let $P=MAN$ be a minimal
parabolic subgroup which is cuspidal. In addition
we make the following

\smallskip
\noindent {\bf assumption}: For each root
$\alpha\in \Sigma^+$ the group
$\Gamma\cap \exp(\gf^\alpha)$ is a lattice
in $\exp (\gf^\alpha)$.

\par \smallskip
For each $\alpha\in \Sigma^+$ and $E_\alpha\in \gf^\alpha$
we set $F_\alpha=\theta(E_\alpha)$ and
$H_\alpha=[E_\alpha,F_\alpha]$. We always normalize $E_\alpha$ in such
a way that $\{ E_\alpha, F_\alpha, H_\alpha\}$ forms an
$\sl(2)$-triplet.

\par For $\alpha\in \Pi$ we define an ideal in $\Sigma$ by

$$\Sigma_\alpha=\left\{ \beta=\sum_{\gamma\in\Pi} n_\gamma \gamma\mid
n_\alpha>0\right\}\subset \Sigma^+, $$
and write $\u_\alpha=\bigoplus_{\beta\in \Sigma_\alpha} \gf^\beta$
for the corresponding ideal in $\nf$. We set $U_\alpha=\exp(\u_\alpha)$
and notice that
$U_\alpha$ is the nilradical of the cuspidal parabolic
subgroup $P_\alpha$ attached to $\alpha\in\Pi$.

\par Associated to
$\alpha\in\Pi $ we define positive constants

$$r_{\alpha, \Gamma}=\max_{\beta\in\Sigma_\alpha}
\min\{ c>0\mid \exp(cE_\beta)\in \Gamma, \
E_\beta\in \gf^\beta \ \hbox{normalized}\}\, , $$
and

$$c_\alpha=\min\left\{ c>0\mid {c\over 2} H_\beta\in \partial\Omega
,\  \beta\in\Sigma_\alpha\right\}\, .$$

\begin{rem} The relevance of the number $c_\alpha$ is the following:
it is the maximal number such that $\exp(itE_\beta).x_0\in \Xi$
for all $0\leq t<c_\alpha$ and $\beta\in\Sigma_\alpha$.
\end{rem}

For a subgroup  $U<N$ with $\Gamma_U\backslash U$
and  compact
we define the {\it constant term} of a function
$f\in C^\infty(\Gamma_N \backslash G)$ with respect
to $U$ as
$$\pi_U f (Ug)=\int_{U_N\backslash U} f(\Gamma_N ug) \ d(U_N u)\, .$$
Note that $\pi_U f \in C^\infty (\Gamma_N U\backslash G )$.
For $U=U_\alpha$ we use
the simplifying notation $\pi_\alpha=\pi_{U_\alpha}$
for the constant term with respect to $U_\alpha$.

\par We can now state the holomorphic analog
of the Main Lemma (Lemma 10) in \cite{HC}.

\begin{lem}\label{mlem} {\rm (Main Lemma)} Let $\alpha\in \Pi$ and $0<\e<1$. Let
$f\in C^\infty(\Gamma_N\backslash G)$ such that
$f$ admits a holomorphic continuation
to $\Gamma_N\backslash \tilde \Xi$. Then there exists
a constant $C_\alpha>0$, only depending on $\alpha$,  such that
\begin{equation}\label{ce}
\left |(f-\pi_\alpha f)(\Gamma_N a)\right|\leq
C_\alpha e^{-{a^\alpha(1-\e)\over r_{\alpha, \Gamma}}} \cdot
\sup_{g\in G\atop \beta\in\Sigma_\alpha} \left|f(\Gamma_N
g\exp(i (1-\e)c_\alpha E_\beta))\right|\end{equation}
for all $a\in A^+$.
\end{lem}

\begin{proof} We follow \cite{HC}, Ch. I,  $\S$ 7.
We order the roots of $\Sigma_\alpha$, say

$$\beta_1> \beta_2>\ldots > \beta_s\, ,$$
and form ideals of $\nf$ by
$$\nf_i=\bigoplus_{j=1}^i \gf^{\beta_j}\qquad (0\leq i \leq s)\, .$$
Set $U_i=\exp(\nf_i)$. Note that $\Gamma\cap U_i$ is cocompact
in $U_i$ by our assumption on $\Gamma$. Thus $\pi_{U_i}$ is defined.
Note that $\pi_{U_0}=\mathrm{id}$ and $\pi_{U_s}=\pi_\alpha$.
We now verify the stronger statement (cf. \cite{HC}, Lemma 19):

\begin{equation}\label{rce}
f- \pi_{U_j} f \quad \hbox{satisfies (\ref{ce})}\, .
\end{equation}
We prove (\ref{rce}) by induction, following the arguments
for the proof of \cite{HC}, Lemma 19. The case $i=0$
is clear.
Notice that $\u_i=\u_{i-1} \oplus \gf^{\beta_i}$. Choose
a basis $E_1, \ldots , E_p$ of $\gf^{\beta_i}$ of normalized
elements. We require in addition that $\exp(r_{\alpha, \Gamma}E_j)
\in \Gamma_N$.
For $0\leq j\leq p$ we set

$$\v_j=\u_i \oplus\bigoplus_{k=1}^j \R E_k\quad\hbox{and}\quad
V_j=\exp (\v_j)\, .$$
Observe that each $V_j$ is a normal
subgroup of $N$ with $\Gamma\cap V_j <V_j$ cocompact.
In particular $\pi_{V_j}$ is defined.
Set $\phi_{f,j}=\pi_{V_j} f$ for $0\leq j\leq p$
and
$$\psi_{f,j}=\phi_{f,j-1}-\phi_{f,j}\qquad (1\leq j\leq p)\, .$$
Then
$$\pi_{U_{i-1}} f -\pi_{U_i} f = \phi_{f,0}-\phi_{f,p}
=\sum_{j=1}^p \psi_{f,j}$$
and
$$f-\pi_{U_j}f =\sum_{i=1}^{j} \pi_{U_{i-1}}f -\pi_{U_i}f $$
imply that it is sufficient
to establish for all $1\leq j\leq p$ (cf. \cite{HC}, Lemma 20):

\begin{equation}\label{rce2}
\psi_{f,j} \quad \hbox{satisfies (\ref{ce})}\, .
\end{equation}
For that let us fix $j$ and write
$\phi_f=\phi_{f,j-1}$, $V=V_j$. So $\phi_f=\pi_V f $.
Consider the mapping

$$V_{j-1}\backslash V_j\to\C, \ \ v\mapsto \phi_f(va)$$
and note that this function is left invariant
under $\Gamma\cap V_j$. As $\exp(r_{\alpha,\Gamma}E_j)\in\Gamma\cap
V_j$ we obtain that

$$\phi_f (a)=\sum_{q\in\Z} \vartheta_{f,q} (a)
\quad\hbox{where}
\quad \vartheta_{f,q}(a)={1\over r_{\alpha, \Gamma}}\int_0^{r_{\alpha,\Gamma} }
\phi_f(\exp(t E_j)a)
e^{-{2\pi i q t\over r_{\alpha,\Gamma}}} \ dt\, . $$
We fix $a\in A_+$, $q\in\Z$ and consider the function
$$F_{f,q}(z)=\phi_f(\exp(zE_j) a)
=\phi_f(a\exp(a^{-\beta_j} z E_j))$$
for $z\in\R$. We conclude that $F_{f,q}$ admits a
holomorphic continuation to the strip domain

$$S=\{ z\in \C \mid |\Im z|<  c_\alpha\cdot  a^{\beta_j}\}\, .$$
Thus, as in the preceding two sections, we obtain for
$0<\e<1$ the estimate

$$|\theta_{f,q}(a)| \leq  M_\e \cdot e^{-2\pi q a^{\beta_j} (1-\e) c_\alpha\over {r_{\alpha,\Gamma}}}$$
where
$$M_\e=\sup_{g\in G} |f(\Gamma_N g\exp(i(1-\e) c_\alpha E_j)|\, .$$
We sum up the geometric series and note that $\vartheta_{f,0}=\phi_{f,j}$
and obtain the desired estimate for $\psi_{f,j}=\phi_f-\phi_{f,j}$.
This proves the lemma.
\end{proof}

This lemma has an an immediate consequence the following
important result (compare to \cite{HC}, Corollary to Lemma 10).

\begin{cor} {\rm (Main Estimate)}\label{cor=M} Suppose that $f$ is Maa\ss{ } cusp form. Then there
exist a constant $C>0$, independent from $f$, such that
for all $a\in A^+$
\begin{equation} \label{ME1}|f(\Gamma a)|\leq C\cdot \min_{\alpha\in \Pi} e^{ - {2\pi (1-\e) a^\alpha c_\alpha\over
r_{\alpha, \Gamma}}}\cdot M_\e
\end{equation}
with
\begin{equation}\label{M1} M_\e=\sup_{g\in G} \sup_{\beta\in\Sigma_\alpha\atop \alpha\in\Pi}
|f(\Gamma g \exp(i(1-\e)c_\alpha E_\beta)|\, .\end{equation}
\end{cor}

\begin{ex} It is instructive to consider the following example
$$G=\Sl(n,\R)\qquad \hbox{and} \quad \Gamma=\Sl(n,\Z)\, .$$
In this situation we have $r_{\alpha,\Gamma}=c_\alpha=1$ for all
$\alpha$ and the estimate in the Corollary becomes
$$|f(\Gamma a)|\leq C\cdot e^{ - 2\pi (1-\e) \cdot \max_{1\leq i\leq n-1} {a_{i}\over a_{i+1}}}\cdot M_\e$$
with
$$M_\e=\sup_{g\in G} \sup_{1\leq i<j\leq n} |f(\Gamma g \exp(i(1-\e)E_{ij}))|\, .$$
\end{ex}

\subsection{Refinements of the Main Estimate}\label{ss=me}

It is possible to do a little bit better as in the Main Estimate  once we apply
the more refined geometric results from Subsections \ref{ssh} and \ref{sssl}.
To state the inequalities  in a more compact form we define

$$r_\Gamma=\sup_\alpha r_{\alpha,\Gamma}\,. $$

We begin with the Hermitian cases, i.e. where $\Sigma$ is of type $C_n$.
We restrict our attention to the maximal Siegel parabolic with abelian
nilradical and proceed as in Lemma \ref{mlem} going simultaneously
in the direction of  the strongly orthogonal $E_j$ (cf. the notation in
Subsection \ref{ssh}). Then Theorem \ref{ts} gives
the following result.

\begin{lem} \label{mlemref1} {\rm (Main Estimate refined -- the Hermitian case)} Suppose
that $\Xi$ is a Hermitian symmetric space and $f$ is a Maa\ss {} cusp form on $X$.
Then, with the notation of Subsection \ref{ssh}, there exists
a constant $C>0$, independent from $f$,  such that for all
$t_i\geq 0$
\begin{equation} \label{ME2}|f(\Gamma \exp(\sum_{j=1}^n t_jT_j))|\leq C \cdot e^{-{2\pi(1-\e)(\sum_{j=1}^n t_j)\over r_\Gamma}}\cdot M_\e
\end{equation}
with
\begin{equation}\label{M2} M_\e=\sup_{g\in G} \left|f\left(\Gamma g \exp(i(1-\e)\sum_{j=1}^n E_j))\right)\right
|\, .\end{equation}
\end{lem}

Finally we draw our attention to the case of
$G=\Sl(n,\R)$ and our fine geometric results
in Subsection \ref{sssl}. We will state our result
for the Whittaker functionals of a Maa\ss{} cusp form.
It is no loss of generality to assume that $N$, the group of unipotent upper triangular
matrices, is cuspidal for the lattice $\Gamma$.
Let us fix a unitary character $\chi: N\to {\mathbb S}^1$. As
$\chi$ is necessarily trivial on $[N,N]$, it is
clear that $\chi$ is given by a parameter
${\bf m}=(m_1, \ldots, m_{n-1})\in \R^{n-1}$, namely

$$\chi\begin{pmatrix} 1 & t_1 & *      & \ldots & *\\
                              &  \ddots& \ddots & *      & *\\
                              &     &  1     & t_{n-1}& *\\
                              &     &        &        & 1\end{pmatrix}
=e^{2\pi i \sum_{j=1}^{n-1} t_j m_j}\, .$$
In the sequel we assume that $\chi$ is trivial
on $\Gamma_N$.
\par For a cusp form $f$ we then define the
{\it Whittaker function with respect to $\chi$}
by

$$W(f,\chi)(g)=\int_{\Gamma_N \backslash N} f(\Gamma n g) \chi(n) \ d(\Gamma_N n)
\qquad (g\in G)$$
and note that $W(f,\chi)\in C^\infty (G/N, \chi)$.
The obvious application of our
standard technique
yields:

\begin{lem} \label{mlemref2} {\rm (Main Estimate refined -- Whittaker
functionals for the special linear group)}
Let $G=\Sl(n,\R)$ and $f$ be a Maa\ss {} cusp form on $X$.
Then, with the notation of Subsection \ref{sssl}, there exists
a constant $C>0$, independent from $f$,  such that for all
$a=\diag(a_1,\ldots, a_n)\in A^+$
\begin{equation}\label{ME3}
|W(f,\chi) (a)|\leq C \cdot e^{-{2\pi(1-\e)\over r_\Gamma}\sum_{j=1}^{n-1} |m_j|\cdot {a_j\over a_{j+1}}}
\cdot M_\e \end{equation}
with
\begin{equation} \label{M3} M_\e=\sup_{g\in G} \left|f(\Gamma g z(1-\e))\right
|\, .\end{equation}
\end{lem}

\subsection{A Bergman estimate on the local crown domains}

To proceed with our estimates on Maa\ss {} cusp forms
we need to control the quantitities

$$M_\e=\sup_{g\in G} |f(\gamma gn_\e)|$$
for certain $n_\e\in N_\C$.
In order to do so we estimate $M_\e$ against an
$L^2$-norm, which can be controlled in terms of representation theory.

\par We state the result.

\begin{prop}\label{prop=berg} Let $\Gamma$ be a lattice in the
semi-simple group $G$. Fix an element $Z_0\in \partial \pi \Omega/2$
and a constant $0<\e<1$. Let
$f$ be a $\Gamma$-invariant holomorphic
function on $\Xi$. Then there exists a constant
$C>0$, independent from $f$, such that
for all $g\in G$
\begin{align*} |f(\Gamma &g \exp(i(1-\e)Z_0).x_0)| \leq  C \cdot
\e^{-\dim X +{1\over 2} \mathrm{rank} X +{1\over 2}}\cdot \\
&  \cdot \sup_{\{Z\in\Omega\mid \|Z-(1-\e)Z_0\|<\e/2\}}
\left(\int_{\Gamma\backslash G} |f(\Gamma g\exp(iZ).x_0)|^2 \ d(\Gamma g) \right)^{1\over 2}
\, .\end{align*}
\end{prop}

Before we start with the proof let us recall the
basic Bergman estimate for polydiscs in $\C^n$.
Fix $z_0\in \C^n$. For $r>0$ let us define the polyydisc centered at $z_0$  with radius $r$ by
$$P(z_0, r)=\{ z\in \C^n \mid \|z-z_0\|_\infty <r\}\, .$$
One expands a holomorphic function $f\in \O(P(z_0,r))$ in a power series
at $z_0$, and uses orthogonality of the monomials; the result
is the {\it Bergman estimate}

\begin{equation}\label{besti}
|f(z_0)|\leq {1\over \pi^n\cdot  r ^n} \cdot \|f\|_{L^2 (P(z_0, r))}\, .
\end{equation}

We turn to the proof of the proposition.

\begin{proof}
We normalize the Killing norm  $\|\cdot\|$ on $\pf$ such that
$\|Z_0 \|=1$. Let $Z_\e=(1-\e) Z_0$. We define  various balls
in $\pf$ and $\pf_\C$:

\begin{align*} B_1 & =\{ U\in \pf\mid \|U-Z_\e\|<\e/2\},\\
 B_2 & =\{ V\in \pf \mid \|V\|<\e/2\},\\
  B & =\{Z=U+iV\mid U\in B_1, V\in B_2\}.
\end{align*}
If necessary we may replace $\e$ by $c\e$ for some positive constant in the definition
of $B_1$ and henceforth assume that $B_1\subset \pi\hat\Omega/2$.
Then it is clear that $\exp(B_2)\exp(B_1).x_0\subset\Xi$, but what
about $\exp(B).x_0$ ? This is not clear, but after
some controlled shrinking we are in good shape:

\begin{lem} There exists $c>0$ such that for all $0<\e < 1$
$$\exp(B_2)\exp(iB_1).x_0\supset \exp(\{ Z=U+iV\in \pf_\C: \|Z-Z_\e\|< c\e/2\})\, .$$
\end{lem}

\begin{proof} We remark that $\pi \overline {\hat \Omega}/2$ is compact
and that

$$d\exp(iZ): \pf_\C \to T_{\exp(iZ).x_0} X_\C$$
is invertible for all $Z\in \pi \overline {\hat \Omega}/2$.
In fact, the Jacobian of $\exp$ at $iZ$ is given by

$$|\det d\exp(iZ)|=\left | \prod_{\alpha\in\Sigma^+}
{\sinh \alpha(iZ)\over \alpha(iZ)}\right|\, . $$
 The assertion follows from the implicit function theorem.
\end{proof}
It is no loss of generality to assume that
the constant $c$ in the previous lemma is $1$.
At any rate the previous lemma combined with
the Bergman estimate yields

\begin{equation}\label{besti2}
|\phi(\exp(iZ_\e).x_0)|\leq C \cdot{1\over \e^{\dim X}} \left(\int_{\exp(B_2)\exp(iB_1).x_0} |\phi(z)|^2 \ dz\right)
^{1\over 2}\end{equation}
for a constant $C>0$ and all functions $\phi\in \O(\Xi)$.
Here, $dz$ denotes the Haar measure on $X_\C$.
\par Next we set $B_{\af,1}=B_1\cap \af$
and note that
$$\Ad (K) B_{\af,1}\supseteq B_1'\, .$$
with $B_1'= \{ U\in \pf\mid \|U-Z_\e\|<c\e/2\}$ for some
constant $c>0$. Again it is no loss of generality to assume that
$B_1'=B_1$.
As a consequence we derive from (\ref{besti2}) that
\begin{equation}\label{besti3}
|\phi(\exp(iZ_\e).x_0)|\leq C\cdot {1\over \e^{\dim X}} \left(\int_{\Gamma\backslash G \exp(iB_1).x_0}
|\phi(z)|^2 \ dz\right)
^{1\over 2}\end{equation}
for all $\Gamma$-invariant holomorphic functions $\phi$ on $\Xi$.
Finally we use the integration formula \cite{KSII}, Prop. 4.6, and obtain with
$$J(Y)=\prod_{\alpha\in\Sigma^+} |\sin 2\alpha(Y)|^{m_\alpha} \qquad (Y\in\af)$$
that

\begin{equation}\label{besti4}
\int_{\Gamma\backslash G \exp(iB_1).x_0}
|\phi(z)|^2 \ dz =\int_{\Gamma\backslash G}\int_{B_{\af,1}}
|\phi(\Gamma g\exp(iY).x_0)|^2 \cdot J(Y) \ d(\Gamma g) \ dY
\end{equation}
Notice that $J(Y)\leq 2 \e$ as at least one root is going to vanish
on $Z_\e$ for $\e\to 0$.
Thus after combining (\ref{besti3}) and (\ref{besti4}) we obtain
for all $g\in G$  that

\begin{align*} |f(g\exp(iZ_\e).x_0)|& \leq C \cdot \e^{-\dim X+ {1\over 2}(1+ \mathrm{rank} X)}
\cdot \\
& \cdot  \sup_{Y\in B_{\af,1}} \left(\int_{\Gamma\backslash G}
f(\Gamma g\exp(iY).x_0|^2 \ dz\right)
^{1\over 2}\end{align*}
for all $\Gamma$-invariant $f\in\O(\Xi)$.
\end{proof}

\subsection{Main estimates in final form}
In this concluding subsection we put our previously
obtained results together in order to obtain final version
of our main estimates Corollary \ref{cor=M}, Lemma \ref{mlemref1}
and Lemma \ref{mlemref2}.
The main task is to obtain estimates for the quantities $M_\e$ in (\ref{M1}), (\ref{M2}) and (\ref{M3}).
We give details for the main case in
(\ref{M1}), and confine ourselves with stating the analogous results
for the remaining two cases.
So we wish to control the behavior of
$$M_\e=\sup_{g\in G} \sup_{\beta\in\Sigma_\alpha\atop \alpha\in\Pi}
|f(\Gamma g \exp(i(1-\e)c_\alpha E_\beta)|\, .$$
First we deduce from Lemma \ref{lem=or} that
\begin{equation} \label{m1}
M_\e\leq
\sup_{g\in G}\sup_{Y\in \partial\Omega}
|f(\Gamma g \exp(i(1-2\sqrt{\e})\pi Y/2)|\, .\end{equation}
Set
$$r^X:=-\dim X + {1\over 2} \mathrm{rank} X + {1\over 2}\, .$$
Then it follows from Proposition \ref{prop=berg} and
(\ref{m1}) that
\begin{equation} \label{m2}
M_\e\leq  C \cdot \e^{{r^X}/2}
\sup_{Y\in \partial\Omega}
  \cdot
\left(\int_{\Gamma\backslash G} |f(\Gamma g\exp(i(1-c \sqrt{\e})\pi Y/2).x_0)|^2 \ d(\Gamma g) \right)^{1\over 2}
\end{equation}
for constants $C,c>0$ only depending on $X$.
Assume that $f$ corresponds to the spherical representation $\pi_\mu$.
Recall the exponents $s^X$ and $d^X$ from (\ref{def=ds}).
Now, Theorem \ref{thm=se} applies and we arrive at
\begin{equation} \label{m3}
M_\e\leq  C(\mu) \cdot \e^{{r^X}/2+ s^X/4}| \log \e|^{d^X/2}
\end{equation}
for a constant $C=C(\mu)$ depending on $\mu$ and the geometry of $X$.
If we specialize in Corollary \ref{cor=M} to $\e=\min_{\alpha\in \Pi}  a^{-\alpha}$
we get from (\ref{m3}) the following

\begin{thm}\label{mef1}{\rm (Main Estimate)} Suppose that $f$ is Maa\ss{ } cusp form corresponding to
$\pi_\mu$. Then there
exist a constant $C=C(\mu)>0$ such that
for all $a\in A^+$
\begin{equation}\label{ME11} |f(\Gamma a)|\leq C\cdot \min_{\alpha\in \Pi} e^{ - {2\pi a^\alpha c_\alpha\over
r_{\alpha, \Gamma}}}\cdot \max_{\alpha\in \Pi} a^{-\alpha(r^X/2+ s^X/4)}\cdot |\alpha(\log a)|^{d^X/2}\, .
\end{equation}
\end{thm}

In similar manner we obtain a more concrete version of Lemma \ref{mlemref1}:

\begin{thm} \label{mef2} {\rm (Main Estimate refined -- the Hermitian case)} Suppose
that $\Xi$ is a Hermitian symmetric space and $f$ is a Maa\ss {} cusp form on $X$ corresponding
to $\pi_\mu$.
Then, with the notation of Subsection \ref{ssh}, there exists
a constant $C=C(\mu)>0$
 such that for all
$t_i\geq 0$
\begin{equation} \label{ME22}|f(\Gamma \exp(\sum_{j=1}^n t_jT_j))|\leq
C \cdot e^{-{2\pi(\sum_{j=1}^n t_j)\over r_\Gamma}}\cdot
\left(\sum_{j=1}^n t_j\right)^{ -r^X/2 - s^X/4} \, .
\end{equation}
\end{thm}

Finally we state a new version of Lemma \ref{mlemref2}
(for which one also needs to employ the estimate in Proposition \ref{prost}) :

\begin{thm} \label{mef3} {\rm (Main Estimate refined -- Whittaker
functionals for the special linear group)}
Let $G=\Sl(n,\R)$  and $f$ be a Maa\ss {} cusp form on $X$
corresponding to $\pi_\mu$.
Then there exists
a constant $C=C(\mu)$  such that for all
$a=\diag(a_1,\ldots, a_n)\in A^+$
\begin{align*}
|W(f,\chi) (a)|\leq & C \cdot e^{-{2\pi\over r_\Gamma} (\sum_{j=1}^{n-1} |m_j|\cdot {a_j\over a_{j+1}})}\cdot
\left(\sum_{j=1}^{n-1} |m_j|\cdot {a_j\over a_{j+1}}\right)^{-r^X - s^X/2}\\
& \cdot
\left|\log \left(\sum_{j=1}^{n-1} |m_j|\cdot {a_j\over a_{j+1}}\right)\right|
\, .\end{align*}
\end{thm}

\begin{rem}{\rm (Some generalizations of the Main Estimates)}
Let $F$ either denote a Maa\ss{} cusp form $f$ or a Whittaker function
$W(f,\chi)$ in case $G=\Sl(n,\R)$. The general
form of our estimates in Theorems \ref{mef1}, \ref{mef2}, \ref{mef3}, then is

\begin{equation}\label{111}
(\forall a\in A^+) \qquad |F(\Gamma a)|\leq C e^{-\phi(a)} \cdot P(a)
\end{equation}
where $C>0$ is a constant only depending on the representation
$\pi_\mu$ associated to $F$,
$$\phi(a)=\sum_{\alpha\in\Sigma^+} c_\alpha  a^\alpha \qquad (c_\alpha\geq 0)$$
is a "positive" linear functional and $P(a)$ is a polynomial
in the variables $a^\alpha, \alpha(\log a)$ with $\alpha\in\Sigma^+$.

\par\noindent (a) {\rm (Extension to a Siegel domain)} We restricted ourselves
to estimates on $A^+$. However, for certain applications in number
theory one needs estimates which are uniform on a  Siegel domain
$${\mathfrak S}_t=\omega A_t K \qquad (t>1)$$
where $\omega\subset N$ is a fixed compact and $A_t=\{ a\in A\mid
a^\alpha> t\  \forall \alpha \in \Pi\}$.
The version of (\ref{111}) on the whole Siegel domain ${\mathfrak S}_t$
is
\begin{equation}\label{222}
(\forall n\in\omega) (\forall a\in A_t) \qquad |F(\Gamma na)|\leq C_t e^{-c_t\phi(a)} \cdot P(a)
\end{equation}
where $c_t,C_t>0$ are such that $c_t\to 1^- $ for $t\to \infty$.
Let us explain how this is derived from (\ref{111}).
For the proof of (\ref{111}) we use certain subsets $\Lambda_0\subset \Lambda$
and applied the fact that $F$ extends to a function on
$\Gamma \exp(i\Ad(a) \Lambda_0) a.x_0$; recall, the precise rate of exponential
decay was directly linked to the geometry of $\Lambda_0$.
If one wants estimates on ${\mathfrak S} _t$ one needs to bring in
$\omega$-variables, i.e. we look for maximal subsets  $\Lambda_t\subset\Lambda_0$
such that $F$ extends to $\Gamma \exp(i \Ad(a)\Lambda_t) \omega a.x_0$.
This is equivalent to the requirement of
$$\exp(i\Lambda_t) a\omega a^{-1}.x_0\subset \Xi$$
Now $a\omega a^{-1}$ shrinks to $\{\bf 1\}$ for $t\to \infty$, meaning $\bigcup_{t>0} \Lambda_t=\Lambda_0$.

\par\noindent (b) {\rm (Extension to other $K$-types)}
We only considered Maa\ss{} cusp forms, i.e. cusp forms associated
to a trivial $K$-type. However, it is possible to extend
to other $K$-types $\sigma\in\hat K$. For that one needs uniform quantitative
control of the projection $\kappa: N_\C A_\C K_\C \to K_\C$
on $K$-orbits through $\exp(i\pi \Omega/2)$. This will be defered
to another paper. In any case,
the result then is

\begin{equation}\label{333}
(\forall nak\in{\mathfrak S}_t) \qquad |F(\Gamma nak)|\leq C_t e^{-c_t\phi(a)} \cdot P_\sigma (a)
\end{equation}
with $C_t=C_t(\mu,\sigma)$ and $P=P_\sigma$ now depending on $\sigma$.
Further explanation is given in item (c) below.

\par\noindent (c) {\rm (Extension to non-spherical representations)} So far
we only considered spherical representations $\pi=\pi_\mu$.
But estimate (\ref{333}) remains true for an arbitrary irreducible
unitary representation $\pi$ and arbitrary $K$-types $\sigma$.
Here is the reason. We can embed $\pi$ into a principal series
representation induced off a minimal parabolic (subrepresentation theorem).
One uses the fact that the smooth structures are unique (Casselman-Wallach).
Now for a principal series one can look at the corresponding
Eisenstein integrals for the $K$-types (Harish-Chandra) and everything
boils down to estimate the spherical function and the sup-norm of a
holomorphically extended $K$-type (for that one needs the quantative
control of $\kappa$). The other details can be found in the
proof of \cite{KSI}, Th. 3.1.
\end{rem}

\begin{rem}\label{rem11} {\rm (Quality of the estimates)} The rate of exponential
decay given in above three Theorems are sharp. We provide some evidence
in the section below. Concerning the polynomial part, one could likely
replace $r^X$ by zero. However that would require to prove
Conjecture C in \cite{KSI}; something which is
out of reach with currently available techniques.
\end{rem}

\begin{rem} {\rm (Applications to automorphic forms)}
\par\noindent (a) {\rm ($L$-functions)}
For various reasons on wants to know whether certain
automorphic $L$-functions are meromorphic of finite order.
For instance this information is required if one wants
to exhibit zero-free
regions (in the spirit of de la Vall\'ee Poussin)  for those $L$-functions.
We refer to
\cite{GLS}, \cite{GL}, \cite{GS} for results
in this direction. We wish to point out that our estimates help
to establish that $L$-functions with appropriate
integral representations are in fact of finite
order.
\par\noindent (b) {\rm (Voronoi summation)} Recently Voronoi
summation was established for ${\mathrm {Gl}}(3,\R)$, cf. \cite {MS}.
Shortly after it was extended to ${\mathrm {Gl}}(n,\R)$ in
\cite{GLi}. In the approach of \cite{GLi} it becomes visible that
exponential decay is an important analytical ingredient
to establish Voronoi summation.
\end{rem}

\section{Final Remarks}

\subsection{Estimates on Whittaker functionals for $\Gl(n)$ are sharp}

We show that the rate of exponential decay for Whittaker functionals for $G=\Gl(n,\R)$ proved in Theorem \ref{mef3}
is optimal.
\par To begin with we recall the Whittaker expansion of Piatetski-Shapiro and Shalika
for a cuspidal Maa\ss{} form of the group
$\Gl(n,\R)$. For simplicity let us restrict ourselves to the case $n=3$. Our arithmetic subgroup of choice will be
$\Gamma=\Gl(3,\Z)$. Let us define subgroups of $\Gamma$ by

$$\Gamma^2=\begin{pmatrix} \Gl(2,\Z) & 0\\ 0 & {\bf 1}
\end{pmatrix}\quad\hbox{and}\quad \Gamma^2_N=\Gamma^2\cap N\, .$$

In the sequel we use the notation introduced in Subsection \ref{ss=me}. For $\chi$ corresponding to ${\bf m}=(1,1)$ and
$f$ a Maa\ss {} cusp form we set $W(f)=W(f,\chi)$. For $n_1,n_2\in\N$ we put

$$W_{n_1,n_2}(f)(z)=W(f)\left( \begin{pmatrix} n_1 n_2 & & \\ & n_1 & \\
& & 1\end{pmatrix} z\right)$$
where $z\in X=\Gl(3,\R)/ {\rm O}(3,\R)$. The Whittaker expansion of $f$ reads as

$$f(z)=\sum_{\gamma\in \Gamma_N^2\backslash \Gamma^2}\sum_{n,1,n_2\in\N}
{a_{n_1,n_2}\over n_1 n_2} \cdot W_{n_1,n_2}(f)(z)$$ for complex coefficients $a_{n_1, n_2}$, \cite{Sh}, Th. 5.9. We
normalize $f$ such that $a_{1,1}=1$ and draw our attention to the main result in \cite{Bu}, (10.1), which gives a
formula for the Mellin transform of $W(f)$:

\begin{align}\label{eq=bump}& \int_0^\infty \int_0^\infty W \begin{pmatrix}y_1 y_2 & & \\
& y_1 &\\  & & 1\end{pmatrix}\cdot   y_1^{s_1-1}\cdot  y_2^{s_2-1}
\ {dy_1\over y_1} {dy_2\over y_2}=\\
\nonumber&\quad = {1\over 4} \pi^{-s_1-s_2} \cdot{\Gamma\left({s_1+\alpha\over 2}\right)
\Gamma\left({s_1+\beta\over 2}\right)\Gamma\left({s_1+\gamma\over 2}\right)
\cdot\Gamma\left({s_2-\alpha\over 2}\right)
\Gamma\left({s_2-\beta\over 2}\right)\Gamma\left({s_2-\gamma\over 2}\right)
\over \Gamma\left({s_1+s_2\over 2}\right)}\, .
\end{align}
Here $s_1,s_2$ are sufficiently large real numbers and $\alpha,\beta,\gamma=-\alpha-\beta$ are complex numbers related
to the parameter of the principal series representation associated to $f$
 (see \cite{Bu}, p. 161).
We perform a Stirling approximation of the right hand side (RHS) of (\ref {eq=bump})
 and obtain

$$(RHS)(s_1,s_2)\sim {1\over 4} \pi^{-s_1-s_2} \cdot \sqrt{2\pi}^5 \cdot
e^{-(s_1+s_2)}\cdot { \left({s_1\over 2}\right)^{3({s_1\over 2} -{1\over 2})}
\left({s_2\over 2}\right)^{3({s_2\over 2} -{1\over 2})} \over
 \left({{s_1+s_2}\over 2}\right)^{({{s_1+s_2}\over 2} -{1\over 2})}}\, .$$
We specialize to $s_1=s_2=s$ and get the simpler expression

\begin{equation}\label{eq=a1}
(RHS)(s,s)\sim {1\over 2} (2\pi)^{-2s+5/2}  \cdot e^{-2s}\cdot {s^{2s-1/2} \over 2^s}\, .
\end{equation}

Similarly, if we keep one variable fixed to be zero we get

\begin{equation}\label{eq=a2}
(RHS)(s,0)= RHS(0, s)\sim C
(2\pi)^{-s}  \cdot e^{-s}\cdot {s^{s-1}}\, .
\end{equation}

We wish to compare these asymptotics with what we obtain by applying the estimate for $W(f)$ from
Theorem \ref{mef3}
 With $N>3$ we get
 \begin{equation} \label{eq=a3}
 \left|W(f)\begin{pmatrix} y_1y_2 & & \\ & y_1 & \\ & & 1\end{pmatrix}\right|
 \leq C (1+ y_1^N + y_2^N ) e^{-2\pi (y_1+y_2)} \, .\end{equation}

We insert the estimate (\ref{eq=a3}) into the left hand side (LHS) of \ref{eq=bump} and  arrive at the inequality
\begin{equation}\label{eq=a4}
LHS(s,s)\leq C (2\pi)^{-2 s} e^{-2s} (s+N-1)^{2s+2N-3}
\end{equation}
and likewise
\begin{equation}\label{eq=a5}
LHS(s,0)\leq C (2\pi)^{- s} e^{-s} (s+N-1)^{s+N-3/2}\, .\end{equation}

Conjecturally we could even take any $N>0$ (cf.\  Remark \ref{rem11}). In any case, if
we compare (\ref{eq=a5}) with (\ref{eq=a2}) we see that the  exponential decay for the Whittaker functional established
in Theorem \ref{mef3}  is optimal. In fact, any better exponential decay rater would
lead to decrease of $(2\pi)^{-2 s}$ to $(2\pi +a)^{-2s}$ for some $a>0$ on the right of
(\ref{eq=a5}); this would contradict the asymptotics in (\ref{eq=a2}).

\subsection{Automorphic holomorphic triple products}
We introduce a holomorphic version of triple products and raise some natural questions. The setting here is: $G$ a
semisimple noncompact Lie group and $\Gamma<G$ a cocompact lattice. For three automorphic forms
 $\phi_1, \phi_2,\phi_3$
one $\Gamma\backslash G$ one forms the automorphic triple product, or automorphic trilinear functional in the
terminology of J. Bernstein and A. Reznikov,

$$\ell_{\rm aut}(\phi_1, \phi_2,\phi_3)=\int_{\Gamma\backslash G} \phi_1(\Gamma g) \phi_2(\Gamma g )
\phi_3(\Gamma g) \ d(\Gamma g)\, . $$

Assume now that the $\phi_i$ are Maa\ss{} forms so that the integral defining $\ell_{\rm aut}$ is effectively over the
locally symmetric space $\Gamma\backslash X$. From the general theory we know that the $\phi_i$ extend to holomorphic
functions $\tilde\phi_i$ on the local crown domain $\Gamma\backslash \Xi$. For the moment  we restrict ourselves to the
basic case of $G=\Sl(2,\R)$ with comments on the general situation thereafter.

We form the {\it holomorphic automorphic triple product} by

$$\ell_{\rm aut}^{\rm hol}(\phi_1, \phi_2,\phi_3)=\int_{\Gamma\backslash \Xi} \tilde \phi_1(\Gamma z)
\tilde \phi_2(\Gamma z ) \tilde \phi_3(\Gamma z) \ d(\Gamma z)\,  $$
where $d(\Gamma z)$ is the measure on $\Gamma\backslash \Xi$ induced from the Haar measure on $X_\C$. That $\ell_{\rm
aut}^{\rm hol}$ is actually defined is content of the next lemma.

\begin{lem}  Let $G=\Sl(2,\R)$ and $\Gamma<G$ be a cocompact lattice.
Let $\phi_1,\phi_2,\phi_3$ be Maa\ss{} automorphic forms. Then the intergral defining $\ell_{\rm aut}^{\rm hol}$
converges absolutely.
\end{lem}

\begin{proof} In view of \cite{KSI}, Th. 5.1 and Th. 6.17, there exists a constant
$C>0$ such that we have for all $Y\in\Omega$

$$\sup_{g\in G} |\tilde \phi_i(\Gamma g \exp(i\pi Y/2).x_0)|\leq C \left |\log \cos \alpha(\pi Y/2)\right| $$
Thus in view of the polar decomposition of the measure $\Xi$ (see \cite{KSII}, Prop. 4.6), we get
$$|\ell_{\rm aut}^{\rm hol}(\phi_1, \phi_2,\phi_3)|\leq C^3 \int_0^1   \left |\log \cos \alpha(\pi Y/2)\right|^3
\cdot \sin \alpha(\pi Y) \ dY <\infty$$ and this proves the lemma.
\end{proof}

\begin{problem} Determine the relation between $\ell_{\rm aut}$ and $\ell_{\rm aut}^{\rm hol}$
and explain its significance.
\end{problem}

\begin{rem} For a general semisimple Lie group the integrals defining
$\ell_{\rm aut}^{\rm hol}$ are not absolutely convergent. However, we have some freedom in the choice of the $G$-invariant
measure on $\Xi$. Under the parameterization map $p: G/M\times \Omega^+\to \Xi$ the pull back of the Haar measure $dz$
on $\Xi$ is given by
$$ p^*(dz)=d(gM)\times J(Z)dZ $$
with $J(Z)=\prod_{\alpha\in\Sigma^+} [\sin \alpha(\pi Z)]^{m_\alpha}$ (see \cite{KSII}, Prop. 4.6). For example if we
replace $J$ by sufficiently high power $J^k$, then $\ell_{\rm aut}^{\rm hol}$ is defined with regard to this measure.
It is possible to carry out the details using the results obtained in this article.
\end{rem}

\subsection{Exponential decay of automorphic triple products}

With the methods of analytic continuation one can prove exponential decay of automorphic
triple products (see \cite{Sa}
and \cite{BR} for the first results). To be a more precise, consider a compact locally symmetric space
$\Gamma\backslash X$ and fix a Maa\ss{} form $\phi$. Then for a Maa\ss{} form $\phi_\pi$ corresponding to an
automorphic representation $\pi$ there is interest in finding the precise exponential decay of

$$|\ell_{\rm aut}(\phi,\phi_\pi,\overline{\phi_\pi})|$$
in terms of the parameter $\lambda(\pi)$ of $\pi$. This was first determined by Sarnak for $G=\Sl(2,\C)$ \cite{Sa} and
then by Petridis \cite{Pet} and Bernstein-Reznikov \cite{BR} for $G=\Sl(2,\R)$. Optimal bounds for all rank one groups
were established in \cite{KSI}. For higher rank groups , such as $\Sl(n,\R)$,  partial results were obtained in
\cite{KSI}. These bounds however fail to be optimal in general. The results in this paper combined with the methods of
\cite{BR} and \cite {KSI} allow to establish non-trivial (although) non-optimal bounds for the exponential decay of
automorphic triple products.
\section{Appendix: Leading exponents of holomorphically extended
elementary spherical functions}\label{app:exp}
In this appendix we prove Theorem \ref{thm:esteta} and the table
of Theorem \ref{thm:table} for the asymptotic behavior of norm of
the holomorphic extension of the orbit map $G/K\ni gK\to\pi(g)v$ of a
spherical vector $v\in\H$ in an irreducible spherical
representation $(H,\pi)$ of
$G$ when the argument approaches the distinguished boundary of the crown
domain $\Xi$.
The key property equation (\ref{eq:orbnorm}) translates this to the problem
of finding the asymptotic behavior of certain solutions of a system of
differential equations when approaching the singular locus of the
system.
In the theory of ordinary Fuchsian differential equations this
boils down to the study of characteristic exponents at its singular
points of $\mathbb{P}^1$ and their relation to the monodromy of the
system. A beautiful application (closely related to our problem
in fact, via (\ref{eq:l2norm})) of this classical theory is the
study of the asymptotic behavior of certain classes of oscillatory
integrals via the monodromy of the Gau\ss-Manin connection of the
Milnor fibration of the phase function \cite{Malg}.

In the case of ``Fuchsian systems'' of differential equations in several
complex variables we first need to develop some fundamental facts
on exponents and their first properties. In the case of a regular holonomic
$\mathcal{D}$-module of the form $\mathcal{D}/\mathcal{J}$ on $\C^n$ which
is $\mathcal{O}$-coherent on the complement of a \emph{hyperplane arrangement}
in $\C^n$ we propose a definition of the set of exponents of local solutions
of the system of equations $D\phi=0,\ \forall D\in\mathcal{J}$ at
any point $\eta\in\C^n$. This translates our original problem to
that of determining the set of exponents of a special solution to
Harish-Chandra's radial system of differential equations on $A_C.x_0$
(namely the holomorphic extension of the restriction of the elementary
spherical function to $A.x_0$) at the extremal boundary points
$t(\eta)^2 x_0$ of $T_\Omega$.

What turns this into a successful method is
the fact that there exists a well behaved parameter deformation of
Harish-Chandra's radial system of differential equations for which
we have rather explicit knowledge of the monodromy representation
of its solutions \emph{for generic parameters}. This deformation
is the hypergeometric system of differential equations
\cite{HO1},\cite{HS},\cite{Op4}. Its monodromy factors through an
affine Hecke algebra, thus bringing the representation theory of
affine Hecke algebras into play. In the spirit of the study of
the Bessel function equations \cite{Op2} this leads to the
description of the set of all exponents of the hypergeometric
system at $t(\eta)^2.x_0$. Using that and the relation
between exponents and monodromy (which we will carefully establish below)
we can compute the leading exponents of the holomorphically extended
hypergeometric function at the points $t(\eta)^2 x_0$.

Specialization of the parameters then leads to the desired lower bounds
for the leading exponents of holomorphically extended elementary
spherical functions on a Riemannian symmetric space $X$,
leading to the proof of Theorem \ref{thm:esteta} and
Theorem \ref{thm:table}.

One remarkable phenomenon that comes out of these considerations
is that the leading exponent of the hypergeometric function at
an extremal point $t(\eta)^2.x_0$ is related to a \emph{leading
character} $\s_\eta$ of the isotropy group $W^a_\eta\subset W$
of $t(\eta)^2.x_0$ which depends only on the geometry of $\Omega$
locally at the extremal point $\eta$, but not on the multiplicity $m$ (if
$m$ is real and satisfies certain inequalities which hold for the
multiplicity functions of Riemannian symmetric spaces).
\subsection{Exponents and hyperplane
arrangements}\label{sub:expAR}
In this subsection we propose a definition of exponents of local Nilsson
class functions \cite[Chapter 6.4]{Bj} on the complement of a hyperplane arrangement of
$\C^n$ at points $\eta\in\C^n$. The main results are that the exponents at
$\eta$ are invariant for local monodromy at $\eta$ and the relation between
exponents and monodromy.

Let $\eta\in\C^n$ and let $\phi$ be a local Nilsson class function
at $\eta$. By this we mean a multivalued holomorphic function
$\phi$ on the complement $N\backslash Y:=N^{\operatorname{reg}}$ of an
analytic hypersurface $Y\subset \C^n$ inside a small open ball
$N\subset \C^n$ centered at $\eta$ such that
\begin{enumerate}
\item[LN1:] $\phi$ has finite determination order in
$N^{\operatorname{reg}}$.
\item[LN2:] The pull back of any branch
of $\phi$ via any holomorphic map $j:\mathbb{D}\to N$ with the property
that $j^{-1}(Y)\subset\{0\}$ has moderate growth at $0\in\mathbb{D}$.
\end{enumerate}
Suppose that $j$ as in LN2 is an embedding such that
$j(0)=\eta$. Then the pull back of (any branch of) $\phi$ via $j$
has a singular expansion at $\e=0$ (where $\e$ denotes the standard
coordinate in the unit disk $\mathbb{D}$) of the form
\begin{equation}\label{eq:singexp}
\phi(j(\e))=\sum_{s,l}\e^s\log^l(\e)f_{s,l}(\e),
\end{equation}
a sum over a finite set of pairs $(s,l)$ with $s\in\C$ and
$l\in\mathbb{Z}_{\geq 0}$, such that for each pair $(s,l)$ in
this sum the function $f_{s,l}(\e)$ is holomorphic on $\mathbb{D}^\times$
with at most a pole at $\e=0$. This expansion is obviously not unique, and even
if one tries to make it unique by imposing additional requirements one will find
that the set $S$ which enters in (\ref{eq:zsingexp}) will in general depend on
the chosen embedding $j$ and of the chosen branch of $\phi$ (with respect
to local monodromy in $N^{\operatorname{reg}}$) in an essential way.
In order to define exponents of $\phi$ at $\eta$ we assume from now on
the following.
\begin{enumerate}
\item[AR:] For sufficiently small $N$ we may take $Y=Y^\eta$ to be
a linear hyperplane arrangement centered at $\eta$.
\end{enumerate}
We call a holomorphic map
$i:\mathbb{D}^{n}\to N$ a standard coordinate map if
\begin{enumerate}
\item $i(0,z)=\eta$ for all $z\in\mathbb{D}^{n-1}$.
\item $i(\mathbb{D}^\times\times\mathbb{D}^{n-1})
\subset N^{\operatorname{reg}}$.
\item The lift of the map
$i:\mathbb{D}^\times\times\mathbb{D}^{n-1}\to N^{\operatorname{reg}}$
to the blow-up $X_\eta\to \C^n$ of $N$ at the point $\eta$
extends to a coordinate map $i:\mathbb{D}^n\to X_\eta$ such that
$i(\mathbb{D}^n)\cap Z=\emptyset$, where $Z$ denotes the strict
transform of $Y\cap N$.
\end{enumerate}
Let $i$ be a standard coordinate map. Choose a base point
$p=i(P)\in i(\mathbb{D}^\times\times\mathbb{D}^{n-1})$ and fix a germ
$\phi_p$ of a branch of $\phi$ at $p$. Let $C\subset\mathbb{D}^\times$
be a cut disk (the complement in $\mathbb{D}$ of a ray emerging
from $0$) such that $P\in C\times\mathbb{D}^{n-1}$.
The pull back of $\phi_p$ via $i$ to $P\in C\times\mathbb{D}^{n-1}$ is
the germ of a Nilsson class function on
$\mathbb{D}^\times\times\mathbb{D}^{n-1}$. Hence we have the
following standard result \cite[Proposition 4.4.2]{Bj}:
\begin{prop} There exists a finite set $S$ of pairs $(s,l)$ with
$s\in\C$ and $l\in\mathbb{Z}_{\geq0}$ such that the unique analytic
continuation of $i^*(\phi_p)$ to $C\times\mathbb{D}^{n-1}$ admits an
expansion of the form
\begin{equation}\label{eq:zsingexp}
\phi_p(i(\e,z^\prime))=\sum_{(s,l)\in S}
\e^s\log^l(\e)f_{s,l}(\e,z^\prime),
\end{equation}
where each $f_{s,l}$ extends meromorphically to
$\mathbb{D}^\times\times\mathbb{D}^{n-1}$.
\end{prop}

As before this expansion is not unique for obvious reasons
but we can rearrange (\ref{eq:zsingexp}) in such a way that
\begin{enumerate}
\item[(a)] all $f_{s,l}$ extend holomorphically on $\mathbb{D}^n$,
\item[(b)] if the pairs $(s,l)$
and $(s^\prime,l^\prime)$ occur in (\ref{eq:zsingexp}) then $s-s^\prime$
is not equal to a nonzero integer, and
\item[(c)] if the pair $(s,l)$ occurs in
(\ref{eq:zsingexp}) then there exists an $l^\prime\in\Z_{\geq 0}$ such that
$f_{s,l^\prime}(0,\cdot)\not\equiv 0$.
\end{enumerate}
That makes the expansion unique.
\begin{dfn}\label{dfn:expphi}
Let $\phi$ be a local Nilsson class function at $\eta\in\C^n$
and let $i:\mathbb{D}^\times\times\mathbb{D}^{n-1}\to
N^{\operatorname{reg}}$ be a standard coordinate map.
Choose a base point $p$ in the image of $i$, and
choose a germ $\phi_p$ of a branch of $\phi$ at $p$.
We define the finite set $E^{\eta,i,p}(\phi_p)\subset\C\cup\{\infty\}$
of exponents of $\phi_p$ at $\eta$ as the projection of the finite set
$S\subset\C\times\mathbb{Z}_{\geq0}$ defined above to the first component
if $\phi\not=0$. We put $E^{\eta,i,p}(0)=\{\infty\}$.
\end{dfn}
\begin{prop}\label{prop:indep}
The set $E^{\eta,i,p}(\phi_p)$ is independent of the choice of $i$
(satisfying the requirements (i),(ii) and (iii) above) and
is independent of analytic continuation of $\phi_p$ within
$N^{\operatorname{reg}}$. Hence we may
speak about the set of exponents $E^{\eta}(\phi)$ without
referring to a specific branch of $\phi$ and coordinate map $i$.
If $w:N\to N$ is a linear automorphism of the hyperplane
arrangement $Y^\eta$ then $E^\eta(\phi^w)=E^\eta(\phi)$.
\end{prop}
\begin{proof}
By equation (\ref{eq:zsingexp}) it is clear that
$E^{\eta,i,p}(\phi_p)$ is independent of analytic continuation
of $\phi_p$ along paths inside
$i(\mathbb{D}^\times\times\mathbb{D}^{n-1})$.
Suppose that $i,i^\prime$ both satisfy the requirements above,
and suppose that $i(\{0\}\times\mathbb{D}^{n-1})\cap
i^\prime(\{0\}\times\mathbb{D}^{n-1})\not=\emptyset$. Let
$V\subset i(\{0\}\times\mathbb{D}^{n-1})\cap
i^\prime(\{0\}\times\mathbb{D}^{n-1})\subset E$
be a connected contractible open set
(where $E$ denotes the exceptional divisor).
By the properties of $i$ and $i^\prime$ we have
$(i^\prime)^{-1}(i(\e,z^\prime))=(\e^\prime,w^\prime)$
with $\e^\prime(\e,z^\prime)=\e h(\e,z^\prime)$
where $h$ is holomorphic and nonzero on
$i^{-1}(V)$. If we plug this in the expansion (\ref{eq:zsingexp})
we see that the exponents defined by $i$ and by $i^\prime$
are equal if we use branches of $\phi$ on the image of
$i$ and on the image of $i^\prime$ which are related by
analytic continuation via the connected component of
the intersection of the images of $i$ and $i^\prime$ which
contains $V$. By AR we see that
that any path in $N^{\operatorname{reg}}$
is homotopic to a path which is contained in a finite
union of coordinate patches of the form
$i(\mathbb{D}^\times\times\mathbb{D}^{n-1})$.
With the above this
shows at once that the set of exponents does not depend
on the choice of the coordinate map $i$ and is independent
for analytic continuation of $\phi_p$ within
$N^{\operatorname{reg}}$.
For the last assertion
we remark that $i_w=w\circ i$ is also a standard coordinate map,
hence $E^\eta(\phi^w)=
E^{\eta,i_w,w(p)}((\phi^{w})^{w(p)})=E^{\eta,i,p}(\phi_p)=E^\eta(\phi)$.
\end{proof}
What lies behind this notion of exponents is the well
known ``decone construction'' on a central hyperplane arrangement.
This elementary construction implies that if $Y^\eta$ is nonempty then
$N^{\operatorname{reg}}$ is a isomorphic to a product
\begin{equation}
N^{\operatorname{reg}}\simeq\mathbb{D}^\times\times E^{\operatorname{reg}}
\end{equation}
Indeed, the restriction of the Hopf fibration
$p:\C^n\backslash\{0\}\to E=\mathbb{P}(\C^n)$ to the complement
of one of the hyperplanes $H$ of $Y$ is a trivial fibration since
$E\backslash\mathbb{P}(H)\simeq \C^{n-1}$ is contractible. Hence
the further restriction of this fibration to
$N^{\operatorname{reg}}$ is a fortiori trivial.
Thus we have a decomposition
\begin{equation}\label{eq:funddec}
\pi_1(N^{\operatorname{reg}},p)\simeq\Z\times
\pi_1(E^{\operatorname{reg}},[p])
\end{equation}
Now let $\L\subset\O(N^{\operatorname{reg}})$ be a local system
of finite rank $r$ of germs of Nilsson class function on
$N^{\operatorname{reg}}$.
We remark that, as a result of LN1, the germs of any local
Nilsson class function $\phi$ on $N^{\operatorname{reg}}$ are contained
in such a local system.
\begin{dfn}\label{dfn:gamma}
We denote by $T_p^\eta$ the monodromy map on $\L_p\simeq\C^r$
which corresponds to analytic continuation along the loop
$\g_p^\eta:t\to \exp(2i\pi t)p$.
Observe that $[\g_p^\eta]$ is a generator
of $\Z$ in (\ref{eq:funddec}).
\end{dfn}
In view of
(\ref{eq:funddec}) we may use $T_p^\eta$ to split the sheaf
$\L$ as a direct sum
\begin{equation}\label{eq:block}
\L=\bigoplus_{t\in\C^\times}\L^\eta(t)
\end{equation}
of generalized eigensheaves of $T^\eta$. In view of
(\ref{eq:zsingexp}) it is clear that if $\phi\not=0$
then $\phi\in\L_p^\eta(t)$ iff
$E^\eta(\phi)=\{s\}$ for some exponent $s\in\C$ such
that $t=\exp(2i\pi s)$.

Given $s\in\C\cup\{\infty\}$ we define a subsheaf
$F_s^\eta(\L)\subset \L^\eta(\exp(2i\pi s))$ of $\L$ by setting for each
$p\in N^{\operatorname{reg}}$:
\begin{equation}
F_s^\eta(\L)_p=\{\phi\in\L_p\mid E^\eta(\phi)
=\{\kappa\}\mathrm{\ with\ }
\kappa-s\in\mathbb{Z}_{\geq 0}\cup\{\infty\}\}
\end{equation}
One checks easily that this is a linear subspace of $\L_p$.
By Proposition \ref{prop:indep} it
is invariant for the parallel transport in the local system $\L$,
hence it defines a subsheaf. Moreover, these subsheaves of $\L$
are invariant for the action of the group of automorphisms of $\L$
which are induced by linear automorphisms of the arrangement
$Y^\eta$. For each $t\in\C^\times$ the subsheaves $F^\eta_s(\L)\subset\L$
with $s\in\C$ such that $\exp(2i\pi s)=t$
define a descending filtration
\begin{equation}
\dots\supset F^\eta_s(\L)\supset F^\eta_{s+1}(\L)\supset\dots
\end{equation}
of the direct summand $\L^\eta(t)$ of $\L$.
\begin{dfn}
We define a local system $\operatorname{Gr}^\eta(\L)$ by
\begin{equation}
\operatorname{Gr}^\eta(\L)=\bigoplus_{s\in\C}
\operatorname{Gr}^\eta_s(\L),\mathrm{\ with\ }
\operatorname{Gr}^\eta_s(\L)=
F^\eta_s(\L)/F^\eta_{s+1}(\L)
\end{equation}
\end{dfn}
For each $s\in\C$ we define the
multiplicity $\operatorname{mult}^\eta(\L,s)$ of $s$ as an exponent
at $\eta$ of the system of $\L$ by
\begin{equation}
\operatorname{mult}^\eta(\L,s)=\dim(\operatorname{Gr}^\eta_s(\L))
\end{equation}
\begin{dfn}
The (multi-)set $E^\eta\subset\C$ of exponents of $\L$ at $\eta$ are the complex
numbers $s\in\C$ such that $\operatorname{mult}^\eta(\L,s)>0$.
\end{dfn}
\begin{cor}\label{cor:expT}
The (multi-)set $\exp(2i\pi E^\eta)\subset\C^\times$
is the generalized eigenvalue spectrum of $T^\eta$ acting on $\L$.
\end{cor}
\begin{ex}\label{ex:bench}
Consider for $\mu\in\af_\C^*$ the sheaf $\L$ of local solutions
of the set of equations
\begin{equation}
\partial(p)\phi=p(\mu)\phi,\ \forall p\in\C[\af_\C^*]^W
\end{equation}
and let $\eta\in\af_\C$ be any point. Any local solution $\phi$
is holomorphic at $\eta$ and is completely determined by its
harmonic derivatives $\partial(q)(\phi)(\eta)$ at $\eta$.
Hence the set of exponents of $\L$ at $\eta$ is independent
of $\eta$ and $\mu$, and is equal to the set
$0,1,\dots,|\Sigma^l_+|$ where $\operatorname{mult}^\eta(\L,s)
=\dim{\operatorname{Harm}_s}(W)$, the dimension of the space of
$W$-harmonic polynomials of homogeneous degree $s$.
\end{ex}
\begin{ex}
Consider the sheaf $\L$ of local solutions of (\ref{eq:cpx}).
Suppose that $\eta\in\af_\C^{\operatorname{reg}}$ is a regular point.
Again a local solution $\phi$ of (\ref{eq:cpx}) near $\eta$ is
holomorphic at $\eta$ and is completely determined by its harmonic
derivatives $\partial(q)(\phi)(\eta)$ at $\eta$. Hence the answer
is the same as in the previous example.
\end{ex}
\begin{ex}\label{ex:cpx}
Let $\L$ be as in the previous example, but now we take
$\eta=i\pi\omega_j/k_j$ as in subsection \ref{sub:distbdy}.
The exponents of $\L$ at $\eta$ are equal to
$-|\Sigma^a_{\eta,+}|,\dots,|\Sigma^l_+|-|\Sigma^a_{\eta,+}|$, and
if $\phi_\mu$ denotes the holomorphic extension to $\af+i\pi\Omega$
of the spherical function, then for generic $\mu$ we have
$E^{\eta}(\phi_\mu)=\{|\Sigma_{\eta,+}|-|\Sigma^a_{\eta,+}|\}$
(see Proposition \ref{cpx}).
\end{ex}
\subsection{Harish-Chandra's radial system of differential equations}
In this subsection we describe the system of differential equation we
are mainly interested in, the radial differential equations
for an elementary spherical functions $\phi_\mu^X$ on a Riemannian
symmetric space $X=G/K$ restricted to a maximal flat, totally geodesic
subspace $A_X=A.x_0\subset X$.

The elementary spherical function $\phi_\mu^X$ (with
$\mu\in\af_\C^*$) on $X=G/K$ is a $K$-invariant solution of the
$G$-invariant system of differential equations
\begin{equation}
(\Delta-\gamma_X(\Delta)(\mu))\phi=0\ \forall\Delta\in\mathbf{D}(X)
\end{equation}
where $\mathbf{D}(X)$ denotes the ring of $G$-invariant
differential operators on $X$, and where
$\gamma_X:\mathbf{D}(X)\to\mathbb{C}[\af^*]^{W}$ is the
Harish-Chandra isomorphism.

By separation of variables we see that the
restriction of $\phi^X_\mu$ to $A_X$ is a $W$-invariant solution of
the system of differential equations
\begin{equation}\label{eq:sphsys}
(D-\gamma_X(D)(\mu))\phi=0\ \forall D\in\Ri_X
\end{equation}
on $A_X$ or on its complexification $A_{X,\C}=A_\C.x_0$, where
$\Ri_X\simeq \mathbf{D}(X)$ is the algebra of radial parts of the
operators $\Delta\in\mathbf{D}(X)$. Notice that we use the same
notation $\g_X$ for the Harish-Chandra isomorphism defined on
$\Ri_X$.

Let $T_X=T.x_0\subset A_{X,\C}$ be the compact form of $A_{X,\C}$.
It is a maximal flat totally geodesic subspace of a compact dual
symmetric space $U/K$ (which is by our choices simply connected).
The restrictions to $T_X$ of the zonal spherical functions of
$U/K$ are $W$-invariant simultaneous eigenfunctions of $\Ri_X$.
Since these zonal polynomials constitute a linear basis of the
space of $W$-invariant Laurent polynomials on $A_{X,\C}=A_\C/F$
this implies that the operators in $\Ri_X$ descend to
\emph{polynomial} differential operators on the complex affine
quotient space $W\backslash A_\C/F$.
\subsection{The hypergeometric system of differential equations}
In this subsection we describe a parameter deformation of the
Harish-Chandra system (\ref{eq:sphsys}) of differential equations
that we will use to study properties of solutions
of (\ref{eq:sphsys}). This parameter family of systems of differential
equations is called the system of hypergeometric equations associated with
root systems.
As was explained, this deformation is an essential ingredient for the
computation of the leading exponents (\ref{eq:leadexp}) of the spherical
functions at extremal points of $T_{\Omega}$.

We need to introduce some notations.
Let $\Sigma$ be a (not
necessarily reduced) irreducible root system in $\af^*$.
We consider indeterminates $\mathbf{m}_{\a}$ which are labeled by
the \emph{$W$-orbits} of the roots $\a\in\Sigma$ (in other words,
$\mathbf{m}_\a=\mathbf{m}_\b$ if $W\a=W\b$). Let
$\C[\mathbf{m}_\a]$ be the complex polynomial algebra over these
indeterminates $\mathbf{m}_\a$. If $X$ is a Riemannian symmetric
space with restricted root system isomorphic to $\Sigma$ then
$m_\a^X\in\mathbb{N}$ denote the root multiplicities of $X$.

The following result is one of the cornerstones of the theory of
hypergeometric functions for root systems.
\begin{thm}[\cite{Op0}, \cite{H}]
Let $\mathbb{A}$ denote the Weyl algebra of polynomial
differential operators on the complex affine space $W\backslash
A_\mathbb{C}/F\simeq \C^n$ with coefficients in the polynomial ring
$\C[\mathbf{m}_\a]$. There exists a unique subalgebra
$\Ri\subset\C[\mathbf{m}_\a]\otimes\mathbb{A}$ with the following
properties:

(1) The algebra $\Ri$ is isomorphic to the polynomial ring
$\C[{\mathbf{m}}_\a][\af^*_\C]^W$ via a Harish-Chandra isomorphism
$\g$ of algebras. This isomorphism $\g$ has the characterizing
property that any element $D\in\Ri$ is asymptotically equal to
the constant
coefficient operator $\g(D)(\cdot-\rho(\mathbf{m}))$ on
$A_\mathbb{C}$ (viewed as an element of the symmetric algebra on
$\af_\C$) along regular directions
towards infinity in $A_+$.

(2) If we specialize $\mathbf{m}$ at the multiplicity function
$m^X$ for a Riemannian symmetric space $X=G/K$ with restricted
root system $\Sigma_X\subset\af^*$ such that
$\Sigma^l=\Sigma^l_X$, such that $A_\C$ is the maximal torus of
$G_\C$, then $\Ri$ specializes to $\Ri_X$ and $\g$ to the
Harish-Chandra isomorphism $\g_X$.
\end{thm}
It is remarkable that the theory of Dunkl operators provides a
proof of this theorem which is both elementary and simple
\cite{H}.
\begin{prop}[\cite{HO1},\cite{Op4}, Remark 6.10]\label{prop:RS}
Let $\mathbb{D}^{\operatorname{reg}}$ be the ring of algebraic
differential
operators on the affine variety $A^{\operatorname{reg}}_\C/F$.
For each multiplicity parameter $m=(m_\a)$ and $\mu\in\af^*_\C$
let
$\mathcal{I}_{m,\mu}\subset\mathbb{D}^{\operatorname{reg}}$
denote the $W$-invariant left ideal
\begin{equation}\label{eq:ideal}
\mathcal{I}_{m,\mu}:=\sum_{D\in\Ri_m}
\mathbb{D}^{\operatorname{reg}}(D-\g_m(D)(\mu))
\end{equation}
Here
$\Ri_m$ is the specialization of $\Ri$ at $m$, and $\g_m$
the corresponding Harish-Chandra homomorphism.

Consider the $\mathbb{D}^{\operatorname{reg}}$-module
$\mathcal{M}_{\mu,m}=\mathbb{D}^{\operatorname{reg}}/
\mathcal{I}_{m,\mu}$
on $A_\C^{\operatorname{reg}}/F$.
Then in the terminology of Chapter IV, section 7 of
\cite{Bo}, $\mathcal{M}_{\mu,m}$
is an algebraic connection on $A_\C^{\operatorname{reg}}/F=
A_\C/F-\{\delta=0\}$ of rank $|W|$ which is regular.
Moreover $\mathcal{M}_{\mu,m}$ is $W$-equivariant.
\end{prop}
\begin{proof}
The elements of $\Ri_m$ are algebraic and the coefficients
are known to be regular on $A_\C^{\operatorname{reg}}/F$ (the simplest
way to see this is to use Dunkl-Cherednik operators \cite{H}, \cite{Op3}).
It is known that
$\mathcal{M}_{\mu,m}$ is $\mathcal{O}(A_\C^{\operatorname{reg}}/F)$-
free of rank $|W|$ by \cite{HO1}, and it is clear that
$\mathcal{M}_{\mu,m}$ is $W$-equivariant. It remains to prove the
regularity.

The elements of
$\Ri_m$ descend to the regular part of the adjoint torus
$A_\C^{\operatorname{adj,reg}}/F$ with character lattice
$Q=\mathbb{Z}\Sigma\subset\af^*$. We view this as an open subset
of the toric completion of $A_\C^{\operatorname{adj}}/F$ associated
with the decomposition of $\af^*$ in Weyl chambers. This is a
projective variety. It clearly suffices to prove the regularity
on $A_\C^{\operatorname{adj,reg}}/F$.

On $A_\C^{\operatorname{adj,reg}}/F$ one can explicitly rewrite the
module $\mathcal{M}_{\mu,m}$ as a connection of rank $|W|$ with
logarithmic singularities at infinity (see \cite{HO1}).
According to \cite{HO1} the connection matrix depends polynomially
on the parameters $\mu$ and $m$. It remains to
show that the connection is also regular singular at the
components of the discriminant locus $\delta=0$. Since the connection
depends polynomially on the parameters $\mu,m$ it is easy to see
that the set of parameters $\mu,m$ for which the connection is
regular singular is a Zariski-closed set. If $m=0$ the system is
trivially regular singular. If $\mu$ is sufficiently generic then
the theory of shift operators gives equivalences between the
modules $\mathcal{M}_{\mu,m}$ and $\mathcal{M}_{\mu,m^\prime}$ if
$m-m^\prime$ belong to the ``lattice of integral shifts'' (see
e.g. \cite{Op0} or \cite{HS}) in the space of multiplicity
parameters. The result follows.
\end{proof}
\begin{rem}\label{rem:cycl}
The element $u=1\in \M_{\mu,m}$ is a cyclic vector. Via $u$ the complex
vector space of $D$-module homomorphisms of $\M_{\mu,m}$ to $\mathcal{O}_p$
correspond to the space $\L_p(\mu,m)$ of solutions
in $\mathcal{O}_p$ of the $W$-invariant system of
differential equations
\begin{equation}\label{eq:hypsys}
(D-\gamma_m(D)(\mu))\phi=0\ \forall D\in\Ri_m
\end{equation}
on $A_\C^{\operatorname{reg}}/F$.
\end{rem}
\begin{cor}
The local system $\L(\mu,m)$ of germs of solutions of
(\ref{eq:hypsys}) is a local system of germs of Nilsson class functions
on $\af_\C^{\operatorname{reg}}$. Hence the results of Subsection
\ref{sub:expAR} are applicable to $\L(\mu,m)$.
\end{cor}
\begin{proof}
By \cite[Proposition 4.6.6]{Bj} it is sufficient to check the
moderate growth conditions for solutions of (\ref{eq:hypsys})
on the dense open set of subregular points of $\d=0$.
Since we can rewrite the system (\ref{eq:hypsys}) as a meromorphic
connection on $\af_\C^{\operatorname{reg}}$ which is regular
singular along $\d=0$ according to Proposition \ref{prop:RS}
this follows from \cite[Remark (5.9)]{Bo} (see also \cite{D}).
\end{proof}
Let $X$ be a Riemannian symmetric space with maximal flat geodesic subspace
$A.x_0$. The holomorphic extension of the restriction to $A.x_0$ of the
spherical function $\phi^X_\mu$ to $AT_\Omega^2 x_0$ is a
holomorphic $W$-invariant solution of (\ref{eq:sphsys}) on
$AT_\Omega^2.x_0$. This function is the
specialization of a holomorphic family (in the parameter $m$) of
solutions of (\ref{eq:hypsys}) by virtue of the following theorem:
\begin{thm}\label{thm:hypfun}(\cite{HO1},\cite{HS},\cite{Op4})
There exists an $\e>0$ such that
for all multiplicity parameters $m\in\Qc(-\e)$, the space of multiplicity
parameters such that $\operatorname{Re}(m_\a)\geq -\e\
\forall\a\in\Sigma$, the hypergeometric system (\ref{eq:hypsys})
has a unique solution $\phi_{\mu,m}$, the hypergeometric function,
which extends to a $W$-invariant and holomorphic function on
$AT_\Omega^2.x_0$. The function
$(t,\mu,m)\to\phi_{\mu,m}(t)$ is holomorphic on
$(AT_\Omega^2.x_0)\times\Qc(-\e)\times\af_\C^*$.
\end{thm}
Recall the covering map
$\pi:\af_\C\to A_\C/F\simeq A_\C.x_0$
of (\ref{eq:cov}) which is given by the exponential map
$\pi(X)=\exp(\pi i X)F$.
Via this map we will lift the differential
equations (\ref{eq:hypsys}) to $\af_\C^{\operatorname{reg}}$ and work
on $\af_\C$ rather than $A_\C/F$. On this space the system of
differential equations (\ref{eq:hypsys}) is invariant for the
action of the affine Weyl group $W^a=W\ltimes {Q^\vee}$.
In particular, we will work on the tube domain $i\af+\Omega\subset \af_\C$
instead of $AT_\Omega^2/F\subset A_\C/F$
(recall that the logarithm is well defined on
$AT_\Omega^2$). It is well known \cite{HO1} that the spherical
system of eigenfunction equations can be cast in the form of an
integrable connection on $\af_\C$ with singularities along the
collection of affine hyperplanes $\a(H)\in\Z$ (not $\in\pi
i \Z$ as in \cite{HO1}, since we have multiplied everything by
$(\pi i)^{-1}$).
\subsection{The indicial equation}\label{sub:ind}
We will show in this subsection that the
exponents of the hypergeometric equations (\ref{eq:hypsys})
at $\eta\in\af_\C$ coincide with the eigenvalues of the residue
matrix of a specially chosen integrable connection with
simple poles which is equivalent to (\ref{eq:hypsys}).
The characteristic equation of the residue matrix has coefficients
which are polynomials in the parameters $m_\a$. This equation is
called the indicial equation of (\ref{eq:hypsys}) at $\eta$.

Let us first construct an explicit standard coordinate map $i$
as used in the definition of the set of exponents.
Consider a parameterized line $x\to \eta+xV_1$ through $\eta$,
where $V_1$ is small and chosen in such a way that
this line is not contained in the union of the singular affine
hyperplanes. We choose coordinates $(z_1=\e,z_2,\dots z_n)$
(with $z_1=\e\in\mathbb{D}^\times$ and for $i>1$: $z_i\in \mathbb{D}$),
which we will often write as $z=(\e,z^\prime)\in\mathbb{D}^\times
\times\mathbb{D}^{n-1}$ with $z^\prime=(z_2,\dots,z_n)$.
First we choose $V_2,\dots,V_n$ in $\af$ such that $\Vert
V_i\Vert$ is small for all $i$, and such that
$(V_1,V_2,\dots,V_n)$ is a basis of the real vector space $\af$.
Then our coordinate map $i$ is given by
\begin{equation}
i(\e,z^\prime)=\eta+\e(V_1+\sum_{i\geq 2} z_iV_i)\in \af.
\end{equation}
If we lift this coordinate map to the blow-up of $\af_\C$ at
$\eta$ then the coordinates can be naturally extended to the
polydisk $\mathbb{D}^n$, and this is
then a coordinate neighborhood of a regular point of the
exceptional divisor $E$. The intersection of this neighborhood
with $E$ is described by the equation $z_1=0$. The complement of
$z_1=0$ in $\mathbb{D}^n$ is $\mathbb{D}^\times\times\mathbb{D}^{n-1}$,
the ``punctured polydisk''.
The Euler vector field $\E^\eta$ is given in these
coordinates by $z_1\partial/\partial{z_1}=\e\partial/\partial{\e}$.

Let $p$ be a point in the punctured polydisk
$\mathbb{D}^\times\times\mathbb{D}^{n-1}$
and let $\O_p$
denote the ring of holomorphic germs at $p$. Consider a subspace $U^*$
of dimension $|W|$ of the ring of holomorphic linear partial differential
operators on $\mathbb{D}^\times\times\mathbb{D}^{n-1}$
such that at all points $p\in
\mathbb{D}^\times\times\mathbb{D}^{n-1}$,
the free $\O_p$-module $\O_p\otimes U^*$ is a complement
for the left ideal $\mathcal{I}_{\mu,m}$.
We require further that the elements of $U^*$ commute with the Euler
vector field $\E^\eta$ (in other words, they are homogeneous of
degree $0$), and that $1\in U^*$. Such linear subspaces $U^*$ exist,
for instance one could take as a basis $b_i=
\e^{\operatorname{deg}(q_i)}\partial(q_i)$, where $q_i$ runs
over a homogeneous basis of $W$-harmonic polynomials
on $\af^*_\C$ with $q_1=1$.

We rewrite the differential equations (\ref{eq:hypsys})
(with $\mu\in\af_\C^*$)
in connection form with respect to the above basis $\{b_i\}$ and
coordinates $\{z_i\}$. We define matrices
$A^i_{\mu,m}\in\operatorname{End}_{\O_p}(\O_p\otimes U)$ (where $U$
denotes the dual of $U^*$, with dual basis $b_i^*$) which are
characterized by the requirement that
\begin{equation}\label{eq:A}
\frac{\partial}{\partial z_i}\circ b_k\in\sum_{j}
(A^{i}_{\mu,m})^{\mathrm{tr}}_{jk}b_j+\I_{\mu,m}.
\end{equation}
As an $\O_p$-module, the cyclic $D$-module $(M_{\mu,m},u)$ is equal
to $\O_p\otimes U^*u$, with basis $\overline{b_i}=b_i.u$.
Then the desired (flat) connection form of (\ref{eq:hypsys}) is
defined on the free $\O_p$-module $\O_p\otimes U$ by
\begin{equation}\label{con}
\frac{\partial \Phi}{\partial z_i}=A^i_{\mu,m}\Phi\ \ (\Phi\in\O_p\otimes U).
\end{equation}
By construction, if $\phi$ is a solution of (\ref{eq:hypsys}) then
\begin{equation}\label{eq:iso}
\Phi(\phi):=\sum_i b_i(\phi)b_i^*
\end{equation}
is a solution vector of (\ref{con}).
Conversely, if $\Phi$ is a solution vector of (\ref{con}) then the
first coordinate $\phi=\langle b_1,\Phi\rangle$ is a solution of
(\ref{eq:hypsys}). Since the linear map $\phi\to\Phi=\sum_i b_i(\phi)b_i^*$
is clearly injective we see by a dimension count that
these linear maps are inverse isomorphisms between the
solution spaces of these two systems of differential equation.
\begin{rem}\label{rem:loc}
Since the local solution space of an integrable connection at a
regular point $p$ can be identified with the fiber of the
underlying vector bundle at $p$, the above gives an isomorphism
(depending on $p$) between the local solution space
$\L_p(\mu,m)$ of (\ref{eq:hypsys}) at $p$ and the complex vector
space $U$.
\end{rem}
We claim that the system (\ref{con}) has simple singularities at
$\e=0$. The basis vectors $b_i=\e^{\operatorname{deg}(q_i)}\partial(q_i)$
have homogeneous degree zero and thus belong to the ring
$\mathcal{D}_0$ of holomorphic differential operators
on $\mathbb{D}^{n}$ generated by vector fields tangent to $\e=0$
(i.e. by $\partial/\partial{z_i}$ with $i>1$ and by
$\e\partial/\partial\e$). Therefore our claim is easily implied by
(also compare to \cite[Proposition 3.2]{HO1}):
\begin{lem}
Given $B\in\mathcal{D}_0$ there exists a unique section
$u(B)_{\mu,m}=\sum_{j}u(B)^j_{\mu,m}b_j\in
\mathcal{O}(\mathbb{D}^n)\otimes U^*$ such that
\begin{equation}\label{eq:Bmat}
B\in u(B)_{\mu,m}+\I_{\mu,m}
\end{equation}
The map $\mathcal{D}_0\ni B\to u(B)_{\mu,m}$ is an
$\mathcal{O}(\mathbb{D}^n)$-module morphism
which depends polynomially on $\mu$ and $m$.
For all $B\in\mathcal{D}_0$,
$u(B)_{\mu,m}|_{\{0\}\times \mathbb{D}^{n-1}}$ is
independent of $\mu$.
\end{lem}
\begin{proof}
We use induction on the order $d$ of $B$.
Using the well known theorem that $\C[\af_\C^*]$ is the
free $\C[\af_\C^*]^W$-module generated by the $W$-harmonic
polynomials, we have a unique decomposition
\begin{equation}
B=\sum_{i,j} f_{i,j}(\e,z^\prime)b_i \e^{d_{i,j}}\partial(p_{i,j})
\end{equation}
with $p_{i,j}\in\C[\af_\C]^W$ a homogeneous polynomial of
degree $d_{i,j}$ such that $d_{i,j}+
\operatorname{deg}(b_i)\leq d$,
and where $f_{i,j}(\e,z^\prime)$ is holomorphic
for all $i,j$.
Now $\e^{d_{i,j}}\partial(p_{i,j})=
\e^{d_{i,j}}(D_{p_{i,j}}-\g(D_{p_{i,j}})(\mu))$ modulo lower order
operators in $\mathcal{D}_0$, where we have used the fact that
for $p\in\C[\af_\C]^W$ homogeneous,
the lowest homogeneous part $h^\eta(D_{p})$ at $\eta$ of
$D_{p}\in\Ri_m$ contains the highest order term $\partial(p)$
of $D_p$ (see \cite{Op0}). By the induction hypothesis we conclude
the existence $u(B)_{\mu,m}$. Using the independence of the
$W$-harmonic polynomials over the ring $\C[\af_\C^*]^W$ and the
induction hypothesis the uniqueness of $u(B)_{\mu,m}$ follows too.
By induction and using the fact that the operators $D_p$ depend polynomially
on $\mu$ and $m$ we conclude that $u(B)_{\mu,m}$
is holomorphic on $\mathbb{D}^n$ and polynomial in $\mu$ and $m$.
Since in the induction step $\mu$ only occurs via the terms of the
form $\e^{d_{i,j}}\g(D_{p_{i,j}})(\mu)$ we see that $\mu$ does not
influence the evaluation at $\e=0$ of $u(B)_{\mu,m}$.
\end{proof}
Let $R_m$ be the residue matrix of $A^1_{\mu,m}$ at $z_1=0$.
By the previous lemma $R_m$ is independent of $\mu$ and is polynomial
in $m$. As is well known (cf. \cite{Bo}, Chapter IV, section 4, or
\cite{D}) $R_m$ is
independent of the coordinate map $i$. Moreover, let us consider
on $V=i(\{0\}\times\mathbb{D}^{n-1})$ the integrable connection defined
by the restrictions $B^i_{\mu,m}:={A^i_{\mu,m}}|_{V}$ for $i>1$.
Then the residue $R_m$ is known to be flat for this integrable
connection on $V$. In particular, its
characteristic equation is independent of $z^\prime$.
\begin{thm}\label{lem:eigen}
The exponents of (\ref{eq:hypsys}) at $\eta$ are the eigenvalues
of the residue matrix $R_m$ of $A^1_{\mu,m}$ at $z_1=0$. The
characteristic polynomial of $R_m$ is independent of $\mu$ and of
$z^\prime$ and has polynomial coefficients in the $m_\a$.
\end{thm}
\begin{proof}
By changing the basis of the trivial vector bundle
(with fiber $U$) on $i(\mathbb{D}^n)$ by a suitable invertible
matrix depending on $z^\prime$ only we may assume that
$B^i_m=0$ on $V$ for all $i>1$. We denote the finite dimensional
complex vector space of sections spanned by this basis of flat
sections $\mathcal{U}$.
By the flatness of $R_m$ for the restricted connection on $V$ as above,
$R_m$ is constant in this new basis (i.e. independent of $z^\prime$).
Let $s$ be an eigenvalue of $R_m$, and let $v$ be a generalized
$R_m$-eigenvector with eigenvalue $s$. Put
$u(\e)=\exp(\log(\e)R_m)v=\e^s\exp(\log(\e)(R_m-s\operatorname{Id}))v$,
and observe that
\begin{equation}\label{eq:q0}
q_0^{^{(s,v)}}(\e,z^\prime):=\exp(\log(\e)(R_m-s\operatorname{Id}_\mathcal{U}))v
\end{equation}
is a $\mathcal{U}$-valued polynomial in $\log(\e)$.
We denote the series expansion of $\e A^1_{\mu,m}$ in $\e$
with respect to a fixed basis of $\mathcal{U}$ by
\begin{equation}
\e A^1_{\mu,m}(\e,z^\prime)=R_m+\sum_{k>1}\e^k
A^1_{\mu,m,k}(z^\prime)
\end{equation}
with $A_{\mu,m,k}^1(z^\prime)$ holomorphic for
$z^\prime\in\mathbb{D}^{n-1}$.
Now we use the following relative version of
\cite[Ch. IV, \S 24, Hilfssatz XI]{Wa}: If
$q_i(\e,z^\prime)$ ($i<k$) are $\mathcal{U}$-valued polynomials in
$\log(\e)$ of degree $\leq N$ with coefficients in the
ring of holomorphic functions on $\mathbb{D}^{n-1}$,
then the equation
\begin{equation}\label{eq:rec}
\e\frac{\partial{q_k}}{\partial{\e}}+((s+k)\operatorname{Id}_\mathcal{U}-R_m)q_k=
\sum_{i=0}^{k-1}A_{\mu,m,k-i}^1(z^\prime)q_i
\end{equation}
has at least one solution $q_k$ which is polynomial in $\log(\e)$
and has coefficients in the ring of holomorphic functions in
$z^\prime\in\mathbb{D}^{n-1}$. The solution $q_k$ is unique and
has degree $\leq N$ if $(s+k)$ is not an eigenvalue of $R_m$. In
general there exist several solutions $q_k$ which are polynomial in
$\log(\e)$ and these solutions are all of degree $\leq N+r$
in $\log(\e)$, where $r$ is the maximal length of a Jordan block
of $R_m$ with eigenvalue $(s+k)$.

Given a set $\{q_k^{(s,v)}\}$ of solutions of the recurrence relations
(\ref{eq:rec}) (with $q_k^{(s,v)}$ polynomial in $\log(\e)$ for
all $k$, and $q_0^{(s,v)}$ given by (\ref{eq:q0}))
there exists a convergent (but multivalued) series
solution $\Phi^{(s,v)}$ of (\ref{con}) on
$i(\mathbb{D}^\times\times\mathbb{D}^{n-1})$ of the form
\begin{equation}\label{eq:series}
\Phi^{(s,v)}(\e,z^\prime)=\e^s\sum_{k\geq 0}\e^k q_k^{(s,v)}(\e,z^\prime)
\end{equation}
(see e.g. \cite[Ch. IV, \S 24, XII]{Wa}). Notice that the degree of
$q_k^{(s,v)}(\e,z^\prime)$ ($k\geq 0$) as a polynomial in $\log(\e)$
with coefficients in the ring of holomorphic functions in
$z^\prime\in\mathbb{D}^{n-1}$ is
uniformly bounded.

Such series expansion is not necessarily unique, but by
choosing such a series solutions $\Phi^{(s,v)}$ for a set of pairs $(s,v)$
where $s$ runs through the set of eigenvalues of $R_m$ and for
each $s$, $v$ runs through a basis of the generalized $s$-eigenspace
of $R_m$ then the collection of multivalued solutions $\Phi^{(s,v)}$
on $i(\mathbb{D}^\times\times\mathbb{D}^{n-1})$ constitutes a
basis for the space of multivalued solutions of (\ref{con}).
On the other hand we have seen above that the flat sections on
$i(\mathbb{D}^\times\times\mathbb{D}^{n-1})$
all are of the form $\Phi=\sum_i b_i(\phi)b_i^*$ where
$\phi=\langle b_1,\Phi\rangle$ is a solution of (\ref{eq:hypsys}).
Hence the set of exponents of (\ref{eq:hypsys}) must coincide with
the set of eigenvalues of $R_m$, counted with multiplicity.
\end{proof}
As a result of the above theorem the following definition makes
sense.
\begin{dfn} Let $R_{\bf{m}}=R_{\bf{m}}^\eta$ denote the
$|W|\times|W|$-matrix
with coefficients in the ring
$\C[{\bf{m}}_\a]\otimes\O(\mathbb{D}^{n-1})$ such that
$R_m=R^\eta_m$ is the specialization of $R_{\bf{m}}^\eta$ at ${\bf{m}}=m$
(this matrix depends on the coordinate map $i$).
We call the characteristic polynomial $I_{\bf{m}}^\eta\in\C[{\bf{m}_\a}][X]$
of $R_{\bf{m}}^\eta$ the ``indicial polynomial'' of (\ref{eq:hypsys}) at $\eta$.
\end{dfn}
\begin{cor}\label{cor:ind}
The (multi-)set $E^\eta$ of exponents of (\ref{eq:hypsys}) at $\eta$
is equal to the (multi-)set of roots of the indicial
polynomial $I_m^\eta$ of (\ref{eq:hypsys}) at $\eta$.
\end{cor}
\subsection{Hecke algebras and exponents}
We now bring into play well known results on the monodromy of the
system of hypergeometric differential equations. We have quite  good
control for generic parameters as a consequence of the main
result, the fact that this representation of the affine braid group
factors through an affine Hecke algebra. We apply these results to
prove that the indicial polynomial $I^\eta$ at $\eta$ factors completely
over the ring of rational polynomials in the indeterminates $\bf{m}_\a$
with roots that are affine linear functions in the $\bf{m}_\a$
with half integral coefficients.

By affine Weyl group symmetry we may assume without loss
of generality that $\eta\in \overline{\Omega}\cap C$,
the fundamental alcove. From now on we will make this assumption.

By Corollary \ref{cor:expT} the generalized eigenvalue spectrum of
$\L_p(\mu,m)$ under the action of $T^\eta_p$ contains information
on the set $E^\eta$ of exponents of (\ref{eq:hypsys}). Since by
\ref{eq:funddec} $T^\eta_p$ is certainly central in
$\Pi_1(N^{\operatorname{reg}},p)$, the decomposition of $\L_p(\mu,m)$
in indecomposable blocks for the monodromy action of
$\Pi_1(N^{\operatorname{reg}},p)$ on $\L_p(\mu,m)$ refines the
decomposition in generalized $T^\eta_p$ eigenspaces
(by virtue of Schur's Lemma).

Therefore we now recall some fundamental facts on the monodromy representation of
the fundamental group $\Pi_1(W^a\backslash\af_\C^{\operatorname{reg}},p)$
(at a regular base point $p\in\af_\C^{\operatorname{reg}}$ in the fundamental
alcove $\overline{\Omega}\cap C$) on the local solution space
$\L_p=\L_p(\mu,m)$ of (\ref{eq:hypsys}).
By a well known result of Looijenga and Van der Lek (\cite{L},
also see \cite{HO1}, \cite{HS}, \cite{Op4}) the group
$\Pi_1(W^a\backslash\af_\C^{\operatorname{reg}},p)$
is isomorphic to the affine braid group $B^a$ of
$W^a=W\ltimes Q(\Sigma^\vee)$, the affine Weyl group of the affine root
system $\Sigma^a=\Sigma^l\times\mathbb{Z}$. In order to formulate
the result we need to define an affine root multiplicity function
$m^a$ on the affine roots in $\Sigma^a$ as follows. For the affine
simple roots $a_0=1-\theta,a_1=\a_1,\dots,a_n=\a_n$ we define
\begin{align}\label{eq:maff}
m_{a_0}^a&=m_\theta\\ \nonumber m_{a_i}^a&=m_{\a_i}+m_{\a_i/2}
\end{align}
and then we extend this to $\Sigma^a$ by $W^a$-invariance.
\begin{thm}(cf. \cite{HO1}, \cite{HS}, \cite{Op4})
The monodromy action on the $W^a$-equivariant local system
$\L_p(\mu,m)$ on $\af_\C^{\operatorname{reg}}$ factors
through an affine Hecke algebra $H(W^a,q^a)$ in the following sense.

Let $q^a$ be the label function \emph{on the affine root system}
$\Sigma^a=\Sigma^l\times\mathbb{Z}$ defined by
$q_b^a=\exp(-i\pi(m^a_b)$ for all $b\in\Sigma^a$. For the simple
affine roots $a_i$ we write $q^a_{a_i}:=q^a_i$.
The monodromy matrices $M_{\mu,m}(b_i)$ ($i=0,\dots,n$) of the
generators $b_i$ of $B^a$ satisfy
$(M_{\mu,m}(b_i)-1)(M_{\mu,m}(b_i)+q^a_i)=0$.
The monodromy representation $M_{\mu,m}$ of
$\Pi_1(W^a\backslash\af_\C^{\operatorname{reg}},p)$
depends analytically on the
parameters $m$ and $\mu$.
\end{thm}
Recall that $W_\eta^a$ is the isotropy
subgroup of $\eta$ in $W^a$, which is a finite reflection group,
and let $\Sigma_\eta^a$ be the corresponding root system.
There is a natural monomorphism $W_\eta^a\to W$ with
image $\tilde{W}^a_\eta\subset W$. We put $N_\eta=[W:\tilde{W}_\eta^a]$
for the index of this subgroup.

Let us denote by $B_\eta^a\subset B^a$ the braid group of
$W_\eta^a$, which we can identify, by Brieskorn's theorem on the
fundamental group of the regular orbit space of a finite
reflection group, with the fundamental group of the ``local
regular orbit space'' at $\eta$, namely
$\Pi_1(W^a_\eta\backslash N^{\operatorname{reg}},p)$.

Let $m_\eta^a$ be the restriction of $m^a$ to $\Sigma_\eta^a$, and
let $q_\eta^a$ be corresponding the corresponding root
multiplicity function on $\Sigma_\eta^a$. Let $\Qc$ denote the
finite dimensional complex vector space of complex multiplicity
functions $m$ on (the possibly non-reduced) root system $\Sigma$.
In a dense, open set
$\Qc_\eta^{\operatorname{reg}}\subset\Qc$ of values of the parameter $m$,
the finite dimensional Hecke algebra $H(W_\eta^a,q_\eta^a)$
(with $q_\eta^a=q(m_\eta^a)$) is a semisimple algebra.
If we assume that $m\in \Qc_\eta^{\operatorname{reg}}$ then, by Tits'
deformation lemma, we can index its set of irreducible modules
by $\widehat{W_\eta^a}$, the set
of irreducible representations of $W_\eta^a$. Given $\tau\in
\widehat{W_\eta^a}$ and $m\in \Qc_\eta^{\operatorname{reg}}$ we will
write $\pi_\tau^\eta(q_\eta^a)$ for the corresponding irreducible
$H(W_\eta^a,q_\eta^a)$-module. Upon restriction of the
monodromy action of $B^a$ on $\L_p(\mu,m)$ to $B_\eta^a$
we have:
\begin{cor}\label{cor:monfin}
Let $q=q(m)$ and $q_\eta^a=q(m_\eta^a)$ for $m\in
\Qc_\eta^{\operatorname{reg}}$. The monodromy action of
$\Pi_1(W^a_\eta\backslash N^{\operatorname{reg}},p)$ on
$\L_p(\mu,m)$ factors through the semisimple finite type Hecke algebra
$H(W_\eta^a,q_\eta^a)$, and
the local solution space $\L_p(\mu,m)$ decomposes under this
action in isotypical components
\begin{equation}\label{eq:mondec}
\L_p(\mu,m)=\bigoplus_{\tau\in\widehat{W_\eta^a}}
\L_p(\mu,m)(\tau)
\end{equation}
such that for each $\tau\in\widehat{W_\eta^a}$,
$\L_p(\mu,m)(\tau)\simeq K(\tau,m)\otimes\pi_\tau^\eta(q_\eta^a)$
with $\dim(K(\tau,m))=N_\eta\operatorname{deg}_\tau$
(independent of $m\in\Qc_\eta^{\operatorname{reg}}$).
\end{cor}
\begin{proof}
Using the rigidity of semisimple finite dimensional
algebras (Tits' deformation lemma, \cite{Ca}, Proposition 10.11.4)
the multiplicity of $\pi_\tau^\eta(q_\eta^a)$ is constant
in $(\mu,m)\in\af_\C^*\times\Qc_\eta^{\operatorname{reg}}$.
We may therefore compute the multiplicity by evaluating
at $(\mu,m)=(0,0)$. Hence it is equal to the multiplicity of
$\tau$ in the restriction of the regular representation of $W$ to
$\tilde{W}_\eta^a$, which is $N_\eta\operatorname{deg}_\tau$.
\end{proof}
The following topological observation due to Deligne
\cite{D1} is crucial for our purpose:
\begin{lem}\label{lem:delcent}
Let $\b^\eta_p\in B_\eta^a$ denote the local braid
in $\Pi_1(W^a_\eta\backslash N^{\operatorname{reg}},p)$
which corresponds to a reduced expression of the longest element
of $W_\eta^a$. Then $(\b^\eta_p)^2=[\g^\eta_p]$
(see Definition \ref{dfn:gamma}). In particular this
element is central in $B_\eta^a$.
\end{lem}
Given $\tau\in\widehat{W_\eta^a}$ we denote by
$p^i_{\eta,\tau}$ (with $i=1,\dots,N_\eta\operatorname{deg}_\tau$)
the embedding degrees of $\tau$ in the graded vector space of
$W$-harmonic polynomials.
We choose these $W$-harmonic embedding degrees so that
$i\to p^i_{\eta,\tau}$ is a non-decreasing sequence.
In particular, $p^1_{\eta,\tau}$ is the ``harmonic birthday''
of $\tau$ in the $W$-harmonic polynomials.
\begin{thm}\label{prop:ev}
Let $\tau\in\widehat{W_\eta^a}$
and let $m\in\Qc_\eta^{\operatorname{reg}}$.
The multiset $E^\eta(\tau,m)$ of exponents of
$\L_p(\mu,m)(\tau)$ consists of the complex numbers
\begin{equation}
s^i_{\eta,\tau}(m)=
p^i_{\eta,\tau}-\frac{1}{2}c_{\eta,\tau}(m),
\end{equation}
where $i$ runs from $1$ to $N_\eta\operatorname{deg}_\tau$,
each $s^i_{\eta,\tau}(m)$
occurring with multiplicity $\operatorname{deg}_\tau$.
Here $c_{\eta,\tau}(m)$ is the affine linear function of the
multiplicity parameters $m_\a$ with nonnegative integral
coefficients defined by (cf. (\ref{eq:maff}) for the
definition of $m^a_{\eta,b}$):
\begin{equation}
c_{\eta,\tau}(m)=\sum_{b\in\Sigma_{\eta,+}^a}
(1-\frac{\chi_\tau(s_b)}{\operatorname{deg}_\tau})m_{\eta,b}^a
\end{equation}
\end{thm}
\begin{proof}
Since $T_p^\eta$ is the monodromy action of the (locally) central
braid $(b^\eta_p)^2$ (by Lemma \ref{lem:delcent}) we see that $T_p^\eta$
acts trivially in the multiplicity space $K(\tau,m)$ and acts by
scalar multiplication in the irreducible representation
$\pi^\eta_\tau(q_\eta^a)$ of the Hecke algebra $H(W^a_\eta,q^a_\eta)$
by some scalar $C$. This $C$ is an element of the ring of Laurent polynomial
in the Hecke algebra labels $(q_{\eta,b}^a)^{1/2}$ (with $b\in\Sigma^a_\eta$)
since this is the splitting ring of the Hecke algebra.
By taking the determinant of $T_p^\eta$ in
$\pi^\eta_\tau(q_\eta^a)$ we find easily that
$C^{\operatorname{deg}_\tau}=
\exp(-i\pi\operatorname{deg}_\tau c_{\eta,\tau}(m))$.
This implies that $C$ is a root of $1$ times a monomial in the
$(q^a_{\eta,b})^{\pm 1/2}$. For $(q^a_{\eta,b})^{1/2}=1$ we have $C=1$, hence
\begin{equation}
C=\exp(-i\pi c(m))
\end{equation}
Let $\mathcal{N}$
denote the collection of functions $\nu$ on the set
$\tau\in\widehat{W^a_\eta}$ which associate to each
$\tau$ a finite multiset
$\nu(\tau)=\{\nu_{\tau,j}\mid j=1,\dots,N_\eta\operatorname{deg}_\tau^2\}$
of $N_\eta\operatorname{deg}_\tau^2$ integers $\nu_{\tau,j}\in\Z$.
By Corollary \ref{cor:expT} and Corollary \ref{cor:ind} it follows that
for each $m\in\Qc_\eta^{\operatorname{reg}}$ the set of roots of the
indicial polynomial $I^\eta_m$ is a multiset of the form
$\rho_{\tau,\nu,j}(m)=\nu_{\tau,j}-1/2c_{\eta,\tau}(m)$
for some $\nu\in\mathcal{N}$. For each $\nu\in\mathcal{N}$ the set
$\Qc(\nu)\subset\Qc$ of multiplicity parameters $m\in\Qc$ for which
the multiset of roots of $I^\eta_m$ is equal to the multiset
$\{\rho_{\tau,\nu,j}(m)\}$ is Zariski-closed
(since $I^\eta_m$ is a polynomial in $m$, by Corollary \ref{cor:ind}).
Moreover the union of these sets contains
$\Qc_\eta^{\operatorname{reg}}$. Since
$\mathcal{N}$ is countable, Baire's category
theorem implies that there must exist at least one $\nu_0\in\mathcal{N}$
such that the interior (in the analytic topology) of $\Qc(\nu_0)$ is
nonempty, and hence such that $\Qc=\Qc(\nu_0)$.

In the situation
$m\in\Qc_\eta^{\operatorname{reg}}$ we have the splitting in the
isotypical components (\ref{eq:mondec}). It follows that
the set $E^\eta(\tau,m)$ consists of the subset
$\rho_{\tau,\nu_0,j}(m)$
($j=1,\dots,N_\eta\operatorname{deg}_\tau^2$)
of roots of the indicial equation. Finally we need to
determine $\nu_0(\tau)$. This is resolved by taking
$m=0\in\Qc_\eta^{\operatorname{reg}}$ and comparing to
Example \ref{ex:bench}, after
making the additional remark that the monodromy representation
$\pi^\eta_\tau(q^a_\eta)$ is by definition equal to $\tau$
if $q^a_\eta=1$.
\end{proof}
\begin{cor}\label{cor:nolog}
For $m\in\Qc^{\operatorname{reg}}_\eta$ the action
of $T^\eta_p$ is semisimple. In particular, there are no
logarithmic terms in the decomposition (\ref{eq:zsingexp})
if $m\in\Qc^{\operatorname{reg}}_\eta$ and if
$\phi_p\in\L_p(\mu,m)$.
\end{cor}
So we conclude this subsection with the following remarkable result:
\begin{cor}
The indicial polynomial factorizes as
\begin{equation}
I_m(X)=\prod_{\tau\in\widehat{W^a_\eta}}
\prod_{i=1}^{N_\eta\operatorname{deg}_\tau}
(X-s^i_{\eta,\tau}(m))^{\operatorname{deg}_\tau}
\end{equation}
For $m\in\Qc^{\operatorname{reg}}$ this factorization
is compatible with the decomposition of $\L_p(\mu,m)$
in blocks of the form $\L_p(\mu,m)(\tau)$
as in Corollary \ref{cor:monfin}.
\end{cor}
\subsection{Computation of the leading exponents}
Let $\phi_{\mu,m}\in\L_p(\mu,m)$ with $p\in\overline{\Omega}\cap C$
denote the hypergeometric function, the solution of
(\ref{eq:hypsys}) whose germ at points of the fundamental
alcove $\overline{\Omega}\cap C$ we define by analytic
continuation along a path in $\overline{\Omega}\cap C$ of the
unique normalized $W$-invariant solution of (\ref{eq:hypsys})
which extends holomorphically to a neighborhood of $0\in \af_\C$.
\begin{cor}\label{cor:cont}
By definition, $\phi_{\mu,m}$ extends
holomorphically over all finite walls, the walls of $C$, to
a $W$-invariant function on $\Omega$, the interior of $WC$.
In particular, if $\eta\in
\overline{\Omega}\cap C$ and $\theta(\eta)\not=1$ then
$E^\eta(\phi_{\mu,m})=\{\kappa\}$ with $\kappa\in\Z_{\geq 0}$
(generically $\kappa=0$, of course).
\end{cor}
We will be interested in this section in the case where
$\eta=\omega_j/k_k$ as in Theorem \ref{th=omega}.
In this case we know that $\Sigma^a_\eta$ is an irreducible
root system. From the
definition of $\phi_{\mu,m}$ we see that
\begin{cor}\label{cor:trivind}
Let $m\in\Qc_\eta^{\operatorname{reg}}$.
Let $W_\eta\subset W_\eta^a$ be the maximal parabolic subgroup
of $W_\eta^a$ generated by the simple reflections $s_i$ of
$W$ which fix $\eta$. Then $\phi_{\mu,m}\in\L_p(\mu,m)$
belongs to the subspace $\L_p^\eta(\mu,m)$ defined by
\begin{equation}
\L_p^\eta(\mu,m):=
\bigoplus_{\tau\in J_\eta}\L_p(\mu,m)(\tau)
\end{equation}
where $J_\eta\subset\widehat{W^a_\eta}$ is the
subset consisting of irreducible representations which occur
in the induction of the trivial representation of $W_\eta$ to
$W_\eta^a$.
\end{cor}
The above fact restricts the $T^\eta_p$-spectrum of $\phi_{\mu,m}$,
and thus the set of exponents $E^\eta(\phi_{\mu,m})$, drastically
for $m\in\Qc^{\operatorname{reg}}$. We assume from now on that
$m$ is real valued, which we denote by
$m\in\Qc(\R)$. By Theorem \ref{prop:ev} the multiset
$E^\eta(\phi_{\mu,m})$ consists of real numbers now.
\begin{dfn} Let $m\in\Qc(\R)$.
We call the smallest element in the the multiset
$E^\eta(\phi_{\mu,m})$ the \emph{leading exponent} of
$\phi_{\mu,m}$ at $\eta$. The irreducible characters
$\tau\in J_\eta\subset W^a_\eta$ affording the leading exponent
are called \emph{leading characters}.
\end{dfn}
\begin{thm}
If $\Sigma^a_\eta$ is reduced and simply laced we denote the
root multiplicity by $m=m_1\geq 1$. In general $m_1$ denotes the
root multiplicity of the longest roots. The multiplicity of
half a long root is denoted by $m_{1/2}\geq 0$ (i.e. we consider
$C_n$ as the special case of $BC_n$ where $m_{1/2}=0$; since
the geometry of $\Omega$ depends on $\Sigma^l$ only this is
allowed).
Let $\eta\in\partial{\Omega}\cap C$ be an extremal boundary point
of $\Omega$ and assume that $m\in \Qc(\R)$
is in the cone $\mathcal{C}\in\Qc(\R)$ defined by the
inequalities
\begin{align}\label{eq:condm}
1&\leq m_1\leq m_2
\end{align}
(these inequalities are obviously satisfied by
the multiplicity function $m^X$ of a Riemannian symmetric
space $X$ with restricted root system $\Sigma_X$
such that $\Sigma_X^l=\Sigma^l$).
The leading exponent $s_\eta(m)$ of
$\phi_{\mu,m}$ at $\eta$ satisfies
\begin{equation}\label{eq:ineq}
s_\eta(m)\geq s^1_{\eta,\tau}(m)
\end{equation}
where
\begin{equation}
\tau=\sigma_\eta=\operatorname{det}_\eta^a\otimes
j_{W_\eta}^{W_\eta^a}(\operatorname{det}_\eta)
\end{equation}
Here $\operatorname{det}_\eta^a$ is the determinant representation
of $W_\eta^a$, and $\operatorname{det}_\eta$ its restriction to
$W_\eta$. Moreover, for generic $m\in\mathcal{C}$
the inequality (\ref{eq:ineq}) is an equality and
$\sigma_\eta$ is a
leading character.
\end{thm}
\begin{proof}
This is based on a case-by-case analysis.
We first assume that $m\in\Qc^{\operatorname{reg}}_\eta(\R)$ is
regular. We compute in all cases the
set $J_\eta$ of irreducible components $\tau$ of the induction
of the trivial representation of $W_\eta$ to $W^a_\eta$ (which is
relatively easy, as $W_\eta$ is a rather large subgroup of $W^a_\eta$).
In the classical cases we use the Littlewood-Richardson rule,
and in the exceptional cases we refer to the character tables in the
computer algebra packet CHEVIE. We use below the notations for
the irreducible characters as used in \cite{Ca}).

We have luck: if we consider for each $\pi\in J_\eta$
the smallest associated exponent $s^1_{\eta,\tau}(m)$
(using Theorem \ref{prop:ev}) we can simply check that these are
indeed all greater than or equal to $s^1_{\eta,\tau}(m)$, where
$\tau=\s_\eta$ and if $m\in\mathcal{C}$. Recall that
$\s_\eta\in J_\eta$ was the term which gave the unique leading
exponent in the complex case (see Example \ref{ex:cpx},
Proposition \ref{cpx}), which corresponds only to one
interior point $m^{X_\C}\in\mathcal{C}$
of $\mathcal{C}$. In any case, this surprising fact is enough to
prove that for generic $m\in\mathcal{C}$ the value $s^1_{\eta,\tau}(m)$
really is the leading exponent of $\phi_{\mu,m}$ at $\eta$ by the fact that
$\phi_{\mu,m}$ is holomorphic in the parameter $m$
(see Theorem \ref{thm:hypfun}).

Below will now show these claims in a case-by-case analysis:

\emph{Type $A_{l-1}(l\geq3)$:}
For $\Sigma=A_{l-1}$ all the nodes of the Dynkin diagram are
minuscule and thus extremal according to Theorem \ref{th=omega}.
Let $\omega_j$ be the $j$-th node of the
Dynkin diagram. By symmetry we may assume without loss of generality
that $1\leq j\leq l/2$. Recall that the irreducible characters $\chi_\l$
of $S_l$ are parameterized by the partitions $\l$ of $l$ in such a
way that $\chi_l=1$ and  $\chi_{1^l}=\e$ (the determinant representation).
We denote the i-th exponent corresponding to $\chi_{\l}$ by
$\s_\l^i(m)$.

By the Littlewood-Richardson rule \cite[Section I.9]{M} we have
\begin{equation}\label{eq:LR}
\operatorname{Ind}_{S_j\times S_{l-j}}^{S_l}(\chi_j\times \chi_{(l-j)})=
\bigoplus_{0\leq i\leq j}\chi_{(l-i,i)}
\end{equation}
and we have that (see \cite[Sections 11.2, 11.4]{Ca}):
\begin{equation}
\s_{\omega_j}:=\e\otimes
j_{S_{l-j}\times S_j}^{S_l}(\e_{l-j}\times\e_j)
=\chi_{(l-j,j)}
\end{equation}
Using Theorem \ref{prop:ev} and standard facts on
representations of $S^l$ we find that:
\begin{equation}\label{eq:formexp}
s_{(l-i,i)}^1(m)=i(1-(l+1-i)m/2)
\end{equation}
Under the condition (\ref{eq:condm}) (namely $m\geq 1$)
we see that among the exponents $s_{(l-i,i)}(m)$
at $\omega_j$ (thus with $i\leq j$) indeed
\begin{equation}
s_{\omega_j}(m):=s_{(l-j,j)}^1(m)=j(1-(l+1-j)m/2)
\end{equation}
is the unique minimal one, unless $l$ is even, $m=1$ and $j=l/2$.
In this last case the two components $(l/2,l/2)$
and $(l/2+1,l/2-1)$ of (\ref{eq:LR}) both have the same
exponent $l(l-2)/8$.

\emph{Type $B_l(l\geq3)$ ($\eta=\omega_1$):}
Recall that the irreducible characters $\chi_{(\l,\mu)}$ of
$B_l$ are parameterized by ordered pairs $(\lambda,\mu)$
of partitions of total weight $l$. Here $\chi_{(l,0)}=1$
and $\chi_{(0,1^l)}=\e$. We have (using the LR rule again for
wreath products, see \cite[I.Appendix B]{M}):
\begin{align}\label{eq:bn}
\operatorname{Ind}_{B_{l-1}}^{B_l}&(\chi_{(l-1,0)})\\
\nonumber&=\operatorname{Ind}_{B_{l-1}\times B_1}^{B_l}(\chi_{(l-1,0)}\times\chi_{(1,0)})+
\operatorname{Ind}_{B_{l-1}\times B_1}^{B_l}(\chi_{(l-1,0)}\times\chi_{(0,1)})=\\
\nonumber&=\chi_{(l,0)}+\chi_{(l-1,1)}+\chi_{((l-1,1),-)}
\end{align}
and thus
\begin{equation}
J_{\omega_1}=\{\chi_{(l,0)},\chi_{(l-1,1)},\chi_{((l-1,1),-)}\}
\end{equation}
{}From \cite[Proposition 11.4.2]{Ca} we find
\begin{equation}
\s_{\omega_1}:=
\e_l\otimes j_{B_{l-1}}^{B_l}(\chi_{\e_{l-1}})=
\e_l\otimes\chi_{(1,1^{l-1})}=\chi_{(l-1,1)}
\end{equation}
and the birthday of $\chi_{(1,1^{l-1})}$ is
$|\Sigma(B_{l-1})_+|=(l-1)^2$.

Using Theorem \ref{prop:ev} and standard results on
representations of $W(B_l)$ (e.g. \cite[Chapter 11]{Ca}) we
find
\begin{align}\label{eq:bnj}
s_{(l-1,1)}^1(m_2,m_1)&=1-m_2-(l-1)m_1\\
\nonumber s_{((l-1,1),-)}^1(m_2,m_1)&=2-lm_1
\end{align}
Under the condition (\ref{eq:condm}) (i.e. if
$1\leq m_1\leq m_2$) then we see that the first one
is indeed always smaller than the second one.

\emph{Type $B_l(l\geq3)$ ($\eta=\omega_l/2$):} Not minuscule, with
$W^a_\eta=W(D_l)$ and $W_\eta=S_l$, so this reduces to the
minuscule case $D_l, \eta=\omega_l$ (with $m=m_1$) if $l\geq4$, or
to $A_3, \eta=\omega_1$ if $l=3$.

\emph{Types $BC_l(l\geq 1)$ and $C_l(l\geq2)$:}
We treat these cases together, since the geometry of
$\Omega$ is the same.

We have one boundary orbit to consider, namely $\eta=\omega_l$,
a minuscule case. We have $W^a_\eta=W(C_l)$
and $W_\eta=W(A_{l-1})=S_l$, with root multiplicities $m^a_{a_0}=m_1$
for the long roots of $C_l$, and $m_2$ for the short roots of $C_l$.

In the construction of \cite[Proposition 11.4.2]{Ca} it is easy to see
that the irreducible character $\chi_{(i,l-i)}$ of $W(C_l)$ is realized
on the space of polynomials in $\C[x_1,\dots,x_n]$ by the action
of $W(C_l)$ on the monomial $x_1\dots x_{l-i}$ ($i=0,\dots,l$).
Hence this character contains the trivial character of $S_l$
and has dimension $\operatorname{binomial}(l,i)$, and
has its birthday in degree $l-i$.  By dimension
count we find that
\begin{equation}\label{eq:cl}
\operatorname{Ind}_{S_l}^{W(C_l)}(\chi_l)=\oplus_{i=0}^l\chi_{(i,l-i)}
\end{equation}
and so
\begin{equation}
J_{\omega_l}=\{(\chi_{(i,l-i)}\}_{i=0}^l
\end{equation}
We also see easily from the above realization that
\begin{equation}\label{eq:clno}
s_{(i,l-i)}^1(m)=(l-i)(1-m_1-im_2)
\end{equation}
The characters $\e\otimes\chi_{(i,l-i)}=
\chi_{(1^{l-i},1^i)}$ all contain the sign representation of $S_l$
($i=0,\dots,l$); thus together they fill up (multiplicity free) the character
of $W(C_l)$ induced from the sign representation of $S_l$.
According to \cite[Proposition 11.4.2]{Ca} the birthday
of $\chi_{(1^{l-i},1^i)}$ is in degree
$|\Sigma(D_{l-1})_+|+|\Sigma(C_i)_+|=l(l-1)+i((2l-1)-2i)$.
We see that the minimum is attained for $i=r$ if $l=2r$ or if
$l=2r+1$. Hence if $l=2r$ we get
\begin{equation}
\s_{\omega_{2r}}=\chi_{(r,r)}
\end{equation}
whereas in the case $l=2r+1$ we have
\begin{equation}
\s_{\omega_{2r+1}}=\chi_{(r,r+1)}
\end{equation}

Case $l=2r(r\geq1)$ even: One checks that
\begin{equation}
s_{(i,2r-i)}^1(m)-s_{(r,r)}^1(m)=(r-i)((r-i)m_2-m_1+1)
\end{equation}
which is strictly positive on $\mathcal{C}$ for $0\leq i\leq 2r$
and $i\not=r$.

Case $l=2r+1(r\geq0)$ odd: One checks that
\begin{equation}
s_{(i,2r-i+1)}^1(m)-s_{(r,r+1)}^1(m)=(r-i)((r-i+1)m_2-m_1+1)
\end{equation}
For $0\leq i\leq 2r+1$ and $i\not=r$ this is nonnegative on $\mathcal{C}$,
and it is zero precisely when $m_1=1$ and $i=r+1$.
Observe that this is also true if $r=0$.

\emph{Type $D_l(l\geq4)$, $\eta=\omega_1$:} This is a minuscule case.
Recall that the irreducible characters $\chi_{(\l,\mu)}$ of
$D_l$ are parameterized by unordered pairs $(\lambda,\mu)$
of partitions of total weight $l$ where $\l\not=\mu$, and
characters $\chi_{(\l,\l)}^\prime,\chi_{(\l,\l)}^{\prime\prime}$
if $l$ is even (weight of $\l$ is $l/2$).
The character $\chi_{(\l,\mu)}$ is the restriction of the
character of $W(B_l)$ with the same label $(\l,\mu)$ to $W(D_l)$.
This restriction stays irreducible unless $\l=\mu$, in which case
the character splits as a sum of two irreducible characters
which we distinguish by ${}^\prime$ and ${}^{\prime\prime}$. Thus
$\chi_{(l,0)}=1$ and $\chi_{(1^l,0)}=\e$.

By restriction of (\ref{eq:bn}) to $W(D_l)$ we find
\begin{equation}
J_{\omega_1}=\{\chi_{(l,0)},\chi_{(l-1,1)},\chi_{((l-1,1),-)}\}
\end{equation}
{}From \cite[Proposition 11.4.2]{Ca} we find
\begin{equation}
\s_{\omega_1}:=
\e_l\otimes j_{D_{l-1}}^{D_l}(\chi_{\e_{l-1}})=
\e_l\otimes\chi_{((2,1^{l-1}),-)}=\chi_{((l-1,1),-)}
\end{equation}
where the birthday of $\chi_{((2,1^{l-1}),-)}$ is
in $|\Sigma(D_{l-1})_+|=(l-1)(l-2)$.

Using Theorem \ref{prop:ev} and standard results on
representations of $W(D_l)$ (e.g. \cite[Chapter 11]{Ca}) we
find (one should compare this to (\ref{eq:bnj}))
\begin{align*}
s_{(l-1,1)}^1(m)&=1-(l-1)m\\
s_{((l-1,1),-)}^1(m)&=2-lm
\end{align*}
Under the condition (\ref{eq:condm}) (i.e. if
$1\leq m$) then we see indeed that the second one
is smaller than the first one, except in the case
$m=1$ when they coincide.

\emph{Type $D_l(l\geq4)$, $\eta=\omega_l$:} This is minuscule too.
For the computation of $J_{\omega_l}$ we recall the realizations for
the characters $\chi_{(l-i,i)}$ as described in the text above
(\ref{eq:cl}). We introduce an intertwining operator $\J$ for
the restriction of these representations to $W(D_l)$.
If $\Omega\subset\{1,\dots,l\}$ we denote by $x_\Omega$ the product
of the $x_i$ with $i\in\Omega$.
We now define $\J(x_{\Omega})=x_{\Omega^c}$ and extend by
linearity. Then $\J$ is an intertwining isomorphism
$\J:\pi_{(\a,\b)}|_{W(D_l)}\to\pi_{(\b,\a)}|_{W(D_l)}$,
and if $\a=\b$ then $\J$ splits $\pi_{(\a,\a)}|_{W(D_l)}$ in
$\pi_{(\a,\a)}^\prime$ (the $+1$-eigenspace of $\J$)
and $\pi_{(\a,\a)}^{\prime\prime}$ (the $-1$-eigenspace of $\J$).
Thus $\pi_{(\a,\a)}^\prime$ contains the $S_l$-spherical vector
with this convention. Hence if $l=2r$ then
\begin{equation}
J_{\omega_l}=\{\chi_{(r,r)}^{\prime}\}\cup\{\chi_{(i,2r-i)}\}_{i=r+1}^{2r}
\end{equation}
and if $l=2r+1$ then
\begin{equation}
J_{\omega_l}=\{\chi_{(i,2r-i+1)}\}_{i=r+1}^{2r+1}
\end{equation}
As in the text below (\ref{eq:cl}) we find that
if $l=2r+1$ then
\begin{equation}
\s_{\omega_{2r+1}}=\chi_{(r+1,r)}
\end{equation}
whereas if $l=2r$ then
\begin{equation}\label{eq:Dodd}
\s_{\omega_{2r}}=\chi_{(r,r)}^\prime
\end{equation}

In the odd case $l=2r+1$ we thus get the specialization of the result
(\ref{eq:cl}) for $C_l$ at $m_1=0$, namely:
\begin{equation}
s_{\chi_{(2r+1-i,i)}}^1(m)=(2r+1-i)(1-im)
\end{equation}
but this time this has a unique minimal value among
$J_{\omega_l}$ at $i=r+1$ (which proves our claim in this
case, in view of (\ref{eq:Dodd})).
Hence $s_{\omega_{2r+1}}=r(1-(r+1)m)$.

In the even case we need to look more closely at our model for
$\chi_{(r,r)}^\prime$ first. The degree of this representation is
$\operatorname{binomial}(2r,r)/2=\operatorname{binomial}(2r-1,r-1)$.
The dimension of the $-1$ eigenspace of a reflection is equal to
$\operatorname{binomial}(2r-2,r-1)/2=\operatorname{binomial}(2r-3,r-2)$.
This leads to
\begin{equation}
s_{\chi_{(r,r)}}^1(m)=r(1-rm)
\end{equation}
which is still same the same answer as we had in $C_{2r}$
when substituting $m_1=0$ (cf. (\ref{eq:clno})). Therefore this
exponent indeed represents the unique minimal exponent among
those associated with the characters in $J_{\omega_1}$
proving the claim in this case as well.

\emph{Type $E_6$, $\eta=\omega_1$:}
This is minuscule. By the character tables in ``CHEVIE'' we
find that
\begin{equation}
J_{\omega_1}=\{\chi_{1,0},\chi_{6,1},\chi_{20,2}\}
\end{equation}
and
\begin{align*}
s^1_{6,1}(m)&=1-6m\\
s^1_{20,2}(m)&=2-9m
\end{align*}
The second one is the unique minimal exponent,
and we check that
$s^1_{20,2}(2)=-16=|\Sigma(D_5)_+|-|\Sigma(E_6)_+|$.
In view of Proposition \ref{cpx} this proves the claims
in this case.

\emph{Type $E_7$, $\eta=\omega_7$:}
This is minuscule. By the character tables in ``CHEVIE'' we
find that
\begin{equation}
J_{\omega_1}=\{\chi_{1,0},\chi_{7,1},\chi_{27,2},\chi_{21,3}\}
\end{equation}
and
\begin{align*}
s^1_{7,1}(m)&=1-9m\\
s^1_{27,2}(m)&=2-14m\\
s^1_{21,3}(m)&=3-15m
\end{align*}
The last one is the unique minimal exponent,
except when $m=1$ when it coincides with the second one.
We check that
$s^1_{21,3}(2)=-27=|\Sigma(E_6)_+|-|\Sigma(E_7)_+|$.
In view of Proposition \ref{cpx} this proves the claims.

\emph{Type $E_7$, $\eta=\omega_2/2$:}
This is not minuscule, and reduces to the case ($A_7$,
$\eta=\omega_1$).

\emph{Type $E_8$, $\eta=\omega_1/2$:}
This is not minuscule, and reduces to the case
($D_8$, $\omega_1$).

\emph{Type $E_8$, $\eta=\omega_2/3$:}
This is not minuscule, and reduces to the case
($A_8$, $\omega_1$).

\emph{Type $F_4$:}
This is not minuscule, and reduces to the case
($B_4$, $\omega_1$).

\emph{Type $G_2$:}
This is not minuscule, and reduces to the case
($A_2$, $\omega_1$).
\end{proof}
\begin{cor}(of the proof of the previous Theorem)
For $m\in\partial{\mathcal{C}}$ (the boundary
of $\mathcal{C}$ there are at most two inequivalent irreducibles
$\tau,\pi\in J^\eta$ such that $s^1_{\eta,\tau}(m)=
s^1_{\eta,\pi}(m)$. The cases where this occurs
are indicated in the last column of the table of Theorem \ref{thm:table}.
\end{cor}
\begin{cor} Let $m\in\partial{\mathcal{C}}$ be such that
there are two inequivalent irreducible representations
$\tau,\pi\in J^\eta$ with coinciding exponents
$s^1_{\eta,\tau}(m)=s^1_{\eta,\pi}(m)$. Then the term
of (\ref{eq:zsingexp}) corresponding to the leading exponent
$s^1_{\eta,\tau}(m)$ contains possibly a $\log(\e)$ term
of degree at most one. Otherwise the leading
term in (\ref{eq:zsingexp}) has no logarithmic term.
\end{cor}
\begin{proof}
This is an easy consequence of Corollary \ref{cor:nolog} and
the fact that $\phi_{\mu,m}$ is holomorphic in $m$.
Indeed, suppose that an expression of the form
(with $s\in\mathbb{R}\backslash\{0\}$ fixed)
\begin{equation}
f(m,\e)=a(m,\e)+b(m,\e)\e^{sm}
\end{equation}
is a local Nilsson class function of
$(m,\e)\in\mathbb{D}\times\mathbb{D}^\times$
with $a(m,\e)$, $b(m,\e)$ both holomorphic for
$\e\in\mathbb{D}$ for all fixed $m\in\mathbb{D}^\times$.
Then analytic continuation around $\e=0$ implies that
\begin{equation}
f^\prime(m,\e)=a(m,\e)+\exp(2i\pi sm)b(m,\e)\e^{sm}
\end{equation}
is in the local Nilsson class on
$\mathbb{D}\times\mathbb{D}^\times$
too. Hence $a$ and $b$ have poles in $m$ of order at most $1$
and their residues at $m=0$ cancel.
Now use $\log(\e)=s^{-1}\lim_{m\to 0}(m^{-1}(\e^{sm}-1))$.
\end{proof}
This finishes the proofs of Theorem \ref{thm:esteta} and Theorem
\ref{thm:table}.

\end{document}